\documentclass[12pt]{amsart}
\usepackage{amsmath}
\usepackage{amsfonts}
\usepackage{amssymb}
\usepackage{amsthm}
\usepackage[all]{xy}

\newtheorem{theorem}{Theorem}[section]
\newtheorem{lemma}[theorem]{Lemma}
\newtheorem{proposition}[theorem]{Proposition}
\newtheorem{corollary}[theorem]{Corollary}
\newtheorem{definition}[theorem]{Definition}
\newtheorem{remark}[theorem]{Remark}

\numberwithin{equation}{section}

\newcommand{\E}{\mathcal{E}}

\begin{document}

\baselineskip=15pt

\title[Hitchin maps on curves]
{Hitchin maps and parabolic Hitchin maps on the moduli spaces of Hitchin sheaves on nodal curves}
\author[Usha]{Usha N. Bhosle}

\address{ Indian Statistical Institute , Bangalore 560059, India}

\email{usnabh07@gmail.com}

\subjclass[2010]{Primary 14H60; Secondary 14D20}

\keywords{nodal curves , Hitchin pairs, stability, Hitchin maps}

\thanks{ This work was done during my tenure as an INSA Senior Scientist  at the Indian Statistical Institute, Bangalore.  The online program - Quantum Fields, Geometry and Representation Theory 2021 (code: ICTS/qftgrt2021/7) of the International Centre for Theoretical Sciences (ICTS) was useful for this work. I would like to thank Anna Pe\'on-Nieto for very useful correspondence after this program and for reference to the work of Su, Wang and Wen. I thank the referee for the reference \cite{Pe} and suggestions to improve the exposition of the paper.} 

\date{14 Feb 2024}
\maketitle

\begin{abstract}
             
             We study the Hitchin maps on the moduli spaces of Hitchin sheaves and parabolic Hitchin sheaves on a nodal integral curve $Y$. We study their fibres, the BNR correspondences and the relation of the restriction of the Hitchin map with very stable bundles. As an application, in the "coprime" case we prove that the parabolic Poincar\'e sheaf on $U_{par}(n, \xi)\times Y$ is parabolic stable where $U_{par}(n, \xi)$ is a fixed determinant moduli space of parabolic sheaves on $Y$.             
\end{abstract}

\section{Introduction} 

          The theory of Hitchin bundles or Higgs bundles on smooth curves has flourished by leaps and bounds in the last few decades. We initiated the study of $L$-valued Hitchin sheaves on an integral (nodal) curve $Y$ ($L$ being a line bundle on $Y$) in  \cite{Bh1} and extended it to reducible curves in \cite{Bh2}.  In this paper, we continue the study of $L$-valued Hitchin sheaves on an integral nodal curve $Y$, particularly the Hitchin maps, and generalise it to parabolic Hitchin sheaves. 

          We start with a generalisation of some results of Laumon for $Y$ smooth.  Let $Y$ denote an irreducible nodal projective curve defined over $\mathbb{C}$. Let $Vect_{Y,n}$ denote the algebraic stack over $\mathbb{C}$ of vector bundles $E$ of rank $n$ on $Y$. It is a smooth algebraic stack of pure dimension $n^2(g-1)$.  The cotangent bundle of $Vect_{Y,n}$ is denoted by $T^* Vect_{Y,n}$. For an algebraic stack $S$ over $\mathbb{C}$ and a nilpotent Hitchin bundle $(E,\phi), E$ a vector bundle of rank $n$ on $S$,  there is an exact sequence of $\mathcal{O}_{Y\times_{\mathbb{C}} S}$- modules $0 \rightarrow K_i \rightarrow {E} \stackrel{{\phi}^i}{\rightarrow}  {E}\otimes_{\mathcal{O}_Y} (\omega_Y)^{\otimes i}  \rightarrow C_i \rightarrow 0$  and  a filtration $0 = K_0 \subset K_1 \subset \cdots \subset E$ with $rk(K_i/K_{i-1})_{(s)}) = \nu_i, deg(K_i/K_{i-1}) = \lambda_i$, for $s \in S$. There is a unique stratification of $T^* Vect_{Y,n}$ by locally closed substacks  $(T^* Vect_{Y,n})_ { (\nu_{\bullet}, \lambda_{\bullet}) }$ indexed by $(\nu_{\bullet}, \lambda_{\bullet}) = ( (\nu_1, \lambda_1), \cdots, (\nu_m, \lambda_m) )\, .$  Let $\Lambda'_{Y, (\nu_{\bullet}, \lambda_{\bullet})} \subset \Lambda_{Y, (\nu_{\bullet}, \lambda_{\bullet})}$ be the open algebraic substack consisting of $(E,\phi)$ such that $K_i, E/K_i$ and $K_i/ K_{i-1}$ are vector bundles. 
 
\begin{theorem}  
       For every nilpotent type $(\nu_{\bullet}, \lambda_{\bullet})$, the stack $\Lambda'_{Y, (\nu_{\bullet}, \lambda_{\bullet})} \subset T^* Vect_{Y,n}$ is either empty or smooth over $\mathbb{C}$, of pure dimension $n^2(g-1)$ and Lagrangian in $T^* Vect_{Y,n}$.   
\end{theorem}
Let $\Lambda_{Y,n} \ := \ \bigcup_{(\nu'_{\bullet}, \lambda'_{\bullet}), (1^n) \le \nu'_{\bullet}} \ (T^* Vect_{Y,n})_{(\nu'_{\bullet}, \lambda'_{\bullet})}$
where  $(1^n) = (\nu_1, \cdots, \nu_n)$ with $\nu_i=1,  \forall i=1, \cdots, n$ and the union is taken over all tuples $(\nu'_{\bullet}, \lambda'_{\bullet})$ such that $(1^n) \le \nu'_{\bullet}$.

\begin{corollary}
       The locally closed substack $\Lambda'_{Y, n} \subset T^* Vect_{Y,n}$ is reduced and of pure dimension $n^2(g-1)$ and Lagrangian in $T^* Vect_{Y,n}$.  
\end{corollary}
          
          Let $p : X \to Y$ be the normalisation and $L_0 = p^*L$. Define $A' :=  \bigoplus_{i= 1}^n \ H^0(Y, L^i), A := \bigoplus_{i= 1}^n \ H^0(X, L_0^i)$ and ${\bf A}, {\bf A'}$ the Affine spaces with underlying vector spaces $A, A'$ respectively, ${\bf A} \supset  {\bf A'}\, .$  We had defined a proper Hitchin morphism $h: N_Y(n,d,L) \to {\bf A}$ where $N_Y(n,d,L)$ is a subset of the moduli space $\mathcal{H}_Y(n,d,L)$ of $L$-valued Hitchin sheaves of rank $n$, degree $d$ on $Y$. The open subscheme $\mathcal{H}'_Y(n,d,L)$ of all Hitchin bundles is contained in $N_Y(n,d,L)$ and the morphism $h$ maps it into ${\bf A'}$. We study the  map $h$ in more detail. If $E$ is a stable torsion-free sheaf of rank $n$ and degree $d$, then a Hitchin sheaf $(E, \phi)$ is stable for any Hitchin field $\phi$. Hence $V_E := {\bf H}^0(Y, End E \otimes L)$ embeds in  $\mathcal{H}_Y(n,d,L)$, we identify it with its image. We study the restriction $h_V$ of the Hitchin map to $V_{E,1}:= V_E \cap N_Y(n,d,L)$. We check that $V_{E,1}$  is a closed subset of $N_Y(n,d,L)$ if and only if $h_V$ is a proper map and show that if $h_V$ is proper and $V_{E,1}$ is a vector space then $h_V$ is quasi-finite.  
          
      We define $L$-very stable torsion-free sheaves on an integral nodal curve $Y$, show the existence of very stable vector bundles and determine conditions for the existence of very stable non locally free torsion-free sheaves depending on the local type  (see Propositions \ref{u1} and \ref{u2}).  
      
\begin{theorem} (Theorem \ref{Thm1.1PaPe})
    Denote by $E$ a stable vector bundle of rank $n$ and degree $d$ on $Y$. Then the following statements are equivalent: \\
(1)  $E$ is $L$-very stable.\\
(2)  $V_{E,1}$  is a closed subset of $N_Y(n,d,L)$. \\ 
(3)  $h_V$ is a proper morphism.\\
(4)  $h_V$ is quasi-finite.    
\end{theorem}

           We shift our focus to $L$-twisted parabolic (and strongly parabolic) Hitchin sheaves on an integral (nodal) curve $Y$, with parabolic structure over a divisor $D$ contained in the set of nonsingular points of $Y$. We define parabolic GPH (generalised parabolic Hitchin sheaves with parabolic structure) on  $X$ and study their relation to parabolic Hitchin sheaves on $Y$ (Theorem \ref{mgphcorrespondence}).   
 We use it to define a parabolic Hitchin map. We also define strongly parabolic Hitchin bases ${\bf A}_p \subset {\bf A}$ and ${\bf A}'_p \subset {\bf A}'$, the restriction $h^{str}_p$ (respectively $h^{' str}_p$) of $h_p$ to the moduli space of strongly parabolic Hitchin sheaves (respectively Hitchin bundles) maps to ${\bf A}_p$ (respectively ${\bf A'}_p$).  We prove the parabolic BNR correspondences determining the fibres of the Hitchin map for a general $a \in {\bf A'}_p$ and for a general $a \in {\bf A}_p$. We define the spectral curves $\pi_p : Y^p_a \to Y$ in both cases. 

\begin{theorem}  (Theorem \ref{parBNRcorrA'})

Let $a \in {\bf A}'_p$ be a general element and $\delta = d + \frac{(n^2 - n)}{2}d(L) + \sum_{x\in D} \ {\rm dim} \ (GL(n)/P_x)$. Then we have:\\
(1) The direct image map $(\pi_p)_*$ induces an injective morphism 
$$h_T: \bar{J}^{\delta}(Y^p_a) \longrightarrow \mathcal{H}^{str}_{par}(n,d,L)\, .$$  
(2) It induces an isomorphism $J(Y^p_a) \cong (h_p^{' str})^{-1} (a)$. \\
(3) The parabolic Hitchin fibre $(h^{str}_p)^{-1} (a) = h_T(\bar{J}^{\delta}(Y)) \cap N^{str}_{par}(n,d,L)$. 
\end{theorem}

\begin{theorem} (Theorem \ref{paraBNRcorrA}) Let $m'$ be the number of nodes of $Y$.

Let $a \in {\bf A}_p$ be a general element and $\delta_A = d+ \frac{(n^2 - n)}{2} d(p_*L_0) + \sum_{x\in D} \ dim \ (GL(n) /P_x)\, ,$ where $d(L_0) = d(L)+ m'$. Then we have:\\
(1) The direct image map $(\pi_p)_*$ induces an injective morphism 
$$h_T: \bar{J}^{\delta}(Y^p_a) \longrightarrow \mathcal{H}^{str}_{par}(n,d,p_*L_0)\, .$$  
(2) It induces an isomorphism $J(Y^p_a) \cong (h_p^{' str})^{-1} (a)$. \\
(3) The parabolic Hitchin fibre $(h^{str}_p)^{-1} (a) = h_T(\bar{J}^{\delta}(Y)) \cap N^{str}_{par}(n,d,L)$. 

\end{theorem}
 
                Finally we assume that $L = \omega_Y$. We study the parabolic nilpotent cone and the strongly parabolic nilpotent cone and determine their dimensions (Theorems \ref{P6.7SWW} and \ref{T6.9SWW}). We define very stable parabolic sheaves. We prove their existence (Theorem \ref{verystableexist}) and prove an analogue of Theorem \ref{Thm1.1PaPe} in the parabolic case. We study the forgetful rational map from the moduli space of parabolic Hitchin sheaves of rank $n$ and degree $d$ to the moduli space $U_{par}(n,d)$ of parabolic torsion-free sheaves of the same rank and degree. 
                
                We define very stable (strongly) parabolic Hitchin sheaves.  For a strongly quasi-parabolic sheaf $E$, define  $V^{str}_E := {\bf H}^0(Y, SParEnd(E)\otimes \omega_Y(D))$ where $SParEnd(E)$ denotes the sheaf of strongly parabolic endomorphisms of $E$. 
\begin{theorem}
Let $E_*$ be a stable parabolic bundle. Let 
$$h^{str}_{E} = h^{' str}_p\vert_{V^{str}_E}: V^{str}_E \to {\bf A'}_p$$
 be the restriction of the strongly parabolic Hitchin morphism $h^{' str}_p$. Then the following are equivalent:\\
(1) E is strongly very stable.\\
(2) $h^{str}_E$ is finite.\\
(3) $h^{str}_E$ is quasi-finite.\\
(4) $h^{str}_E$ is proper. \\
(5) The inclusion map $i^{str}_V: V^{str}_E \rightarrow N^{str}(n,d,\omega_Y)$ is proper. 
\end{theorem}

              We now specialise to the fixed determinant Hitchin moduli spaces.  The fixed determinant moduli space ${\mathcal H}^{str}_{par}(n, \xi, L) \subset  {\mathcal H}^{str}_{par}(n, d, L)$ is defined as the closure of ${\mathcal H}^{' str}_{par}(n, \xi, L)$, the moduli space of strongly parabolic Hitchin bundles of rank $n$ and a fixed determinant $\xi$ with the Hitchin field $\phi$ having values in $L(D)$. For $(E_*, \phi) \in {\mathcal H}^{str}_{par}(n, \xi, L)$, one has $tr \phi = 0$. This gives the parabolic Hitchin bases  ${\bf A}_{p, \xi}$ and  ${\bf A'}_{p, \xi}$. We determine the nilpotent cones in this case. 
              
            Let $n_i(x), x \in D$ denote the multiplicities in the parabolic structure.  Assume that the greatest common divisor of $n, d$ and $n_i(x), x \in D$, is $1$.  Then there is a universal family $\mathcal{U}_*$, of parabolic torsion-free sheaves, on  $U_{par}(n, \xi) \times Y$.  
            
            Let $z$ be a closed point of  $Y$ distinct from the nodes and the points in $D$. Let 
$\mathcal{U}_z$ be the  restriction of the universal sheaf to $U_{par}(n, \xi) \times z$. As $z$ is a nonsingular point, $\mathcal{U}_z$ is a vector bundle, note that $\mathcal{U}_z$ has no parabolic structure. 

\begin{theorem}  (Theorem \ref{Uzsemistable})
The vector bundle  $\mathcal{U}_z$ is semistable with respect to the ample theta line bundle $\Theta_{par}$ on $U_{par}(n, \xi)$.
\end{theorem}                        

\begin{theorem} (Theorem  \ref{Uparstable}) 
     The parabolic torsion-free sheaf  $\mathcal{U}_*$ on $U_{par}(n, \xi) \times Y$ is parabolic stable with respect to any integral ample divisor of the form $a \Theta_{par} + b D_Y$ where $a, b >0$ and $D_Y$ is an ample divisor on $Y$. 
\end{theorem}          

                The last two theorems were proved in \cite{BaBiD} under the assumptions that the curve $Y$ is smooth and the parabolic structure has $n_i(x) = 1$ for all $x, i$. We have removed both the assumptions by generalising it to the nodal curve $Y$ and proving it for any type of parabolic structure.

\section{The nilpotent cone for $Y$}

        Let $Y$ be an irreducible nodal curve defined over complex numbers $\mathbb{C}$.  In this section, we study (an open subset of) the nilpotent cone for $Y$ following the work of Laumon \cite{La} in case $Y$ is smooth. We need several modifications.
      
\subsection{The algebraic stack $Vect_{Y,n}$} \hfill

               Let $Vect_{Y,n}$ denote the algebraic stack over $\mathbb{C}$ of vector bundles $E$ of rank $n$ on $Y$. It is a smooth algebraic stack of pure dimension $n^2(g-1)$. We note that the subset of the appropriate quote scheme corresponding to vector bundles $E$ of rank $n$ and degree $d$ is irreducible and smooth \cite[Remark, p. 167]{N}. For each $d \in \mathbb{Z}$, let $Vect^{d}_{Y,n}$ denote the algebraic stack over $\mathbb{C}$ of vector bundles $E$ of rank $n$ and degree $d$, the $Vect^{d}_{Y,n}$'s are the connected components of $Vect_{Y,n}$. Let
               
               $$\tilde{E} \to Y \times _{\mathbb{C}} Vect_{Y,n}$$
 be the universal vector bundle.  The cotangent complex of  $Vect_{Y,n}$ is given by 
 $$L_{Vect_{Y,n}/\mathbb{C}} \ = \ Rp_{2 *} \underline{Hom}(\tilde{E}, \tilde{E}\otimes_{\mathcal{O}_Y} \omega_Y)$$
 \cite[II, 1.2.7]{Il} where $p_1$ and $p_2$ are the projections from  $Y \times _{\mathbb{C}} Vect_{Y,n}$  to $Y$ and  $Vect_{Y,n}$ respectively. We note that the dualising sheaf $\omega_Y$ is a locally free sheaf of rank $1$.  By Grothendieck duality, the tangent sheaf of $Vect_{Y,n}$ is 
 $$T_{Vect_{Y,n}/ \mathbb{C}} \ = \ R^1 p_{2 *} \underline{Hom}(\tilde{E}, \tilde{E})\, .$$
 The cotangent bundle of $Vect_{Y,n}$ is 
 $$ \tau: T^* Vect_{Y,n} \longrightarrow Vect_{Y,n}\ .$$
For every scheme $S$ over $\mathbb{C}$,  $T^* Vect_{Y,n}/ \mathbb{C}$ is the category of couples $(E, u)$ where $E$ is a vector bundle of rank $n$ on $Y \times_{\mathbb{C}}  S$ and $u: E \to E \otimes \omega_Y$ is a homomorphism of $\mathcal{O}_{Y\times_{\mathbb{C}} S}$- modules (as in \cite[Lemma 1.1]{La} ). Thus $T^* Vect_{Y,n}$ is an algebraic stack over $\mathbb{C}$ of couples $(E,u)$ and there is a universal homomorphism of $\mathcal{O}_{Y \times_{\mathbb{C}} T^* Vect_{Y,n} }$-modules 
\begin{equation} \label{La1.2}
\tilde{u}\colon (Id \times \tau)^* \tilde{E} \ \longrightarrow (Id \times \tau)^* \tilde{E} \otimes_{\mathcal{O}_Y} \omega_Y\, .
\end{equation}  
 
          For an algebraic stack $S$ over $\mathbb{C}$ and a couple $(E,u)$ where $E$ is a vector bundle of rank $n$ on $S$ and $u: E \to E \otimes_{\mathcal{O}_Y} \omega_Y$ a homomorphism, define $u^i: E \to E \otimes_{\mathcal{O}_Y} (\omega_Y)^{\otimes i}, i\ge 1$ as the composite $u^i = (u\otimes (Id_{\omega_Y})^{\otimes i-1})\circ \cdots \circ (u\otimes Id_{\omega_Y})\circ u$. 
          
\begin{definition}\label{La0.2}
     The couple $(E,u)$ is called nilpotent if $u^n = 0$.
\end{definition}
          
        For a nilpotent couple $(E,u)$, there is an exact sequence of $\mathcal{O}_{Y\times_{\mathbb{C}} S}$- modules
\begin{equation} \label{La1.3}
0 \rightarrow K_i \rightarrow {E} \stackrel{{u}^i}{\rightarrow}  {E}\otimes_{\mathcal{O}_Y} (\omega_Y)^{\otimes i}  \rightarrow C_i \rightarrow 0\, ,
\end{equation} 
and  a filtration 
\begin{equation} \label{La1.4}
0 = K_0 \subset K_1 \subset \cdots \subset E
\end{equation} 
with  monomorphisms 
\begin{equation} \label{La1.5}
\frac{K_{i+1}}{K_i} \ \stackrel{\bar{u}}{\longrightarrow} \ \frac{K_i}{K_{i-1}} \otimes_{\mathcal{O}_Y} \omega_Y 
\end{equation}  
induced by $u$.            
             If the $\mathcal{O}_{Y\times_{\mathbb{C}} S}$-modules $C_i, 1\le i \le n+1$ are flat over $S$, then $K_i, E/K_i$ and $K_i/ K_{i-1}$ are torsion free  on $Y \times_{\mathbb{C}} s$ for every point $s \in S$ (and are flat over $S$) for $1 \le i \le n$. For $i >n$, $K_i/ K_{i-1}=0$ and $K_i =K_n$.

        Following \cite[Definition 1.7]{La}, we make a definition.     
\begin{definition} \label{La1.7}
Let   $$(\nu_{\bullet}, \lambda_{\bullet}) = ( (\nu_1, \lambda_1), \cdots, (\nu_m, \lambda_m) )\, .$$
 We say that a couple $(E,u)$ is of nilpotent type $(\nu_{\bullet}, \lambda_{\bullet})$  for $Y$ if the following conditions hold. \\
(a) the coherent $\mathcal{O}_{Y\times_{\mathbb{C}} S}$-modules $C_i, 1\le i \le n+1$ are flat over $\mathcal{O}_S$,\\
(b) for all points $s\in S$ and each $i= 1, \cdots, m$, $rk(K_i/K_{i-1})_{(s)}) = \nu_i, deg(K_i/K_{i-1}) = \lambda_i$,\\
(c) for each integer $i > m$, we have $K_i = K_m$.  
\end{definition}             
    
    If $S= Spec \ k, \ k\supset \mathbb{C}$ a field, then every nilpotent couple $(E,u)$ admits a nilpotent type  for $Y$ of length $m \le n$.

    Define an order relation 
$$ \nu_{\bullet} \le \nu'_{\bullet}  \ {\rm if \ and \ only \ if} \ \nu_1+\cdots+\nu_i \le \nu'_1+\cdots+\nu'_i \ \forall i\ge 1,$$
where $\nu_i=0$ if $i>m$ and $\nu'_i=0$ if $i >m'$.

One can show that there is a unique stratification of $T^* Vect_{Y,n}$ by locally closed substacks  $(T^* Vect_{Y,n})_ { (\nu_{\bullet}, \lambda_{\bullet}) }$ indexed by the types (using \cite[Lecture 8]{Mu} as in \cite[Theorem 1.9]{La}).  
The stratification has same properties as in \cite[Remark 1.8]{La}. 
Let
\begin{equation} \label{La1.12}
\Lambda_{Y, (\nu_{\bullet}, \lambda_{\bullet})} \ := \ 
(T^* Vect_{Y,n})_{(\nu_{\bullet}, \lambda_{\bullet})} \, .
\end{equation}
  For every $(\nu_{\bullet}, \lambda_{\bullet})$, the algebraic stack $\Lambda_{Y, (\nu_{\bullet}, \lambda_{\bullet})}$ is the algebraic stack of couples $(E,u)$ of nilpotent type $(\nu_{\bullet}, \lambda_{\bullet})$. 
Define 
\begin{equation} \label{La1.11}
\Lambda_{Y,n} \ := \ \bigcup_{(\nu'_{\bullet}, \lambda'_{\bullet}), (1^n) \le \nu'_{\bullet}} \ 
(T^* Vect_{Y,n})_{(\nu'_{\bullet}, \lambda'_{\bullet})}
\end{equation} 
where  $(1^n) = (\nu_1, \cdots, \nu_n)$ with $\nu_i=1,  \forall i=1, \cdots, n$ and the union is taken over all tuples $(\nu'_{\bullet}, \lambda'_{\bullet})$ such that $(1^n) \le \nu'_{\bullet}$.             
             The  stratification of $\Lambda_{Y,n}$ is called the canonical stratification ( \cite[Rmk1.13]{La}).

 For every integer $d \in \mathbb{Z}$, $\Lambda_{Y,(n,d)}$ is the zero section of $T^* Vect^d_{Y, n}$ (here $(n,d)$ is the type with $m=1, \nu_1=n, \lambda_1=d$ so that $u=0$). 
 Define 
 \begin{equation} \label{lambda'}
      \Lambda'_{Y, (\nu_{\bullet}, \lambda_{\bullet})} \subset \Lambda_{Y, (\nu_{\bullet}, \lambda_{\bullet})},
\end{equation}
    the open algebraic substack consisting of couples $(E,u)$ such that $K_i, E/K_i$ and $K_i/ K_{i-1}$ (defined by equation (\ref{La1.3}) ) are vector bundles  on $Y \times_{\mathbb{C}} S$ (and are flat over $S$) for $1 \le i \le n$.   
      
      For every integer $d \in \mathbb{Z}$, the zero section $\Lambda_{Y,(n,d)}$ of $T^* Vect^d_{Y, n}$ lies in $\Lambda^{' d}_{Y, n}$. The complement of the zero section in $\Lambda^{' d}_{Y, n}$ is nonempty. For example, for $E= \mathcal{O}_Y$, one can define nonzero nilpotent $u : E \to E \otimes \omega_Y$. 

\subsection{The algebraic stacks $Vect^{\lambda_{\bullet}}_{Y, \nu_{\bullet}}$ and $Tf^{\lambda_{\bullet}}_{Y, \nu_{\bullet}}$}   \hfill

          We fix a sequence of positive integers $\nu_{\bullet}= (\nu_1, \cdots, \nu_m)\, , \sum_{i=1}^m \nu_i = n$ and a sequence of integers $\lambda_{\bullet}= (\lambda_1, \cdots, \lambda_m), 
\sum_{i=1}^m \lambda_i = d$.

           The algebraic stack $Vect^{\lambda_{\bullet}}_{Y, \nu_{\bullet}}$  on $\mathbb{C}$ is defined by associating to every scheme $S$ over $\mathbb{C}$ the groupoid  $Vect^{\lambda_{\bullet}}_{Y, \nu_{\bullet}}(S)$. Here $Vect^{\lambda_{\bullet}}_{Y, \nu_{\bullet}}(S)$ is the category whose objects are  flags 
\begin{equation} \label{La2.1} 
E_{\bullet} = ( 0 = E_0 \subset E_1 \subset \cdots \subset E_m = E )
\end{equation}                
where $E_i$ are vector bundles  on 
$Y \times_{\mathbb{C}} S$, $E_i/E_{i-1}$ are flat over $S$
and for every point $s \in S$, the coherent $\mathcal{O}_{Y \times_{\mathbb C} }s $- modules  $(E_i/E_{i-1})_{(s)}$ have generic rank $\nu_i$ and degree $\lambda_i$; the morphisms are isomorphisms of flags. There is a representable morphism of stacks 
\begin{equation} \label{La2.2} 
\pi^{\lambda_{\bullet}}_{Y, \nu_{\bullet}} \colon \ Vect^{\lambda}_{Y, \nu_{\bullet}} \longrightarrow 
Vect^d_{Y,n}
\end{equation}                
which, for $S/ \mathbb{C}$, is the morphism $Vect^{\lambda_{\bullet}}_{Y, \nu_{\bullet}}(S) \rightarrow 
Vect^d_{Y,n}(S)$ defined by $E_{\bullet} \mapsto E$. This is an immediate consequence of the representability of $Quot^P$ functors ( \cite[2.A.1, p.51]{HL} for general case, \cite[Lemma 2.3]{La} for $Y$ a smooth curve). 

\begin{corollary} \label{LaCorollary2.10}
The algebraic stack $Vect^{\lambda_{\bullet}}_{Y, \nu_{\bullet}}$ is smooth on $\mathbb{C}$, of pure dimension $d^{\lambda}_{\nu}$ where 
$$d^{\lambda}_{\nu} \ = \ \sum_{1\le i \le j \le m} \  \nu_i \lambda_j - \nu_j \lambda_i + \nu_i \nu_j (g-1)\, .$$ 
\end{corollary}
\begin{proof}
See \cite[Corollary 2.10]{La}
\end{proof}        

c\begin{theorem} \label{LaTheorem 3.1}
       For every nilpotent type $(\nu_{\bullet}, \lambda_{\bullet})$, $\Lambda'_{Y, (\nu_{\bullet}, \lambda_{\bullet})} \subset T^* Vect_{Y,n}$ is either empty or smooth over $\mathbb{C}$, of pure dimension $n^2(g-1)$ and Lagrangian in $T^* Vect_{Y,n}$.   
\end{theorem}

\begin{proof}
        The proof is on the lines of the proof of \cite[Theorem 3.1]{La}, we briefly sketch it with the modifications needed. As noted after equation \eqref{lambda'}, $\Lambda'_{Y,n}(\mathbb{C})$ is nonempty. Hence there exists a nilpotent type $(\nu_{\bullet}, \lambda_{\bullet})$ such that $\Lambda'_{Y, (\nu_{\bullet}, \lambda_{\bullet})}$ is nonempty. Fix one such nilpotent type. Let $d= \lambda_1+ \cdots + \lambda_m$ so that $\Lambda'_{Y, (\nu_{\bullet}, \lambda_{\bullet})} \subset T^* Vect^d_{Y,n}$. The morphism (defined in \eqref{La2.2})
$$\pi:=  \pi^{\lambda_{\bullet}}_{Y, \nu_{\bullet}} \colon \ Vect^{\lambda_{\bullet}}_{Y, \nu_{\bullet}} \longrightarrow  Vect^d_{Y,n}$$
  induces the cotangent map
 $$\Pi \colon T^* Vect^d_{Y,n} \times_{Vect^d_{Y,n}} Vect^{\lambda_{\bullet}}_{Y, \nu_{\bullet}}  \longrightarrow  T^* Vect^{\lambda_{\bullet}}_{Y, \nu_{\bullet}}\, .$$
 Let 
 $$\bar{\pi} \colon T^* Vect^d_{Y,n} \times_{Vect^d_{Y,n}} Vect^{\lambda_{\bullet}}_{Y, \nu_{\bullet}}  \longrightarrow  T^* Vect^d_{Y,n}$$
 be the (first) projection. 
 
           Let $\Sigma$ be the inverse image by $\Pi$ of  the null section of $T^* Vect^{\lambda_{\bullet}}_{Y, \nu_{\bullet}}, \Sigma \subset T^* Vect^d_{Y,n} \times_{Vect^d_{Y,n}} Vect^{\lambda_{\bullet}}_{Y, \nu_{\bullet}} $. For every scheme $S$ on $\mathbb{C}$, $\Sigma(S)$ is the category of couples $(E_{\bullet}, u: E \to E \otimes_{\mathcal{O}_Y}  \omega_Y)$ where $E_{\bullet} = (0 = E_0 \subset E_1 \cdots \subset E_m= E)$ is an object of $Vect^{\lambda_{\bullet}}_{Y, \nu_{\bullet}} (S), u(E_i) \subset E_{i-1}   \otimes_{\mathcal{O}_Y }\omega_Y,  \forall 1\le i \le m$, in particular $E_i = K_i \forall i$ with the notations of \eqref{La1.3}. One has a unique stratification $\Sigma_{(\nu_{\bullet}, \lambda_{\bullet})}$ of $\Sigma$ indexed by nilpotent types. One sees that the immersion 
           $$\Sigma_{(\nu_{\bullet}, \lambda_{\bullet})} \longrightarrow \Sigma$$
is open and the restriction of $\bar{\pi}$ to $\Sigma_{(\nu_{\bullet}, \lambda_{\bullet})}$ is an isomorphism of $\Sigma_{(\nu_{\bullet}, \lambda_{\bullet})}$ on $\Lambda'_{Y, (\nu_{\bullet}, \lambda_{\bullet})}$.  Then the theorem follows from \cite[Proposition 3.3]{La}. 
\end{proof}

The theorem has an immediate corollary.

\begin{corollary}\label{lagrangian}
       The locally closed substack $\Lambda'_{Y, n} \subset T^* Vect_{Y,n}$ is reduced and of pure dimension $n^2(g-1)$ and Lagrangian in $T^* Vect_{Y,n}$.  
\end{corollary}

\section{The moduli spaces and the Hitchin map}

                We recall some definitions and results needed. We fix a line bundle $L$, with $d(L) \ge 0$, on the integral  curve $Y$.  Let  $L_0 = p^*L$, the pull back of $L$ to the normalisation $X$. Let  $\mathcal{H}_Y(n,d,L)$ (respectively $\mathcal{H}_X(n,d,L_0)$) be the moduli space of semistable Hitchin sheaves $(E, \phi)$ on $Y$ (respectively Hitchin bundles on $X$) of rank $n$ and degree $d$ with the Hitchin field $\phi$ having values in $L$  (respectively $L_0$). Let  $\mathcal{H}'_Y(n,d,L) \subset \mathcal{H}_Y(n,d,L)$ be the open subset corresponding to Hitchin bundles. 
                
                 Let $U(n,d)$ (respectively $U'(n,d)$) denote the moduli space of semistable torsion-free sheaves (respectively vector bundles) of rank $n$ and degree $d$ on $Y$.

\subsection{The Hitchin map} \label{hitchinmap} \hfill

Define 
\begin{equation} \label{A}
A :=  \bigoplus_{i= 1}^n \ H^0(X, L_0^i)\, , A' :=  \bigoplus_{i= 1}^n \ H^0(Y, L^i)\, . 
\end{equation}
Note that $H^0(Y, L^i) \subset H^0(Y, L^i\otimes p_*{\mathcal O}_X) = H^0(Y, p_*(L_0^i))\cong H^0(X, L_0^i)$, by the projection formula. Hence  $A' \subset A$. 

For any line bundle $N$, we write ${\bf H}^0(X, N)$ for the affine variety underlying the vector space $H^0(X, N)$. Let 
\begin{equation} \label{bfA}
{\bf A} :=  \bigoplus_{i= 1}^n \ {\bf H}^0(X, L_0^i)\, , \  {\bf A}' :=  \bigoplus_{i= 1}^n \ {\bf H}^0(Y, L^i)\, , {\bf A}' \subset {\bf A}\, .
\end{equation}

\subsection{The set map $h_{set}$} \label{hset}   \hfill

           In this section, for an integral nodal curve $Y$, we define a set map $h_{set}: \mathcal{H}(n,d,L) \to A\, .$
\begin{definition} \label{localtype}

          For any point $y \in Y$, let  $\mathcal{O}_y$ be the local ring and $m_y$ the maximum ideal. Let $E$ be a torsion-free sheaf on $Y$ and let $E_{(y)}$ denote the stalk of  $E$ at $y \in Y$. Since $Y$ is nodal, there is an isomorphism (not unique) 
            $$E_{(y)} \ \cong  a_y \mathcal{O}_y + b_y m_y\ , a_y \mathcal{O}_y = \mathcal{O}_y^{\oplus a_y}, b_y m_y= m_y^{\oplus b_y}, 0 \le a_y, b_y \le n\, .$$
The integer $b_y$ is called the local type of $E$ at $y$. If $y$ is a nonsingular point, then $b_y= 0$. Hence we define the {\bf local type} of $E$ as the $m'$-tuple $b(E)= \{(b_{y_j})_j\}$ as  $y_j$ varies over singular points $y_j, j= 1, \cdots, m'$, of $Y$.  
\end{definition}          
           
           Let $M(p, q, S)$ denote the set of $p \times q$- matrices with entries in $S$.   
           To every Hitchin sheaf $(E, \phi)$ where $E$ is a torsion-free sheaf of rank $n$, degree $d$ and of local type $b(E)$ (see Definition \ref{localtype}), one can associate an element of $A$ as follows. 
            
          Let $Hom(E, E\otimes L^i)_{(y)}$ be the stalk of the sheaf $Hom(E, E\otimes L^i)$ at $y \in Y$. Then
  \begin{equation} \label{iso1}
  Hom(E, E\otimes L^i)_{(y)} \cong Hom_{\mathcal{O}_y} (a_y \mathcal{O}_y + b_y m_y, (a_y \mathcal{O}_y+b_y m_y)\otimes L^i_{(y)})\, .
  \end{equation}
  Hence any $q \in  Hom(E, E\otimes L^i)_{(y)} $ can be written as a matrix $\mathbf{M}$ of the form 
  $$
  \begin{array}{ccccc}       
 \vert & A & B & \vert  & {}\\
 \vert  & {}  & {}  & \vert & \ \in M(a_y+ b_y, a_y+ b_y, \widetilde{\mathcal{O}}_y\otimes L^i_{(y)}) \, ,\\
 \vert & C & D & \vert & {} 
\end{array} 
$$  
    where $A \in  M(a_y, a_y, \mathcal{O}_y\otimes L^i_{(y)})\, , B \in M(a_y, b_y, \widetilde{\mathcal{O}}_y\otimes L^i_{(y)})\, , C \in M(b_y, a_y, m_y\otimes L^i_{(y)})\, , D \in M(b_y, b_y, \widetilde{\mathcal{O}}_y\otimes L^i_{(y)})\ ,$     
(see \cite[Chapter VIII, Proposition 5, p. 170]{Se} for details). Hence for all $y \in Y$, we have the trace of $\mathbf{M}$ defined by 
tr $\mathbf{M} \in \widetilde{\mathcal{O}}_y\otimes L^i_{(y)}$. Define Tr$_i q = $ tr $\mathbf{M}$. Since trace is invariant under conjugation (by an element of GL$(a_y+ b_y, \widetilde{\mathcal{O}}_y\otimes L^i_{(y)})$ ), it follows that Tr$_i q$ is independent of the choice of the isomorphism  \eqref{iso1}. Hence we have the map 
Tr$_i : Hom(E, E) \otimes L^i \to L^i \otimes p_*\mathcal{O}_X$. Taking  $H^0( )$, we have the map 
$tr_i :  H^0(Y, Hom (E, E) \otimes L^i) \to H^0(Y, L^i \otimes p_*\mathcal{O}_X)$. Let   
$tr:= \oplus_{i=1}^{n} tr_i$. Define 
$$h_{set} (E, \phi) := tr(\phi) \in A\, .$$
This gives the set theoretic Hitchin map 
$$h_{set}: \mathcal{H}(n,d,L) \to A\, .$$
We note that if $E$ is a vector bundle, $b_y =0$ for all $y$. Hence  tr $\mathbf{M} \in L^i_{(y)}$ so that 
$h_{set} (\mathcal{H}'(n,d,L)) \subset A' \, .$

          It is not clear if there is a Hitchin morphism $\mathcal{H}(n,d,L) \to {\bf A}$. For any integral curve $Y$, there is certainly a Hitchin  morphism from the open subset $\mathcal{H}'(n,d,L)$ to ${\bf A}'$ (as in the case $Y$ is smooth). For any integral curve $Y$, I used GPH (Generalised Parabolic Hitchin pairs) to get a Hitchin morphism on a subset $N(n,d,L) = N_Y(n,d,L) \subset \mathcal{H}(n,d,L)$ containing $\mathcal{H}'(n,d,L)$ (\cite{Bh1}, \cite{Bh2}).  The Hitchin morphism and the  subset $N(n,d,L)$ will be described in more detail in the following subsection.

\subsection{ Generalised Parabolic Hitchin pairs (GPH) on an integral curve $Y$} \hfill

Let $Y$ be an integral projective curve. Let $D_j, j = 1, \cdots m$ be mutually disjoint divisors on $Y$ such that for any $j$ the support of $D_j$ does not contain any node of $Y$. Let $L'_0$ be a line bundle on $Y$.

Let $0 < \alpha \le 1$ be a real number. 

\begin{definition}\label{gphdefi}    
A generalised parabolic $L'_0$-twisted Hitchin pair of rank $n$ and degree $d$ on $Y$ is a triple $(E, F(E), \phi)$ 
where \\
(1) $E$ is a vector bundle of rank $n$ and degree $d$ on $Y$,\\
(2) $F(E) = (F_1(E), \cdots, F_{m'}(E))$ where $F_j(E) \subset  E \otimes {\mathcal O}_{D_j}$ is a subspace ($E \otimes {\mathcal O}_{D_j}$ being regarded as a $k$-vector space), \\
(3) $\phi: E \rightarrow E \otimes L'_0$ is a homomorphism. \\
We may call a  generalised parabolic Hitchin pair a GPH in short. \\
 Let 
 $$f_j(E) = \ {\rm dim} \ F_j(E)\, , \ Q_j(E)= (E \otimes {\mathcal O}_{D_j})/F_j(E)\, , \ q_j(E)= \ {\rm dim} \ Q_j(E)\, .$$ 
\end{definition}

\begin{definition}  
    A morphism of GPH $(E_1, F(E_1), \phi_1)$ and $(E_2,F(E_2),\phi_2)$ is a homomorphism of vector bundles $f: E_1 \to E_2$ which is compatible with $\phi_1, \phi_2$ and preserves the generalised parabolic structures, i.e. for the induced map 
    $f_{D_j} \colon E_1 \otimes {\mathcal O}_{D_j} \to E_2 \otimes {\mathcal O}_{D_j}$, one has  
$f_{D_j}(F_j(E_1)) \subset F_j(E_2)$.  
\end{definition}

     We recall that the pair $(E, F(E))$ (together with $\alpha$) is called a generalised parabolic bundle, a GPB in short. By forgetting the Hitchin field $\phi$ in a GPH, we get a generalised parabolic bundle (see \cite{B1} for generalities on GPBs). 

For a subbundle $N \subset E$, define $F_j(N) := F_j(E)\cap (N \otimes {\mathcal O}_{D_j})$ and 
$f_j(N)=$ dim $F_j(N)$. 

$N$ is called $\phi$-invariant if $\phi(N) \subset N \otimes L'_0$. 

\begin{definition}  \label{GPHstability}
A  GPH $(E, F(E), \phi)$ is stable (respectively semistable) if for every proper 
$\phi$-invariant subsheaf $N \subset E$ (with a torsion-free quotient) with induced GPB structure, one has 
$$\frac{d(N)+ \alpha \sum_j f_j(N)}{r(N)} \ \ < \ (\le) \ \ \frac{d(E)+ \alpha \sum_j f_j(E)}{r(E)}\, ,$$ 
\end{definition}

\subsection{Good GPH} \label{goodGPH} \hfill
                
 Let  $Y$ be an integral projective curve  and        
       $$p : X  \rightarrow Y$$ 
 the normalisation map. Let ${\mathcal C}= \ {\rm Ann} \  (p_*{\mathcal O}_X/{\mathcal O}_Y)$ denote the sheaf of conductors of $Y$. Let $y_1, \cdots , y_{m'}$ be the singular points of $Y$. Let ${\tilde{\mathcal O}}_{y_j}$ denote the normalisation of the local ring ${\mathcal O}_{y_j}$ at $y_j$ and 
$$\delta(y_j)= \ {\rm dim}_{\mathbb{C}} \ ({\tilde{\mathcal O}}_{y_j}
/{\mathcal O}_{y_j})\, , \ \delta(Y) = \sum_{j=1}^{m'} \delta (y_j)\, .$$  
For $j=1,\cdots, m'$, let $\{ D_j \}$ be disjoint (Cartier) divisors on $X$ with support of $D_j$ contained in $p^{-1} y_j$ and $p_* {\mathcal O}_X(-D_j) \subset {\mathcal C}$.  Choose 
$$f_j(E)= r(E)(d(D_j) - \delta(y_j))\, .$$ 
We fix a line bundle $L$ on $Y$ and define $L_0:= p^*L$.
 
Note that $\phi: E \to E \otimes L$ induces 
$$p_*\phi \colon p_*E \rightarrow p_*(E\otimes p^*L)\cong p_*E\otimes L\, .$$     

\begin{definition}
A GPH $(E, F(E), \phi)$ is a good GPH if it satisfies the following conditions. 
\begin{equation}\label{c3}
{\rm The} \ k-{\rm subspace} \  p_* F_j(E)  \ {\rm is \ an} \ {\mathcal O}_{y_j}-{\rm submodule \ of} \  p_* (E\otimes{\mathcal O}_{D_j})\, .
\end{equation}
\begin{equation}\label{c4}
(p_* \phi)_{y_j} (p_*(F_j(E))\subset (p_*(F_j(E))\otimes L_{y_j}\, . 
\end{equation}
\end{definition}

\subsection{The Hitchin morphism $h: N_Y(n,d,L) \to {\bf A}$} \hfill

There is a moduli space $H_X(n,d,L_0)$ of $L_0$-twisted semistable GPH of rank $n$ and degree $d$ on $X$ \cite{Bh1}. Let $M_X(n,d,L_0)$ be its closed subset corresponding to good GPH. There is a Hitchin morphism 
$$h_{H_X} : H_X(n,d,L_0) \to {\bf A} = \bigoplus_{i= 1}^n \ {\bf H}^0(X, L_0^i) \ {\rm defined \ by} \ (E, F(E), \phi) \mapsto \sum_{i=1}^n tr \phi^i\, .$$ 
The restriction of the Hitchin morphism to the subset  $M_X(n,d,L_0)$ induces a Hitchin morphism $h$ on its image $N_Y(n,d,L) \subset \mathcal{H}_Y(n,d,L)$ \cite{Bh1}. We note that the pointwise construction of this map can be extended to families of  $L_0$-twisted GPH, hence it defines the morphism $h$. \\
$\bullet$ The Hitchin morphism 
$$h: N_Y(n,d,L) \to {\bf A} $$ 
is a proper morphism (see \cite[Corollary 1.4]{Bh1} ).\\
$\bullet$ $N_Y(n,d,L)$ contains the open subset consisting of all Hitchin bundles  \cite[Proposition 2.8]{Bh1}. On the Hitchin bundles (i.e. Hitchin sheaves with underlying sheaf locally free), the Hitchin map $h$ coincides with the natural Hitchin morphism  and hence has values in $A'$. \\
$\bullet$  Since $N_Y(n,d,L)$ contains $\mathcal{H}'_Y(n,d,L)$,  the image of the rational forgetful map 
$$fr: N_Y(n,d,L) \to U(n,d)$$ 
contains $U'(n,d)$.    \\
$\bullet$  $N_Y(n,d,L)$ is invariant under the action of $\mathbb{C}^*$ defined by $$c(E,\phi) = (E, c \phi), c \in \mathbb{C}\, .$$\\
$\bullet$ Now, let $Y$ be an integral projective nodal curve. Then: \\
(1) The Hitchin map $h$ is surjective by \cite[Corollary 4.5]{Bh2}. \\
(2) The forgetful (rational) map $fr$ is dominant. \\
For $Y$ nodal, $U'(n,d)$ is an open dense subset of the variety $U(n,d)$. Since $fr(N_Y(n,d,L)) \supset  U'(n,d)$, it follows that the forgetful map $fr$ is dominant. In fact, Its restriction to the closure of $\mathcal{H}'_Y(n,d,L)$ in $N_Y(n,d,L)$ is dominant.
    
{\bf Remark : }      The Hitchin map is defined on the whole Hitchin moduli space in some cases. For example, suppose that  $Y$ has only one node $y$ with $x, z \in X$ lying over $y$, also $n=2$ and the determinant $\xi$ of $E$ is fixed. Then any (stable or semistable) torsion-free sheaf in the fixed determinant moduli space $U_{\xi}(n,d)$ is either a vector bundle $F$ or  $F= p_* E'$, the direct image of a vector bundle $E'$ (respectively stable or semistable) of same rank but degree $d-2$ and with a fixed determinant. A Hitchin field $\phi': E' \to E' \otimes p^*(L)$ induces a Hitchin field $\phi := p_* \phi':  p_*E' \to p_*(E' \otimes p^*(L)) = p_*E' \otimes L$. Thus the Hitchin map is defined on the whole fixed determinant moduli space  $\mathcal{H}(n, \xi, L)$.

\subsection{The fibre of the Hitchin map over a general $a \in {\bf A}'$}    \label{hitchinfibre}   \hfill

          For a general $a \in {\bf A'} = \oplus_{i=1}^n {\bf H}^0(Y, L^i)$, the spectral curve $Y_a$ is an $n$-sheeted covering of $Y$ unramified over the nodes of $Y$, it is a nodal curve with nodes precisely at points over the nodes of $Y$ (see \cite{Bh1}, \cite{Bh2} for more details). Let $h_{set}^{-1} (a)$ denote the fibre of $h_{set}$ over $a$ .  

We assume that $d(L) \ge 2g-2$.
               
\begin{proposition} \label{fibreainA'}
For a general $a \in A'$, there is a morphism $h_J$ from the compactified Jacobian $\bar{J}(Y_a)$ of $Y_a$ to $\mathcal{H}_Y(n,d,L)$. It maps bijectively onto $h_{set}^{-1} (a)$. One has $h^{-1}(a) = h_J(\bar{J}(Y_a)) \cap N_Y(n,d,L)$. 

          The  morphism $h_J$ maps $J(Y_a)$ into $\mathcal{H}'(n,d,L)$.  Its restriction to $J(Y_a)$ is an isomorphism on to its image. In particular the general Hitchin fibre $h^{-1}(a)$ is irreducible of dimension equal to  
      $$g(Y_a) =  \frac{n(n-1) deg(L)}{2} + n(g-1) +1\, .$$
\end{proposition}

\begin{proof}
               There is a universal torsion-free sheaf $\mathcal{P} \to \bar{J}(Y_a) \times Y_a$ which is locally free over $J(Y_a) \times Y_a$ and $\bar{J}(Y_a) \times Y'_a$, where $Y'_a$ denotes the subset of $Y_a$ corresponding to nonsingular points. Let $\pi: Y_a \to Y$ be the projection. Define 
              $$\mathbb{V} := (id \times \pi)_* \mathcal{P}\, .$$
It is a torsion-free sheaf of rank $n$ on $\bar{J}(Y_a) \times Y$. In fact for $y \in Y$, $\mathbb{V}_y :=  
\mathbb{V}_{\bar{J}(Y_a) \times y} = \oplus _{i=1}^n \mathcal{P}_{y_i}$ where 
$\mathcal{P}_{y_i}  = \mathcal{P}\vert_{\bar{J}(Y_a)\times y_i}$, where $y_i$ are the points of $Y_a$ lying over $y$.  There is a Hitchin field $ \theta_{\bar{J}}: \mathbb{V} \to L$ induced by the $(id \times \pi)_*  \mathcal{O}_{\bar{J}(Y_a) \times Y_a}$-structure on $\mathbb{V}$.  Thus we have a family 
$(\mathbb{V},  \theta_{\bar{J}})$ of $L$-valued Hitchin sheaves on $Y$ parametrised by $\bar{J}(Y_a)$.  Hence there is a  morphism $h_J$ from the compactified Jacobian $\bar{J}(Y_a)$ (of suitable degree) of the spectral curve $Y_a$ to ${\mathcal H}_Y(n,d,L)$. Since $\bar{J}(Y_a)$ is irreducible, $h_J(\bar{J}(Y_a))$ is irreducible.

                                        The morphism $h_J$ is injective by \cite[Lemma 4.6]{Bh2}.
The first two paragraphs of the proof of \cite[Proposition 4.7]{Bh2} prove that if $Y_a$ is an integral curve, then $h_J$ is a bijective morphism from $\bar{J}(Y_a)$ onto $h_{set}^{-1} (a)$ (the set of closed points of $h_J(\bar{J}(Y_a))$ is $h_{set}^{-1} (a)$). We note that $h_{set}^{-1} (a) \subset \mathcal{H}_Y(n,d,L)$ gets an algebraic structure as the image of $h_J$.  Clearly, $h^{-1}(a) = h_J(\bar{J}(Y_a)) \cap N_Y(n,d,L)$. 

            Since $\mathcal{P}\vert_{J(Y_a) \times Y_a}$ is a vector bundle, $\mathbb{V}\vert_{J(Y_a) \times Y}$ is a vector bundle. Hence $h_J(J(Y_a)) \subset \mathcal{H}'_Y(n,d,L)$ and $h_J(J(Y_a)) \subset h^{-1} (a)$.
 The restriction of the morphism $h_J$ to $J(Y_a)$ is an isomorphism onto image, note that $J(Y_a)$ is smooth. In particular, (for $Y$ an integral curve) the general Hitchin fibre is irreducible of dimension equal to the arithmetic genus of $Y_a$. Hence the dimension of the general fibre of the Hitchin map is given by   
      $$g(Y_a) = \frac{n(n-1) deg(L)}{2} + n(g-1) +1\, ,$$
(see \cite[Remark 3.2]{BNR}).       
\end{proof}

\subsection{Spectral curve and Hitchin Fibre over a general $a \in {\bf A}$}   \label{hitchinfibreA}   \hfill

      For simplicity of exposition, we assume that there is only one node $y$, the general case can be done similarly. Then $p_*(L_0) = I_y L_1 \subset L_1$ where $I_y$ is the ideal sheaf of $y$ and $L_1$ is a line bundle on $Y$ with  $p^*L_1 = L_0(x+z)$.  Let 
      $${\bf A}'_1 := \bigoplus_{i=1}^n {\bf H}^0(Y, L_1^i)\ .$$
Then ${\bf A} \subset {\bf A}'_1$. \\[3mm]
{\bf The Spectral curve}: \\
         Since ${\bf A}'$ is a closed subvariety of ${\bf A}$, a general $a \in {\bf A}$ belongs to ${\bf A} \setminus {\bf A}'$. 

 Now consider $a \in {\bf A} \subset {\bf A}'_1\ \ , a = (a_1, \cdots, a_n)$. Let $p_Y: {\bf P}(\mathcal{O} \oplus L_1^*) \to Y$ be the projection. Let $Y_a \subset {\bf P}(\mathcal{O} \oplus L_1^*)$ be the spectral curve associated to $a \in {\bf A}'_1$. We may assume that for a general $a \in {\bf A}$, the spectral curve $Y_a$ is integral. 

        We consider the Hitchin sheaves $(E, \phi)$ with $h(E, \phi) \in {\bf A}$ as Hitchin sheaves with $h(E, \phi) \in {\bf A}'_1$. Then $\phi$ and all $a_i$ are sections of $I_y (L_1)^i \subset L_1^i$, i.e. they  vanish at $y$. The local equation of $Y_a$ in a neighbourhood of the node $y$ is 
$$(x')^n + p_Y^*a_1 \cdot y' \cdot (x')^{n-1} + \cdots + p_Y^*a_n \cdot (y')^n\, ,$$
(see \cite[p. 172, \$ 3]{BNR} for details). Since all the $a_i$ vanish at  $y$, there is a unique point in $Y_a$ lying over $y$, given by $(x')^n =0$ (the point corresponding to $({\bf P}(\mathcal{O})_y$)) i.e., $Y_a$ is totally  ramified at $y$. \\[3mm]
{\bf The Hitchin fibre}:
 
        Since $Y_a$ is integral, there is a bijective homomorphism $h_J$ from $\bar{J}(Y_a)$ onto its image in $\mathcal{H}_Y(n,d,L_1)$ (by Proposition \ref{fibreainA'}) which gives an isomorphism of $J(Y_a)$ with its image.  To a torsion-free sheaf $M$ on $Y_a$, this morphism associates a  torsion-free sheaf $(p_Y\vert_{Y_a})_* M$ and the Hitchin  field is the natural homomorphism  $\phi: (p_Y\vert_{Y_a})_* M \to (p_Y\vert_{Y_a})_* M \otimes L_1 =  (p_Y\vert_{Y_a})_*((p_Y\vert_{Y_a})^*(L_1)\otimes  M)$ given by the section $x'$ of $(p_Y\vert_{Y_a})^*L_1$. Since this section $x'$ vanishes over $y$, it follows that  $\phi((p_Y\vert_{Y_a})_* M)  \subset  (p_Y\vert_{Y_a})_* M \otimes I_y L_1= (p_Y\vert_{Y_a})_* M \otimes p_*(L_0)$. In fact, $\phi$ has values in $p_*L_0$ if and only if $x'$ vanishes over $y$. Thus $h_J(\bar{J}(Y_a))$ maps into $\mathcal{H}(n,d, p_* L_0) \subset \mathcal{H}_Y(n,d,L_1)$. 
        
         Thus  the fibre  $h_{set}^{-1} (a)$ is in bijective correspondence with $\overline{J}(Y_a)$ and therefore has dimension $\frac{n(n-1) deg(L_1)}{2} + n(g-1) +1$ by \cite[Remark 3.2]{BNR}.  One has $d(L_1) = d(p_* L_0) +1$. 
         
         Let $(E, \phi) \in h_{set}^{-1} (a)$ for a general $a \in {\bf A} \setminus {\bf A}' $. Suppose that $E$ is a vector bundle. Since $E$ is vector bundle and $\phi$ has values in $L$,  $h_{set}(E, \phi) = a \in {\bf A}' $. But  $a \notin {\bf A}' $, it follows that $E$ is not locally free and hence corresponds to an element of $\overline{J}(Y_a) \setminus J(Y_a)$. Thus the fibre $h_{set}^{-1} (a)$ of the Hitchin set map on $\mathcal{H}(n,d, L)$ is contained in the image of $\overline{J}(Y_a) \setminus {J}(Y_a)$, the same is true of the general Hitchin fibre $h^{-1} (a)$ in $N_Y(n,d,L)$ for $a \in  {\bf A} \setminus {\bf A}' $.

Similarly, for a nodal curve $Y$ with $m'$ nodes we have the following result.    

\begin{proposition} \label{fibreainA}
   For a general $a \in {\bf A}$, there is a morphism $h_J: \bar{J}(Y_a) \to \mathcal{H}(n,d, p_* L_0)$, whose restriction is an isomorphism of $J(Y_a)$ onto its image. The general Hitchin fibre $h^{-1} (a)$ in $N_Y(n,d,L)$ is contained in $h_J(\overline{J}(Y_a) \setminus {J}(Y_a))$ and hence has dimension at most  
      $$g(Y_a) -1 = \frac{n(n-1) (deg(p_*L_0) + m')}{2} + n(g-1)\, .$$
      
\end{proposition}

{\bf Notation 2} \label{It}    \\
(1) For a curve with one node, define $\mathcal{H}_i(n,d,L)$ as the subset of $\mathcal{H}(n,d,L)$ corresponding to Hitchin sheaves of local type $n- i$ at the node. Let $N_i(n,d,L) = \mathcal{H}_i(n,d,L) \cap N_Y(n,d,L)$.\\
(2)  Let $y_1, \cdots, y_{m'}$ be the nodes of $Y$. Then a local type (defined in Definition \ref{localtype}) can be written as an $m'$ tuple $\bar{b} = (b_1, \cdots , b_{m'}), b_j$ being the local type at $y_j$. For a local type $\bar{b}$, define $\mathcal{H}_{\bar b}(n,d,L)$ as the subsheme of $\mathcal{H}(n,d,L)$ corresponding to Hitchin sheaves of local type $\bar{b}$.\\
(3)   Let $I_{m'} = [1,\cdots , m']$ denote the set of first $m'$ natural numbers. For $1 \le t \le m'$, let $I_t = [i_1, \cdots, i_t] \subset I_{m'}$ denote a subset with $t$ elements. Let $p_{I_t}: Y_{I_t} \to Y$ be the partial normalisation of $Y$ obtained by blowing up nodes $y_{i_1}, \cdots, y_{i_t}$, so $Y_{I_{m'}} = X$.  One has $g(Y_{I_t}) = g - t,  d((p_{I_t})_* \mathcal{O}_{Y_{I_t}}) 
= t$. Let 
$${\bf A}_{I_t} = \sum_{i=1}^{n}  {\bf H}^0(Y_{I_t}, (p_{I_t})^* L^i) = \sum_{i=1}^{n}  {\bf H}^0(Y, (p_{I_t})_*(p_{I_t})^* L^i) = \sum_{i=1}^{n}  {\bf H}^0(Y,  L^i \otimes (p_{I_t})_*\mathcal{O}_{Y_{I_t}}) \subset {\bf A}\, .$$
   
\begin{proposition} \label{ontoA}   

(1)  Let $\bar{b} = (b_1, \cdots, b_{m'})$ be a local type, with  $b_j= n$ for the nodes $y_j,  j \in I_t$ and $b_j= 0$ for the nodes $y_j,  j \notin I_t$. Then the Hitchin morphism is well defined on $\mathcal{H}_{\bar b}(n,d,L)$ and it maps $\mathcal{H}_{\bar b}(n,d,L)$ into ${\bf A}_{I_t}$. 

(2) Let $\mathcal{H}_0(n,d,L) \subset \mathcal{H}(n,d,L)$ be the subset corresponding to Hitchin sheaves with local type ${\bar b}$ with $b_j= n$ at all the nodes $y_j$.  Then the Hitchin morphism maps $\mathcal{H}_0(n,d,L)$ onto ${\bf A}$.  In particular, Hitchin morphism is surjective. 
\end{proposition}

\begin{proof}
      
(1)   Consider  $(E', \phi') \in \mathcal{H}(n,d,L)$ with $b_j(E') = n$ for all the nodes $y_j,  j \in I_t$ and $b_j(E') = 0$ for all the nodes $y_j,  j \notin I_t$.  Then there is a GPH $(E, F(E), \phi)$ on $X$ mapping to $(E', \phi')$ such that $F_j(E) = E_{x_j}$ for $ j \in I_t$ and  $F_j(E) = E'_y$ for $j \notin I_t$. Hence $F_j(E)$ is $\phi$-invariant for all $j$. It follows that $(E, F(E), \phi)$ is a good GPH and $(E', \phi') \in N(n,d,L)$. 

  We prove that $\mathcal{H}_{\bar b}(n,d,L) \cong \mathcal{H}_{Y_{I_t}}(n, d- t n, p_{I_t}^*L)$ where the latter is the moduli space of $p_{I_t}^*L$-valued Hitchin bundles on $Y_{I_t}$ of rank  $n$, degree $d- n t$ on $Y_{I_t}$. There is a morphism 
$$\mathcal{H}_{Y_{I_t}}(n, d- t n, p_{I_t}^*L) \to \mathcal{H}_{\bar b}(n, d, L) 
\ {\rm defined \ by} \  (E_t, \phi_t) \mapsto (E', \phi')\, ,$$ 
where $E' = (p_{I_t})_*E_t$ and $\phi'$ is defined as follows. The homomorphism $\phi_t: E_t \to E_t \otimes p_{I_t}^*L$  gives a homomorphism $(p_{I_t})_* \phi:  (p_{I_t})_*E_t \to (p_{I_t})_*(E_t \otimes p_{I_t}^*L) = (p_{I_t})_*E _t \otimes L\, ,$ and $\phi' := (p_{I_t})_* \phi_t: E' \to E' \otimes L$. The inverse map is defined as follows. Given $(E' ,\phi') \in  \mathcal{H}_{\bar b}(n, d, L)$, we  define $E_t = (p_{I_t})^*E'/tor$ where $tor$ denotes the torsion subsheaf. We have $(p_{I_t})^*\phi': (p_{I_t})^*E' \to (p_{I_t})^*E' \otimes (p_{I_t})^*L$ which, on going modulo torsion, gives $\phi: E_t \to E_t \otimes (p_{I_t})^*L$. 

           There is a commutative diagram
$$
\begin{array}{ccc}
\mathcal{H}_{Y_{I_t}}(n, d- t n, (p_{I_t})^*L) & \stackrel{\cong}{\longrightarrow} & \mathcal{H}_{\bar b}(n, d, L) \\
{} & {} & {}\\
h_{Y_{I_t}} \downarrow & {} & \downarrow h\\
{} & {} & {}\\
{\bf A}_{I_t} &  \hookrightarrow  & {\bf A}
\end{array}
$$  
from which the Part (1) follows. \\
(2)    For generic $a \in {\bf A} = \oplus \ {\bf H}^0(X, L_0^i)$, there is a nonsingular spectral curve $X_a$ and the fibre of the Hitchin map is isomorphic to the Jacobian $J(X_a)$ of the spectral curve.  Hence the Hitchin map $h_X$ is dominant. Since the Hitchin map is proper, its image is a closed subset. Since ${\bf A}$ is irreducible, this implies that $h_X$ is surjective. In view of the above diagram (for $I_t = I_{m'}$), this implies that $h: \mathcal{H}_0(n, d,L) \to {\bf A}$ is surjective.                 
\end{proof}

\subsection{Examples}   \hfill          

{\bf Example 1:  Rank $n=1$, degree $d$, $d(L) \ge 2g-1\,$ or  $L = \omega_Y$.} 
\\[3mm]
{\bf Case $Y$ nonsingular}:      For $m' =0$,  $\mathcal{H}(1,d,L) = J(Y) \times {\bf H}^0(Y, L)$ and hence it is a nonsingular quasi-projective variety.\\
{\bf Case $Y$ with one node}: 

        The case when $Y$ has only one node was worked out in detail in \cite[Subsection 2.4]{Bh1}.  One has 
$$\mathcal{H}(1,d,L) = [J(Y) \times {\bf H}^0(Y,L)] \coprod [(\bar{J}(Y) \setminus J(Y)) \times {\bf H}^0(Y, L\otimes p_* \mathcal{O}_X)]\, .$$
Thus $\mathcal{H}(1,d,L)$ has two irreducible components, namely $\bar{J}(Y) \times  {\bf H}^0(Y,L)$ and $(\bar{J}(Y) \setminus J(Y)) \times {\bf H}^0(Y, L\otimes p_* \mathcal{O}_X)\, .$  The first component has dimension $d(L)+1$ if $d(L) \ge 2g-1$ and dimension $2g$ if $L = \omega_Y$, it contains all Hitchin line bundles. The second component has dimension $d(L)+1$ if $d(L) \ge 2g-1$ and dimension $2g-1$ if $L = \omega_Y$. It consists of only Hitchin sheaves with underlying sheaves not locally free.

Thus $\mathcal{H}(1,d,L)$ is connected but not irreducible, it has $2$ components. It is equi-dimensional of dimension $d(L) +1$ if $d(L) \ge 2g-1$. It has components of dimension $2g$ and $2g-1$ if $L = \omega_Y$.  

    The Hitchin map $h: \mathcal{H}(1,d,L) \to {\bf H}^0(Y, L\otimes p_* \mathcal{O}_X)$ is defined by $(N, \phi) \mapsto \phi$. For $\phi \in H^0(Y,L)$, the fibre of the Hitchin map is isomorphic to $\bar{J}(Y)$. For $\phi \in  ({\bf H}^0(Y, L\otimes p_* \mathcal{O}_X) \setminus \bf{H}^0(Y,L))$, the fibre of the Hitchin map is isomorphic to $\bar{J}(Y) \setminus J(Y)$ \cite[Lemma 2.13]{Bh1}.\\
{\bf Case $Y$ with $m' \ge 2$ nodes:} 
            
            Let $I_t, p^{I_t},  Y_{I_t}$ be as defined in Notation 2  and let $J(Y_{I_t})$ be the generalised Jacobian of $Y_{I_t}$.  Define 
 $$p^{I_t}_* J(Y_{I_t}) := \{ N \in \bar{J}(Y) \vert N = p^{I_t}_* N', N' \in J(Y_{I_t}) \}\, ,$$
 it has dimension $g-t$. 
 
 One has  
$$
\bar{J}(Y) = J(Y) \coprod (\coprod_{\stackrel{I_t \subset I_m'}{1\le t \le m'}} p^{I_t}_* J(Y_{I_t}) )\, .$$
    Let $\mathcal{H}^{I_t}_Y(n,d,L)$  be the inverse images of $p^{I_t}_* J(Y_{I_t})$ under the forgetful morphism $fr$. 

               Let $N \in  p^{I_t}_* J(Y_{I_t})$ for some $I_t$  i. e.,  $N = p^{I_t}_* N_{I_t}\, , N_{I_t} \in J(Y_{I_t})$. Then 
$$h^0(Y, Hom(N, N\otimes L)) =  h^0(Y, Hom(p^{I_t}_*N_{I_t}, p^{I_t}_*N_{I_t} \otimes L)) = h^0(Y, p^{I_t}_* Hom(N_{I_t}, N_{I_t}) \otimes L))\, .$$
One has $h^0(Y, p^{I_t}_* Hom(N_{I_t}, N_{I_t}) \otimes L))  = h^0(Y, p^{I_t}_*(\mathcal{O}_{Y_{I_t}}) \otimes L) \, .$
Hence $\mathcal{H}^{I_t}_Y(n,d,L) = p^{I_t}_* J(Y_{I_t}) \times {\bf H}^0(Y, L \otimes p^{I_t}_*(\mathcal{O}_{Y_{I_t}}) ) \, .$ It follows that $\mathcal{H}(1,d,L)$ has $2^{m'}$ components which are the closures of 
$\mathcal{H}'(1,d,L)$ and $\mathcal{H}^{I_t}_Y(n,d,L), I_t \subseteq I_{m'}, 1\le t \le m'$ in $\mathcal{H}(1,d,L)$. 
          
         We compute the dimensions for $L = \omega_Y$. For $N \in J(Y), h^0(Y, \omega_Y) = g$ so that dim $\mathcal{H}'(1,d,\omega_Y) = 2 g$. For $N \in  p^{I_t}_* J(Y_{I_t})$, dim $h^0(Y, p^{I_t}_*(\mathcal{O}_{Y_{I_t}}) \otimes \omega_Y) = 2g-2 + t + 1 -g = g +t -1$, dim $\mathcal{H}^{I_t}_Y(n,d, \omega_Y) = (g-t) + (g+t-1) =2g -1$. Thus there is one component of dimension $2g$ which contains all Hitchin line bundles, all other components are of dimension $2g-1$ and contain no Hitchin bundles.
         
         In the case $d(L) \ge 2g-1$, one can similarly see that there is one component of dimension $d(L)+1$ which contains  all Hitchin line bundles, all other components are of dimension $d(L)+1$ and contain no Hitchin bundles.  In this case, $\mathcal{H}(1,d,L)$ is equi-dimensional.

{\bf Example 2: One node, rank $n=2$, degree $d=1\, , L = \omega_Y$.} \hfill

In this case 
$$dim \ A' = h^0(Y, \omega_Y) + h^0(Y, \omega_Y^2) = g +(3g -3) = 4g - 3\, ; $$ 
$$dim \ A =  h^0(Y, \omega_Y\otimes p_*\mathcal{O}_X) + h^0(Y, \omega_Y^2\otimes p_*\mathcal{O}_X) = g+ (3g -2) = 4g - 2\, .$$
            Let $U(2,1)$ denote the moduli space of torsion-free sheaves of rank $2$ and degree $1$ on $Y$, it is a seminormal projective variety. Note that semistability is equivalent to stability.  Hence its subset  $U'(2,1)$, corresponding to vector bundles, is a smooth quasi projective variety. Let $U_1(2,1)$ and $U_0(2,1)$ denote the subsets of $U(2,1)$ corresponding to sheaves of local type $1$ and $2$ respectively. One has
$$dim \ U'(2,1) = 4g -3, \ dim \ U_1(2,1) = 4g - 4, \ dim \ U_0(2,1) = 4g -7\, .$$   

          The moduli space ${\mathcal H}'(2,1, \omega_Y)$ contains an open subset ${\mathcal H}'_s(2,1, \omega_Y)$ consisting of $(E,\phi)$ with $E$ a stable vector bundle. This subset is a fibration  over $U'(2,1)$ with fibre over $E\in U'(2,1)$ isomorphic to $H^0(Y, (End E) \otimes \omega_Y)$. For $E \in U'(2,1)$, we have $h^0(Y, (End E) \otimes \omega_Y) = h^1(Y, End E) = 4g -3$ so that 
$${\rm dim} \  \  {\mathcal H}'(2,1, \omega_Y) = 8g - 6\, .$$
          The dimension of the fibre of the Hitchin map $h \vert_{{\mathcal H}' (2,1, \omega_Y)}$ for a general $a \in A'$ equals $g(Y_a) = 4g -3 > 0$. Hence the Hitchin morphism ${\mathcal H}' (2,1, \omega_Y) \to {\bf A}'$ is dominant.
         
             For $E \in U_0(2,1), E = p_* E_1$ for a unique stable vector bundle $E_1$ on $X$. We have  
$h^0(X, (End E_1) \otimes \omega_X(x+z)) = h^0(X, (End E_1) \otimes p^*\omega_Y) = h^0(Y, p_*(End E_1) \otimes \omega_Y) = h^0(Y, (End E) \otimes \omega_Y)$. Since $End E_1$ is semistable (as $X$ is smooth) vector bundle of degree $0$, $h^1(X, (End E_1) \otimes \omega_X(x+z)) = h^0(X, (End E_1) (-x-z)) = 0$ so that $h^0(X, (End E_1) \otimes \omega_X(x+z)) = 4g(X) + 4 = 4g$ by Riemann-Roch theorem. Hence for the variety ${\mathcal H}_0(2,1,\omega_Y)$ we have 
$$dim \  {\mathcal H}_0(2,1, \omega_Y) = 8g - 7\, .$$ 
  The Hitchin morphism maps ${\mathcal H}_0(2,1, \omega_Y)$ onto $A$ (Proposition \ref{ontoA}). One has 
$$dim \ (h^{-1}(a) \cap {\mathcal H}_0(2,1, \omega_Y)) = g(X_a)= 4g - 5\, .$$ 

       For $E \in U_1(2,1)$, one has $d(End E)^* = -1$ \cite[Lemma 2.5]{Bh2}. The torsion-free sheaf $End E$ is semistable for $E$ in an open subset of $U_1(2,1)$ and hence $(End E)^*$ is semistable for $E$ in this open subset \cite{ABhS}. Therefore $h^1(Y,(End E) \otimes \omega_Y) = h^0(Y, (End(E)^*) = 0$ and  $h^0(Y,(End E) \otimes \omega_Y) = 4g -3$. If ${\mathcal H}'_1(n,d,L) \subset {\mathcal H}_1(n,d,L)$ is the inverse image of the open subset of $U_1(2,1)$ in ${\mathcal H}_1(n,d,L)$, one has 
$$dim \  {\mathcal H}'_1(2,1, \omega_Y) = 8g - 7\, .$$        
          
          For a general $a \in A \setminus A'$, the fibre $h^{-1}(a)$ is contained in ${\mathcal N}_1(2,1, \omega_Y) \cup {\mathcal H}_0(2,1, \omega_Y)$. 
By Proposition \ref{fibreainA}, 
$$dim \ h^{-1} (a) \le g(Y_a) -1 \le  4g -2\, .$$ 

      Alternatively, one can compute $g(Y_a)$ using Hurwitz formula. By \cite[Remark 3.3]{BNR}, the number of ramification points on $Y$ equals $d(p^*L_1\otimes \mathcal{O}(1)) =4g +2$ (the unique node has multiplicity $2$). Let $\widetilde{Y_a}$ be the desingularisation of $Y_a$. It is a double covering over $X$ ramified over $4g +2$ points. By Hurwitz theorem, $2g(\widetilde{Y_a}) - 2 = 2(2 g(X) -2) + 4g +2$ so that  $g(\widetilde{Y_a}) =  4g(X) = 4g - 4.$ Then $g(Y_a) = 4g -1$. 
       
        The case when there is only one node, rank $n=2$, degree $d=1\,, d(L) \ge 2g-1$ can be worked out similarly, we omit the details. One gets $dim \  {\mathcal H}'(2,1, L) = dim \  {\mathcal H}_0(2,1, L) = 4 d(L)
+1$.
     
\subsection{The closed subsets ${\bf H}(n,d,\omega_Y)$ and ${\bf H}(n,d,L), d(L) \ge 2g-1$.} \label{bfH(n,d,L)} \hfill
  
            In the general case, the moduli space $U^{' s}(n,d)$ of stable vector bundles of rank $n$ and degree $d$ on $Y$ is smooth and irreducible of dimension $n^2(g-1)+1$. Any Hitchin bundle $(E, \phi)$ with $E \in U^{' s}(n,d)$ is stable.   We denote by ${\bf H}(n,d, L)$ the closure of $\mathcal{H}'(n,d,L)$ in $N_Y(n,d,L)$.
            
 {\bf Case $L = \omega_Y$:}           One has $h^1(Y, End E \otimes \omega_Y) = h^0(Y, (End E)^*) = h^0(Y, (End E)) =1$ as $E$ is stable. Hence $h^0(Y, End(E) \otimes \omega_Y) = 1 + n^2(2g -2) + n^2(1-g) = n^2(g-1) +1$. It follows that the subset $\mathcal{H}'_s(n,d,\omega_Y)$ of $\mathcal{H}(n,d,\omega_Y)$, corresponding to Hitchin bundles with the underlying bundle stable, is an irreducible smooth open subset with
      $${\rm dim} \ \mathcal{H}'_s(n,d,\omega_Y) = 2n^2(g-1)+2\, . $$ 
 
             We have a proper Hitchin morphism $h: N_Y(n,d,\omega_Y) \to {\bf A}$. It maps $\mathcal{H}'(n,d,\omega_Y)$ into the closed subset ${\bf A}'$ and hence it maps the closure ${\bf H}(n,d, \omega_Y)$ into ${\bf A}'$.  Since  ${\bf H}(n,d, \omega_Y)$ is a closed subset of $N_Y(n,d, \omega_Y)$, the restriction of $h$ to ${\bf H}(n,d, \omega_Y)$ is a proper morphism. For a general $a \in {\bf A}', J(Y)$ maps into $\mathcal{H}'$ under the morphism $h_J$ of Proposition \ref{fibreainA'}. In particular this implies that $h\vert_{\mathcal{H}'}(n,d, \omega_Y): \mathcal{H}' (n,d, \omega_Y) \to {\bf A}'$ is dominant. Since ${\bf H}(n,d, \omega_Y)$ is a closed subset, its image in ${\bf A}'$ is a closed subset, it follows that the restriction 
             $$h_{\bf H}: {\bf H}(n,d, \omega_Y) \rightarrow {\bf A}'$$ 
is surjective. 
    
               Let $U^s_0(n,d) \subset U^s(n,d)$ be the subset consisting of stable torsion-free sheaves $E$ which are direct images of (stable) vector bundles $E_0$ on the normalisation or equivalently $E$ having local type $n$ at all $m'$ nodes.  It is isomorphic to $U^s_X(n,d)$. For $E \in U^s_0(n,d)$, any $(E, \phi)$ is stable, but $d(End(E)) = n^2 m', d((End E)^*) = - n^2 m'$ so that $h^0(Y, (End (E)^*)) = 0$.  By Serre duality, $h^1(Y, End (E) \otimes \omega_Y) = 0$. Then $h^0(Y, End (E) \otimes \omega_Y) = n^2 m' + n^2(2g-2) + n^2(1-g) = n^2(g-1) +n^2 m'$.  It follows that the inverse image $\mathcal{H}_{0, s}(n,d,\omega_Y)$ of  $U^s_0(n,d)$ in $\mathcal{H}(n,d,\omega_Y)$ (under the forgetful map) is an irreducible smooth closed subset of $fr^{-1} U^s_0(n,d)$ with
      $${\rm dim} \ \mathcal{H}_{0, s}(n,d,\omega_Y) = 2n^2(g-1)+1 \, . $$  
      
      We claim that ${\bf H}(n,d, \omega_Y)$ does not contain $\mathcal{H}_0(n,d, \omega_Y)$. The Hitchin morphism $h$ maps ${\bf H}(n,d, \omega_Y)$ into ${\bf A}'$. However, the restriction of $h$ to $\mathcal{H}_0(n,d, \omega_Y)$ surjects onto ${\bf A}$ by Proposition \ref{ontoA}. This proves the claim. 
            
{\bf Case  $L$ is a line bundle with $d(L) \ge 2g -1$:} 
      
      There is a non-empty open dense subset ${\bf U} \subset U(n,d)$ consisting of stable vector bundles $E$ such that the endomorphism bundle $End (E)$ is semistable. This includes the nonempty set of strongly stable vector bundles $E$ i.e., $E$ such that $p^*(E)$ is stable. For $E \in {\bf U}$, one has $h^1(Y, End(E) \otimes L) =  0$ and hence $h^0(Y, End(E) \otimes L) =  n^2(d(L)+ 1-g)$. It follows that the inverse image ${\bf H}'(n,d,L)$ of ${\bf U}$ under $fr$ has dimension  $n^2(d(L)+1)$. 
      
      Similarly, there is a non-empty open dense subset ${\bf U}_0 \subset U_0(n,d)$ consisting of stable torsion-free sheaves $E$ such that the sheaf $Hom(E, E) = End (E)$ is semistable. This includes the nonempty set of direct images $p_* E_0$ of strongly stable vector bundles $E_0$ on $X$. For $E \in {\bf U}_0$, one has $h^1(Y, End(E) \otimes L) = h^1(X, End (E_0) \otimes p^*L) =  0$ and hence $h^0(Y, End(E) \otimes L) =  n^2(d(L)+ m' +1-g)$. It follows that the inverse image ${\bf H}_0(n,d,L)$ of ${\bf U}_0$ under the forgetful map $fr$ has dimension  $(n^2(g-1) - n^2 m') + (n^2(d(L)+m' +1 -g))= n^2(d(L) +1)$. Note that ${\bf H}_0(n,d,L)$ is a nonempty irreducible closed subset of the closed subset $fr^{-1} U_0(n,d)$ of $\mathcal{H}(n,d,L)$ and it has the same dimension as the irreducible nonempty open subset ${\bf H}'(n,d,L)$ of $\mathcal{H}(n,d,L)$. 
It follows that $\mathcal{H}(n,d,L)$ has at least two irreducible components.

\section{Very stable torsion-free sheaves} \hfill

Generalising  Drinfeld's definition, we make the following definition.

\begin{definition}\label{verystable} 
   Let $k \supset \mathbb{C}$ be a field. Let $Y$ be an irreducible complex projective nodal curve defined over $\mathbb{C}$. A torsion-free sheaf $E$ on $Y\otimes_{\mathbb{C}} k$ is called very stable (respectively $L$-very stable) if there is no nonzero nilpotent homomorphism of $\mathcal{O}_{Y\otimes_{\mathbb{C}} k}$-modules $u: E \to E \otimes_{\mathcal{O}_Y} \omega_Y$ (respectively $u: E \to E \otimes_{\mathcal {O}_Y} L)$ .
\end{definition}

\begin{proposition} \label{LaProp3.5}
(1)  For $g \ge 1$, every very stable torsion-free sheaf on $Y$ is semistable.\\
(2)  For $g \ge 2$, every very stable torsion-free sheaf on $Y$ is stable. \\
(3)  For $g =1$, a very stable torsion-free sheaf $E$ is either stable or it has no  subsheaf $N$ non locally free at a node $y$,  with $E/N$ a torsion-free sheaf which is non locally free at the node $y$, contradicting the stability.\\   
\end{proposition} 

\begin{proof}      
(1)    Let $E$ be a very stable torsion-free sheaf. Suppose that $E$ contains a torsion-free subsheaf $N$ of rank $r(N)$ with a torsion-free quotient $Q = E/N$ of rank $r(Q)$ such that their degrees satisfy the inequality $\mu(N) \ge \mu(Q)$. Let $F$ denote the sheaf $Hom(Q, N \otimes \omega_Y)$. Since $\omega_Y$ is locally free, by \cite[Lemma 2.5(B)(1)]{Bh2},
$$
\begin{array}{lll}
d(F) & = & d(N \otimes \omega_Y) r(Q) -  r(N \otimes \omega_Y) d(Q) + \sum_j b_j(N\otimes \omega_Y) b_j(Q)\\ 
{}     & = & r(N) r(Q) (\mu(N) - \mu(Q) + (2g-2)) + \sum_j b_j(N) b_j(Q)\, ,
\end{array}
$$  
where $b_j(N)$ and $b_j(Q)$ denote the local types of $N$ and $Q$ respectively at the node $y_j$, with $j$ varying over the nodes.     
One has $\sum_j b_j(N) b_j(Q) \ge 0$ and 
$$h^0(F) \ge \chi(F) = d(F) + r(N) r(Q) (1-g) = r(N) r(Q) (\mu(N) - \mu(Q) + (g-1)) + \sum_j b_j(N) b_j(Q)\, .$$
 
If $N$ contradicts the semistability of $E$,  $\mu(N) > \mu(Q)$. Then for $g \ge 1$ one has $h^0(F) > 0$ so that there exits a nonzero homomorphism $v: Q \to N \otimes \omega_Y$. Define $u = 
(i\otimes Id_{\omega_Y}) \circ v \circ q$. Then $u$ is a non trivial nilpotent homomorphism, in fact $u^2=0$. This contradicts the very stability of $E$ proving Part (1). \\
(2) We continue as in Part (1). If  $N$ contradicts the stability of $E$ but satisfies the semistability condition,  $\mu(N) = \mu(Q)$. Then for $g \ge 2$, we get $h^0(F) >0$ and we can define $u$ as in the proof of Part (1)  giving a contradiction and proving Part (2).    \\
(3) In the proof of Part (1), suppose that $N$ contradicts the stability of $E$ but satisfies the semistability condition  $\mu(N) = \mu(Q)$. Then for $g =1$, we get $h^0(F) \ge \sum_j b_j(N) b_j(Q) > 0$ if $b_j(N)>0, b_j(Q) > 0$ for some $j$. Then we can define $u$ as in the proof of Part (1)  giving a contradiction and proving Part (3). 

  We remark that if $E$ is locally free  and $N \subset E$ is a sub sheaf with $E/N$  torsion-free, then $N$ is non locally free at a node $y$ if and only if $E/N$ is non locally free at the node $y$.
\end{proof}

 For $g=0$, there are very stable vector bundles  which are not semistable eg. $O_Y \oplus O_Y(1)$ (see \cite{La}).

\begin{proposition} \label{verystablenonempty}
      (1)  For $g \ge 1$, there is an open dense subset $U \subset U(n,d)$ such that any $E_0 \in U$ is very stable i.e., every nilpotent homomorphism $u: E_0 \to E_0 \otimes \omega_Y$ is $0$.  \\
       (2) This open subset $U$ is contained in $U^{'s}(n,d)$ for $g \ge 2$. 
\end{proposition}

\begin{proof} 
 (1)     In the notations of Definition \ref{La1.7} and Equation \eqref{La1.12}, for all nilpotent types $(\nu_{\bullet},\lambda_{\bullet})$, the closures of ($\Sigma_{(\nu_{\bullet},\lambda_{\bullet})} \cong \Lambda_{(\nu_{\bullet},\lambda_{\bullet})} $) determine closed subsets $C_{(\nu_{\bullet},\lambda_{\bullet})}$ of the moduli space of Higgs bundles. These are the (finitely many) components of the nilpotent cone, the inverse image of $0$ under the Hitchin map. The forgetful map $fr$ from the null cone to the moduli space $U'(n,d)$ of vector bundles of rank $n$ and degree $d$ is surjective as it has the null section defined by $(E \to (E,0))$. The forgetful map $fr$ is generically finite as both the spaces are of same dimension $n^2(g-1) +1$. Thus there is an open subset $U \subset U'(n,d)$ such that for $E \in U,  f^{-1} (E)$ is finite.  If $(E, \phi) \in fr^{-1} (E), \phi \neq 0$, then  $(E,\mathbb{C} \phi) \subset fr^{-1} (E)$, a contradiction. Thus  $fr^{-1}{E} = (E, 0)$ for all $E \in U$ and $U \subset U'(n,d)$ is the open subset of very stable vector bundles. Since $U'(n,d)$ is irreducible (by a result of Narasimhan and Newstead \cite{N}), $U$ is dense in $U'(n,d)$. Since $U(n,d)$ is irreducible by a result of Rego (\cite{R}) and $U'(n,d)$ is an open subset of it, it follows that $U$ is open and dense in $U(n,d)$.
        
 (2)  Part (2) follows from the previous part and Proposition \ref{LaProp3.5}. 
      
 \end{proof}
{\bf Remark} Rego is credited to have proved the irreducibility of $U(n,d)$. Rego's Theorem says the following:  \\ 
Let $X$ be a singular integral curve of genus $g \ge 1$ defined over an algebraically closed field. If $X$ is embeddable in a smooth surface $Z$, then $U(n,d)$ is irreducible.   
       
       The irreducibility of $U(n,d)$ for an integral nodal projective curve can also be proved by using the irreducibility of the moduli space of $1$-semistable GPBs (or $1$-semistable GPS, the Generalised parabolic sheaves defined by Narasimhan M.S. and Sun X.) and using the surjective morphism from it to $U(n,d)$.

\subsection{The restriction $h_V$ of the Hitchin map}  \hfill
                          
             If $E$ is a stable torsion-free sheaf of rank $n$ and degree $d$, then $(E, \phi)$ is stable for any $\phi$. Hence for  $V_E := {\bf H}^0(Y, End E \otimes L)$,  there is an embedding
       $$i_V: V_E \hookrightarrow \mathcal{H}_Y(n,d,L) \ ; \  \phi \mapsto (E,\phi)\, ,$$  
we identify $V_E$ with $i_V(V_E)$.  The subsets  $V_E$ and $V_{E,1}:= V_E \cap N_Y(n,d,L)$ are  invariant under the action of $\mathbb{C}^*$ defined by $c(E,\phi) = (E, c \phi), c \in \mathbb{C}$. 

 Let  
$$h_V: V_{E,1} = V_E \cap N_Y(n,d,L) \longrightarrow {\bf A}$$ 
be the restriction of the Hitchin map $h$ to $V_{E,1}$.                  

If $E$ is a vector bundle, then $V_E \subset N_Y(n,d,L)$ 
\cite[Corollary 1.4 (3)]{Bh1} so that $V_{E,1} := V_E \cap N_Y(n,d,L) = V_E$. Moreover, 
$$h_V :   V_E \longrightarrow {\bf A'} \, .$$ 

\begin{proposition} \label{u1}
      Let $g \ge 2$.  If $d(L) > 2g-2$,  then there are no $L$-very stable vector bundles.  
\end{proposition}

\begin{proof}
 Let $d(L) > 2g -2$. Then for $i \ge 1$, one has $h^0(Y, L^i)= \chi(L^i)= i d(L) - (g-1)$, so that 
 $${\rm dim} \  A' \ = \sum_{i=1}^{n} \ h^0(Y, L^i)= \sum_{i=1}^n  i d(L) -n(g-1) = \frac{n(n+1)}{2} d(L) - n(g-1)\, .$$
For a vector bundle $E$, we have   
${\rm dim} \ V_E \ = \ h^0((End E)\otimes L) \ge \chi((End E)\otimes L)\, .$ 
Since $d((End E)\otimes L) = n^2 d(L)$, one has 
$ {\rm dim} \ V_E \ \ge \ n^2 d(L) - n^2(g-1) \, .$
Thus 
${\rm dim} \ V_E - {\rm dim} \ A' \ \ge n^2 d(L) - n^2(g-1) - \frac{n(n+1)}{2} d(L)+ n(g-1)$ so that 
$${\rm dim} \ V_E - {\rm dim} \ A' \ \ge \frac{n(n-1)}{2} (d(L) - 2(g-1)) >0\, .$$
            Since ${\rm dim} \ V_E > {\rm dim} \ A'$, every irreducible component of the fibre of $h_V$ has positive dimension. The fibre over $0 \in A'$ consists of $(E, \phi)$ with $\phi$ nilpotent. It follows that $E$ admits a nonzero nilpotent Hitchin field with values in $L$, i.e. $E$ is not $L$-very stable. 
\end{proof}

\begin{proposition} \label{u2} 
(1)  For $j = 1, \cdots, m'$, let $b_j$ be integers with $0 \le b_j \le n$. Assume that $E$ is a torsion-free sheaf of local type $\bar{b} = (b_1, \cdots, b_{m'})$ such that $V_E = V_{E,1}$. If 
  \begin{equation}  \label{c2}
   d(L) > 2g-2 + \frac{2}{n(n-1)} (m' n - \sum_{j =1}^{m'}  b_j^2)\, ,
 \end{equation}  
 then the torsion-free sheaf $E$ is not $L$-very stable.  \\
 (2) For $1 \le t \le m'$, let $b_j= n$ for the nodes $y_j,  j \in I_t$ and $b_j= 0$ for the nodes $y_j,  j \notin I_t$ ($I_t$ defined in Notation 2).  If $d(L) > 2g - 2 + \frac{2}{n-1} (m' - nt)$ then there is no $L$-very stable torsion-free sheaf of local type $\bar{b}$. 
\end{proposition}

\begin{proof}
(1) Let $R_1 \ =  2g - 2m' -2$ and $R_2 = 2g - 2 + \frac{2}{n(n-1)} (m' n - \sum_{j =1}^{m'}  b_j^2)$. 
Then 
$$
\begin{array}{lll}
R_2 - R_1 & =  & \frac{2}{n(n-1)} (m' n - \sum_{j =1}^{m'}  b_j^2) +2m'\\
           {}    & =  &  \frac{2}{n(n-1)} (m' n^2  - \sum_{j =1}^{m'}  b_j^2) \\
            {}    &  =  &   \frac{2}{n(n-1)} \sum_{j =1}^{m'} (n^2  -  b_j^2)\, .
\end{array}
$$
Since $b_j \le n$, one has $R_2 \ge R_1$. Hence the condition \eqref{c2} implies that $d(p^*L) = d(L) > 2g(X) -2$. Then dim $A = \sum_{i=1}^n h^0(X, p^*L^i) = \frac{n(n+1)}{2} d(L) - n(g - m' -1)$.
For a torsion-free sheaf $E$ of local type $b_j(E) = b_j$ at  a node $y_j$, we have   $d((End E)\otimes L) = n^2 d(L) + \sum_j b_j^2$. By our assumption, $V_E = V_{E,1}$. Hence one has  
$${\rm dim} \ V_{E,1} \ = \ h^0((End E)\otimes L) \ge \chi((End E)\otimes L) = n^2 d(L) - n^2(g-1) + \sum_j b_j^2\, .$$
Then 
$${\rm dim} \ V_{E,1} - {\rm dim} \ A \ \ge \frac{n(n-1)}{2} (d(L) - 2(g-1))  - m' n + \sum_j b_j^2\, .$$
Hence  ${\rm dim} \ V_{E,1} > {\rm dim} \ A$ if  $\frac{n(n-1)}{2} (d(L) - 2(g-1))  - m' n + \sum_j b_j^2 >0$ i.e., 
if $d(L) > 2g-2 + \frac{2}{n(n-1)} (m' n - \sum_{j =1}^{m'}  b_j^2)\, .$          
            Thus when the condition \eqref{c2} holds,  ${\rm dim} \ V_{E,1} > {\rm dim} \ {\bf A}$ and every irreducible component of the fibre of $h_V$ has positive dimension. In particular, the fibre of $h_V$ over $0 \in A$, which consists of $(E, \phi)$ with $\phi$ nilpotent, has a point with $\phi \neq 0$. It follows that $E$ admits a nonzero nilpotent Hitchin field with values in $L$, i.e. $E$ is not $L$-very stable. \\
(2)  If $b_j= n$ for the nodes $y_j,  j \in I_t$ and $b_j= 0$ for the nodes $y_j,  j \notin I_t$,  then  $\mathcal{H}_{\bar b}(n,d,L) \subset N(n, d, L)$ by Proposition \ref{ontoA}(1). Hence for $(E, \phi) \in  \mathcal{H}_{\bar b}(n,d,L)$, one has  $V_E = V_{E,1}$.  Therefore, by Part (1), there are $L$-very stable torsion-free sheaves of type $\bar{b}$ if $d(L) > 2g-2 + \frac{2}{n(n-1)} (m' n - t n^2) =   2g-2 + \frac{2}{(n-1)} (m'  - t n)$.        
\end{proof}

\begin{corollary} \label{u3}
Suppose that $d(L) = 2g-2$. Then there are no very stable torsion-free sheaves in $N_Y(n,d,L)$ of local type $(b_j)_j$ such that $\sum_{j= 1}^{m'} b^2_j \ge m' n + \frac{n(n-1)}{2}$. 

The subset consisting of such sheaves has codimension at least $m' n + \frac{n(n-1)}{2}$.
\end{corollary}
\begin{proof}
Recall that $R_2 = 2g - 2 + \frac{2}{n(n-1)} (m' n - \sum_{j =1}^{m'}  b_j^2)$. Hence $R_2 \le 2g - 3$ if  
$\sum_{j= 1}^{m'} b^2_j \ge m' n + \frac{n(n-1)}{2}$. Then by Proposition \ref{u2}, for $d(L) = 2g -2$, there are no very stable torsion-free sheaves of local type $(b_j)_j$ such that $\sum_{j= 1}^{m'} b^2_j \ge m' n + \frac{n(n-1)}{2}$. 

     By Proposition \cite[Proposition 2.7]{Bh0}, the subset consisting of such sheaves (i.e., with 
 $\sum_{j= 1}^{m'} b^2_j \ge m' n + \frac{n(n-1)}{2}$) has codimension at least $m' n + \frac{n(n-1)}{2}$. We note that $m' n + \frac{n(n-1)}{2} \ge 3$ as $m' \ge 1$.
\end{proof}

\begin{remark} \label{bjEn}
         There exist $L$-very stable torsion-free sheaves $E$ on $Y$ which are not locally free.\\         
(1) Assume that $g \ge m'+1$. Then for $d(L) > 2g- 2m' -2 \ge 0$ there are no $L$-very stable torsion-free sheaves $E$ on $Y$ of local type $b_j(E) = n$ for all $j$. If $d(L) = 2g- 2m' -2 \ge 0$ and $p^*L = \omega_X$,  then there exist $L$-very stable torsion-free sheaves $E$ on $Y$ of local type $b_j(E) = n$ for all $j$. \\
(2) There are $L$-very stable non locally free torsion-free sheaves $E$ on $Y$ which are direct images of $L$-very stable vector bundles on partial normalisations of $Y$. They are of local type $b_j(E) = n$ for $j \in I_t$ and $b_j(E) =0$ for $j \notin I_t$, where $I_t$ is defined in Notation 2.
\end{remark}
\begin{proof} 
(1) Consider the case $b_j = n$ for all $n$. In this case, there is an isomorphism from the moduli space of Hitchin bundles $(\tilde{E}, \tilde{\phi})$ of rank $n$, degree $d-nm'$ on the normalisation $X$ and Hitchin sheaves $(E, \phi)$ on $Y$ of local type $b_j(E) = n$ for all $j$. This isomorphism is given by $E = p_* \tilde{E}, \phi = p_* \tilde{\phi}, \tilde{E} = p^*E /torsion, \tilde{\phi} = p^* \phi/torsion$. 
        
        One has $h^0(Y, End(E) \otimes L) = h^0(Y, p_* (End \tilde{E}) \otimes L) = h^0(Y, p_*(End \tilde{E} \otimes p^*L)) = h^0(X, End \tilde{E} \otimes p^*L)$ so that $V_E = V_{\tilde E}$. Note that $\phi$ is nilpotent if and only if $\tilde{\phi}$ is so. It follows that $E$ is $L$-very stable if and only if $\tilde{E}$ is $L$-very stable. 
        
        Also in this case, $R_1 = R_2$ in Proposition \ref{u2} and the condition \eqref{c2} becomes $d(p^*L) > 2g(X) - 2$. Hence there are no $p^*L$-very stable bundles $\tilde{E}$ if \eqref{c2} holds (by Proposition \ref{u1}). This agrees with the the fact that there are no $L$-very stable torsion-free sheaves $E$ of local type $b_j(E)= n$ for all $j$ (by  Proposition \ref{u2}). Moreover, it is known that for $p^*L = \omega_X$, there are $p^*L$-very stable vector bundles. This implies that for $d(L) = 2g -2m' -2$ with $p^*L = \omega_X$, there are $L$-very stable torsion-free sheaves $E$ on $Y$  of local type $b_j(E) = n$ for all $j$. \\   
(2) This can be proved similarly as (1).
\end{proof}

{\bf A simple Example: $m'=1, n=2, g \ge 3$.} \\
In the case $m'=1, n=2$ one has $0 \le b \le 2$.\\
(1) If $b= 0$, there are no $L$-very stable vector bundles for $d(L) > 2g-2$ (by Proposition \ref{u1}).   \\
(2) If $b =2$, there are no $L$-very stable torsion-free sheaves of local type $b=2$ for $d(L) > 2g-4$ (by Proposition \ref{u2}). For $d(L)= 2g - 4, p^*L= \omega_X$, there are $L$-very stable torsion-free sheaves of local type $b=2$. The sheaves of local type $2$ form a subset of codimension at least $4$.

\subsection{Properness of $h_V$} \hfill

We now investigate properness of $h_V$.   

\begin{proposition} \label{PaPeProp2.2}
(1)  The morphism $i_{V_{E,1}}: V_{E,1}  \to N_Y(n,d,L)$ is proper if and only if $h_V$ is a proper morphism.
In particular, if $V_{E,1}$ is a closed subset of $N_Y(n,d,L)$, then $h_V$ is a proper morphism.\\
(2) If $E$ is a stable vector bundle and $h_V$ is proper, then it is quasi-finite.   
\end{proposition}
\begin{proof}
(1)    The proof is exactly the same as that of \cite[Proposition 2.2]{PaPe} which uses only properties of proper maps \cite[Coro.II4.8(a),(b) and (e)]{Ha}. \\
(2)   If $E$ is a vector bundle, then $h_V: V_E  \to {\bf A'}$. If $h_V$ is proper, then for $a \in {\bf A'}$, the fibre $h_V^{-1} (a)$ is proper over Spec $k$ (by base change of proper morphism). Then $h_V^{-1}(a)$ is a subset of an affine variety $V_E$. Hence  $h_V^{-1} (a)$ must be a finite set i.e., $h_V$ is quasi-finite. 

          We remark that this means that dim $V_E \le $ dim ${\bf A'}$.

\end{proof}

\begin{proposition} \label{PaPeProp2.1}
Let $E$ be a stable vector bundle.\\
(1)  If $h_V$ is proper, then it is finite. \\
(2)  If $h_V$ is quasi-finite, then $E$ is $L$-very stable, i.e. $E$ admits no non-zero nilpotent Hitchin field $\phi$ with values in $L$. \\ 
\end{proposition}
\begin{proof}
(1) This follows from Proposition \ref{PaPeProp2.2} as a proper map between affine varieties is finite \cite{Ha}.

(2) If the vector bundle $E$ admits a non-zero nilpotent Hitchin field $\phi$ with values in $L$, then $h_V(E,\phi) =0, \phi \neq 0$. Since $E$ is stable, $\mathbb{C}^*\phi \cong \mathbb{C}^*$. Hence the nilpotent cone $h_V^{-1}(0)$ will contain the line $\mathbb{C}^*\phi$, contradicting quasi-finiteness of $h_V$. \\
\end{proof}

       Henceforth we assume that $E$ is a vector bundle.  Let $\bar{V_E}$ denote the closure of $V_E$ in $\mathcal{H}_Y(n,d,L)$ and  $\bar{V}_{E,1} = \bar{V_E} \cap N_Y(n,d,L)$. The latter is the closure of $V_E$ in $N_Y(n,d,L)$ and hence a closed subset in $N_Y(n,d,L)$. As $\bar{V}$ and $N_Y(n,d,L)$ are $\mathbb{C}^*$- invariant, so is $\bar{V}_1$.
Putting together Propositions \ref{PaPeProp2.1} and \ref{PaPeProp2.2}, we see that if $E$ is a stable vector bundle such that $V_E$ is closed in $N_Y(n,d,L)$, then $E$ is very stable i.e. $E$ admits no non-zero nilpotent Hitchin field $\phi$ with values in $L$.

\begin{proposition} \label{PaPeProp2.3} 
Let $E$ be a stable vector bundle on $Y$ and $(F,\psi) \in \bar{V}_E \setminus V_E$.
Then $F$ is not a semistable torsion-free sheaf.
\end{proposition}

\begin{proof}
If $(F,\psi) \in \bar{V}_E \setminus V_E$, then there exists a discrete valuation ring $R$ with its fraction field $K$ and a morphism $i_F: Spec R \to N_Y(n,d,L)$ such that the generic point $Spec K$ maps into $V_E$ and the closed point maps to $(F, \psi)$ \cite{StackP}. If $F$ is semistable, composing $i_F$ with the forgetful morphism $fr: N_Y(n,d,L) \to U(n,d)$, we get a morphism which maps $Spec K$ to $E$ and the closed point of $Spec R$ to $F$ contradicting the separatedness of the variety $U(n,d)$. This proof is essentially Lemma 3.11 and Remark 3.12  of \cite{Pe}. 
\end{proof}

\begin{theorem} \label{Thm1.1PaPe}
    Denote by $E$ a stable vector bundle of rank $n$ and degree $d$ on $Y$. Then the following statements are equivalent: \\
(1)  $E$ is very stable.\\
(2)  $V_E$  is a closed subset of $N_Y(n,d,L)$. \\ 
(3)  $h_V$ is a proper morphism.\\
(4)  $h_V$ is quasi-finite.    
\end{theorem}

\begin{proof}
      The proof is on the same lines as that of \cite[Theorem 1.1]{PaPe}, so we omit some details. 
      In view of Propositions \ref{PaPeProp2.1} and \ref{PaPeProp2.2}, we only need to show that if $E$ is very stable, then $V_E$ is closed in $N_Y(n,d,L)$. We shall prove the equivalent statement that if $V_E$ is not closed in $N_Y(n,d,L)$, then $E$ admits a non-zero nilpotent Higgs field $\phi_0$. Assume that $V_E$ is not closed, let 
$(F,\psi) \in \bar{V}_{E,1} - V_E$. \\
{\bf Step 1:} Replacing $(F, \psi)$ by its limit  under $\mathbb{C}^*$-action we may assume that $h(F,\psi)= 0$ and $(F, \psi)$ is a fixed point for $\mathbb{C}^*$-action. 

   This can be  proved exactly as the corresponding statement on \cite[p. 146]{PaPe} (replacing $\bar{V}$ by 
   $\bar{V}_{E,1}$ in the proof in \cite{PaPe}). We note that  we use \ref{PaPeProp2.3} and \cite[Lemma 2.4]{PaPe}. \\ 
{\bf Step 2:}    There exists a smooth curve $C$, a point $c \in C$ and 
a morphism $\Phi: C \to \bar{V}_{E,1}$ with $\Phi(C^*= C \setminus \{c\}) \subset V$ and $\Phi(c) = (F,\psi)$ (by  \cite[Lemma 2.4]{PaPe}).      
Define a morphism
$$\Psi^*:  \mathbb{C}^* \times {C}^* \to V_E \subset \bar{V}_{E,1}$$
by $\Psi(\lambda,p) = \lambda \Phi(p).$

It is easy to see that $\Psi^*$ extends to a morphism $\Psi: \mathbb{C} \times {C}^* \to \bar{V}_{E,1}$  defined by $\Psi \vert_{\mathbb{C}^* \times {C}^*} = \Psi^*$ and $\Psi(0, p) = 0, p\in {C}^*$ and to a morphism $\Psi: \mathbb{C}^* \times {C} \to V_E$ defined by $\Psi(\lambda, c)= \lambda \Phi(c)= (F,\psi)$ (as $(F, \psi)$ is fixed by $\mathbb{C}^*$). 
        Thus 
        $$\Psi: S= \mathbb{C} \times C \dasharrow \bar{V}_{E,1}$$  
 is a rational map with indeterminacy locus the point $(0, c)$. We want to prove that it can be resolved into a morphism (by a finite number of blow ups at points).\\
{\bf Step 3:}  Let $\Gamma$ be the closure of the graph of $\Psi$ in $S \times N_Y(n,d,L)$. Let $\pi_1$ be the projection $S \times N_Y(n,d,L) \to S, \pi = \pi_1\vert_{\Gamma}$. By \cite[Theorem II and discussion on p.140]{Hi}, to resolve the map $\Psi$ it suffices to show that $\pi$ is proper. 

          Let $R$ be a discrete valuation ring with quotient field $K$. By the valuative criterion of properness, to show that $\pi$ is proper we need to show that in any commutative diagram as follows, there is a morphism $Spec R \to \Gamma$ making the resulting diagram commutative. 
$$
\begin{array}{ccc}
    Spec(K) & \longrightarrow & \Gamma \\
    {} & {} & {}\\
    \downarrow & {} & \pi \downarrow\\
    {} & {} & {}\\
    Spec (R) & \longrightarrow & S\, . 
\end{array}
$$           
     Define a rational map $h': h \circ \Psi: S \dasharrow A'$. This map is defined by holomorphic functions on $S \setminus (0,c)$ (as $A'$ is a vector space). These functions extend to $S$ (Hartog's theorem) and $h'(0,c)=0$ (by continuity).

To prove the properness of $\pi$ we use the properness of $h$.  We consider the extended commutative diagram
$$
\begin{array}{ccccc}
    Spec(K) & \longrightarrow & \Gamma & \hookrightarrow & S \times N_Y(n,d,L)\\
    {} & {} & {} & {} & {}\\
    \downarrow & \nearrow & \pi \downarrow & {} &  id\times h \downarrow \\
    {} & {} & {} & {} &{} \\
    Spec (R) & \longrightarrow & S & \stackrel{id\times h'}{\longrightarrow} & S \times A'  . 
\end{array}
$$           
Since $id \times h$ is proper and $\Gamma$ is closed, $(id\times h)\vert_{\Gamma}: \Gamma \to S\times A'$ is proper. Hence by valuative criterion, the composite $Spec (R) \to S \to S \times A'$ lifts to a morphism $e_1: Spec (R) \to \Gamma$ making the diagram commutative.\\
{\bf Step 4}: In view of step 3, the rational map $\Psi$ resolves to a morphism 
$$\hat{\Psi}: \hat{S} \rightarrow \bar{V}_{E,1} \subset  \mathcal{H}_Y(n,d,L)$$
after a finite sequence of blow ups along points \cite[Theorem II and discussion on p.140]{Hi}.  
The exceptional divisor $D$ is a connected union of projective lines $D_i$ and $\Psi$ maps $D$ to a connected curve in $\bar{V}_{E,1}$. Let $p_0 = \stackrel{lim}{p \to c} (0,p), \ p_{\infty}= \stackrel{lim}{\lambda\to 0}(\lambda,c)$. Then $\hat{\Psi}(p_0)= (E,0), \ \hat{\Psi}p_{\infty} = (F, \psi)$ by the separatedness of $\mathcal{H}_Y(n,d,L)$ (the separatedness of $\mathcal{H}_Y(n,d,L)$ follows from the construction of the moduli space as in the case $Y$ is smooth). Moreover, since $h'(0,c) = 0$, one has $h(\phi) = 0$ for $\phi \in \hat{\Psi}(D)$. We can arrange the projective lines $D_i$ in a sequence $D_0, D_1, \cdots , D_{m'}, m' \le m$ such that $p_0 \in D_0, p_{\infty} \in D_{m'}$ and $D_i \cap D_{i+1} \neq \empty$ for $i \le m' -1$.  Let $i_0\ge 0$ be the smallest integer such that 
$\hat{\Psi}(D_{i_0})\neq  \{(E,0)\}$. Then one can show that $\hat{\Psi}(D_{i_0}) \cap V$ contains a point corresponding to a Higgs bundle $(E, \phi_0)$ with a non-trivial Higgs field $\phi_0$, then $h(E, \phi_0) =0$ so that $\phi_0$ is nilpotent.
        
\end{proof}

\begin{remark}    \label{rmk} 

Suppose that $E$ is a $L$-very stable vector bundle.   \\          
(1)   For $L= \omega_Y$, one has dim $V_E = n^2(g-1) +1$, dim $A' = n^2(g-1)+1$ so that if $h_V$ is (quasi-) finite, then $h_V$ is surjective and $V_E$ intersects every Hitchin fibre in finitely many points. Since $h_V^{-1}(0)$ consists of nilpotent Hitchin fields and $E$ is very stable, from Theorem \ref{Thm1.1PaPe} it follows that $h_V^{-1}(0) = 0$.\\
(2)   Let $d(L) = 2g -2, L^n \neq \omega_Y^n$. Then for all $i \ge 1, h^0(Y, L^i)= (2i-1)(g-1)$, so that 
$$\sum_{i=1}^n h^0(Y, L^i)= \sum_{i=1}^n (2i-1)(g-1) = (g-1)n(n+1) - n(g-1) = n^2(g-1)\, .$$  
 Note that $E \cong E\otimes L^*\otimes \omega_Y$ implies that $det E \cong det E \otimes (L^*\otimes \omega_Y)^n$ so that $L^n \cong \omega_Y^n$. Hence under our assumption, $h^1((End E) \otimes L) = $ Hom $(E \ \to \ E\otimes L^*\otimes \omega_Y) = 0$. Then one has dim $V_E = h^0(End E \otimes L) = \chi(End E \otimes L) = n^2(g-1)$. Thus dim $V = $ dim ${\bf A'}$ so the conclusions of Part (1) hold (with the same argument). \\
(3) In case dim $V_E = $ dim ${\bf A'}$, one can give a shorter proof of Theorem \ref{Thm1.1PaPe} using  \cite[Lemma 3.1]{Pe}.   
\end{remark}

\section{Parabolic Hitchin sheaves and Parabolic GPH}

   Let $Y$ be an integral projective algebraic curve. 
Let $I$ be a finite set of smooth points of $Y$. Let $D = \sum_{x\in I} x$ be the Cartier divisor in $Y$.  For a torsion-free sheaf $E$ on $Y$, let $r(E)$ and $d(E)$ denote  the rank and degree of $E$ respectively.
 Let $L_0$ be a (fixed) line bundle on $Y$ with $d(L_0(D)) \ge  0$.

\subsection{Parabolic Hitchin Sheaves}\label{parasheaves}

\begin{definition} \label{parastructure}
A \textit{quasi-parabolic structure} at a point $x \in I$ on a torsion-free sheaf $E$ on $Y$ is a flag of vector subspaces 
on the fibre $E_x$ 
$$E_x = F'_1(E_x) \supset \cdots \supset F'_{l_x}(E_x) \supset F'_{{l_x}+1}(E_x) = 0\, .$$
Let $c_x^i = dim \ F'_i(E_x)$, then $\bar{n}_x= (c^1_x, \cdots, c^{l_x}_x)$ is called the type of the flag at $x$.        

A \textit{parabolic structure} on $E$ over $x \in I$ is a
quasiparabolic structure at $x \in I$ as above
together with a rational number $\alpha_i(x)$ (called a parabolic weight) for each subspace $F'_i(E_x)$ such that
$$0 < \alpha_1(x) < \cdots < \alpha_{l_x}(x) < 1$$

Let $n_i(x) = dim \ F'_{i}(E_x) - dim \ F'_{i+1}(E_x), i= 1, \cdots, l_x$ and $r_i(x) = \sum_{j=1}^{i} n_j(x),  i \ge 1\, .$ The number $n_i(x)$ is called the multiplicity of $\alpha_i(x)$.
 Define 
$$\bar{\alpha}(x) := (\alpha_1(x), \alpha_2(x), \cdots \alpha_{l_x}(x));
~ \bar{n}(x) := (n_1(x), n_2(x), \cdots , n_{l_x}(x))\, .$$

                       A ($L$-twisted) Quasi-parabolic (respectively parabolic) Hitchin sheaf $(E_*, \phi)$ on Y with quasi-parabolic (respectively parabolic) structure over $I$ is a torsion-free sheaf $E$ on $Y$ together with a quasi-parabolic (respectively parabolic) structure at each $x \in I$ and a homomorphism $\phi: E_* \to E_* \otimes L(D)$ of parabolic sheaves i.e.,  $\phi_x (F'_i(E_x)) \subset F'_i(E_x) \otimes L(D)_x$ for all $i$ and $x \in I$.  
  A quasi-parabolic (respectively parabolic) Hitchin sheaf $(E_*, \phi)$ is called a quasi-parabolic (respectively parabolic) Hitchin bundle if the sheaf $E$ is locally free.

                 A strongly quasi-parabolic (respectively parabolic) Hitchin sheaf $(E_*, \phi)$ on Y with quasi-parabolic (respectively parabolic) structure over $I$ is a torsion-free sheaf $E$ on $Y$ together with a quasi-parabolic (respectively parabolic) structure at each $x \in I$ and a strongly quasi-parabolic (respectively parabolic) homomorphism $\phi: E_* \to E_* \otimes L(D)$ of quasi-parabolic (respectively parabolic) sheaves i.e., $\phi _x (F'_i(E_x)) \subset F'_{i+1}(E_x) \otimes L(D)_x$ for all $i$ and $x \in I$.  
  A quasi-parabolic (respectively parabolic) Hitchin sheaf $(E_*, \phi)$ is called a quasi-parabolic (respectively parabolic) Hitchin bundle if the sheaf $E$ is locally free.

\end{definition} 

For a parabolic Hitchin sheaf $(E_*, \phi)$ as above, the parabolic degree is defined as
$$
{par-deg}(E_*, \phi)\, :=\, d(E)+\sum_{x \in I}
\sum_{i=1}^{l_x} \alpha_i (x) n_i(x)\, .
$$
It is independent of $\phi$. The parabolic slope is defined as
$$
{par-}\mu(E_*)\, :=\,
\frac{par-deg(E_*)}{ r(E)}\, \in\, {\mathbb Q}\, .
$$

\begin{definition}
For a parabolic Hitchin sheaf $E_*$, any nonzero subsheaf $F\, \subset\, E$
has an induced parabolic structure. We denote by $F_*$ the sheaf $F$ with the induced parabolic structure.  
A subsheaf $F \subset E$ is called $\phi$-invariant if $\phi(F) \subset F \otimes L(D)$. Then $(F_*, \phi\vert_F)$ is a $\phi$-invariant subsheaf of $(E_*, \phi)$.

               A parabolic Hitchin sheaf $(E_*, \phi)$ is called \textit{stable} (respectively, \textit{semistable}) if 
for every $\phi$-invariant subsheaf $F \subset E$, with $1\, \leq\, r(F)
\, <\, r(E)$ and with a torsion-free quotient, one has
$$
{par-}\mu(F_*)\, <\, {par-}\mu(E_*)
$$
(respectively, ${par-}\mu(F_*)\,\leq\, {par-}\mu(E_*)$)\, .

           The (semi)stability of a strongly parabolic Hitchin sheaf is defined similarly.

\end{definition}

\subsection{Parabolic GPH on $Y$}   \hfill

 For $j=1,\cdots, m'$, let $\{ D_j \}$ be disjoint (Cartier) divisors on $Y$ with supports contained in the set of regular points of $Y$.

\begin{definition}\label{paragphdefi}    
A ($L_0$-twisted) parabolic GPH of rank $n$ and degree $d$ on $Y$ is a triple $(E_*, F(E), \phi)$ 
where \\
(1) $E_*$ is a parabolic vector bundle of rank $n$ and degree $d$ on $Y$,\\
(2) $F(E) = (F_1(E), \cdots, F_{m'}(E))$ is a  generalised parabolic structure on $E$ over divisors $D_j$  as defined in Definition \ref{gphdefi}(2).\\
(3) $\phi: E \rightarrow E \otimes L_0(D)$ is a homomorphism which preserves the parabolic structure. \\

           A parabolic GPH is called a strongly parabolic GPH if $\phi$ is a strongly parabolic homomorphism of GPH.
\end{definition}

For a subbundle $N \subset E$, let $F_j(N)= F_j(E)\cap (N \otimes {\mathcal O}_{D_j})$ and 
$f_j(N)=$ dim $F_j(N)$. 

$N$ is called $\phi$-invariant if $\phi(N) \subset N \otimes L_0(D)$. 

Let $0 < \alpha \le 1$ be a real number. Let $wt(N)$ and $pd(N)$ denote the parabolic weight and parabolic degree  of $N$ i.e.,  $pd(N) := par - deg(N_*, \phi\vert_{N_*})$. Define 
 $$gwt(N) = \alpha \sum_j f_j(N)\, , \ gpd(N) = pd(N) + gwt(N) = d(N) + wt(N) + gwt(N)\, , $$
         $$ gp\mu(N) \ = \ gpd(N)/r(N)\, .$$

\begin{definition}  \label{parGPHstability}
A  parabolic GPH $(E_*, F(E), \phi)$ is (semi)stable if for every proper 
$\phi$-invariant subsheaf $N \subset E$ (with a torsion-free quotient) with induced parabolic and GPB structure, one has 
$$\frac{pd(N)+ \alpha \sum_j f_j(N)}{r(N)} \ \ (\le) \ \ \frac{pd(E)+ \alpha \sum_j f_j(E)}{r(E)}\, ,$$ 
i.e. $gp\mu(N) \le gp\mu(E)$ for semistability and 
$gp\mu(N) < gp\mu(E)$ for stability.
\end{definition}

     We recall that the pair $(E, F(E))$ (together with $\alpha$) is called a generalized parabolic bundle, a GPB in short. By forgetting the Hitchin field $\phi$ in a GPH, we get a generalized parabolic bundle (see \cite{B1} for generalities on GPBs). 

\subsection{Good parabolic GPH} \label{goodparGPH} \hfill
                
       Henceforth, $Y$ is an integral projective curve  and          
       $$p : X  \rightarrow Y$$ 
 the normalisation map. Let ${\mathcal C}= \ {\rm Ann} \  (p_*{\mathcal O}_X/{\mathcal O}_Y)$ denote the sheaf of conductors of $Y$. Let $y_1, \cdots , y_{m'}$ be the singular points of $Y$. Let ${\tilde{\mathcal O}}_{y_j}$ denote the normalisation of the local ring ${\mathcal O}_{y_j}$ at $y_j$ and 
$$\delta(y_j)= \ {\rm dim}_k \ ({\tilde{\mathcal O}}_{y_j}
/{\mathcal O}_{y_j})\, , \ \delta(Y)= \sum_{j=1}^{m'} \delta (y_j)\, .$$  
For $j=1,\cdots, m'$, let $\{ D_j \}$ be disjoint (Cartier) divisors on $X$ with support of $D_j$ contained in $p^{-1} y_j$ and $p_* {\mathcal O}_X(-D_j) \subset {\mathcal C}$.  Choose 
$$f_j(E)= r(E)(d(D_j) - \delta(y_j))\, .$$ 
Let $I \subset X$ be a set of finitely many distinct points disjoint from $p^{-1}y_j$ for all $j$. 
We fix a line bundle $L$ on $Y$ and define $L_0:= p^*L_0$.
 
Note that $\phi$ induces (on taking $p_*$ and using the projection formula) 
$$p_*\phi \colon p_*E \rightarrow p_*(E\otimes p^*L)\cong p_*E\otimes L\, .$$     

\begin{definition}
A GPH $(E, F(E), \phi)$ is a good GPH if it satisfies the following conditions. 
\begin{equation}\label{c3}
{\rm The} \ k-{\rm subspace} \  p_* F_j(E)  \ {\rm is \ an} \ {\mathcal O}_{y_j}-{\rm submodule \ of} \  p_* (E\otimes{\mathcal O}_{D_j})\, .
\end{equation}
\begin{equation}\label{c4}
(p_* \phi)_{y_j} (p_*(F_j(E))\subset (p_*(F_j(E))\otimes L_{y_j}\, . 
\end{equation}
\end{definition}

\begin{proposition}\label{bdd1}
Let $(E_*, F(E),\phi)$ be an $\alpha$-semistable $L_0$-twisted parabolic GPH. \\
(1) If $\phi \neq 0$, then $d(L_0(D)) \ge 0$ and for any subsheaf $F \subset E$ with a torsion-free quotient we have 
$$\mu(F) \le \ {\rm max} \ \{\mu(E)+ \frac{1}{n}\ (gpwt(E)), \  
\mu(E)+  \frac{1}{n} gpwt(E)+ \frac{(n-1)^2}{n}d(L_0(D)) \}\, .$$
(2) If $\phi \neq 0$, then 
$$ gp\mu(E)-\frac{(n-1)^2}{n}d(L_0(D)) \le gp\mu((E_i)_*/ (E_{i-1})_*)
\le gp\mu(E)+\frac{(n-1)^2}{n}d(L_0(D))\, ,$$ 
where $0 = (\underline{E}_0)_* \subset  (\underline{E}_1)_* \subset (\underline{E}_{\ell})_* = E_*$
is the Harder-Narasimhan filtration for the parabolic GPB $\underline{E}_* = (E_*, F(E))$. 
 In particular, if $\phi \neq 0$ and $d(L_0(D))=0$, then $(E_*,F(E))$ is an  $\alpha$-semistable GPB.\\
(3) If $\phi =0$, then $(E_*,F(E))$ is an $\alpha$-semistable parabolic GPB.
\end{proposition}     

\begin{proof}
 (1)      This can be proved on similar lines as \cite[Proposition 3.1(1)]{Bh1}. 
  In Case (1), we assume that the underlying parabolic GPB $(E_*,F(E))$ is $\alpha$-semistable. Then the argument similar to that in the proof of \cite[Proposition 3.1(1)]{Bh1} (with $p\mu$ replaced by $gp\mu$ etc.) gives $\mu(F_*)= p\mu(F) \le p\mu(E) + \frac{1}{n} gwt(E)$. Since $\mu(F) \le p\mu(F), p\mu(E) = \mu(E) + \frac{1}{n} wt(E)$, one has $\mu(F) \le \mu(E)+ \frac{1}{n}\ ( wt (E) + gwt(E))$.  Similarly Case 2 works out as the Case 2 in the proof of \cite[Proposition 3.1]{Bh1}. \\
 (2)  The part (2) also works similarly as in the proof of \cite[Proposition 3.1(2)]{Bh1} .\\
 (3)  This is clear.
\end{proof}

\begin{corollary}\label{bdd2}
Let $(E_*,F(E),\phi)$ be an $\alpha$-semistable $L_0$-twisted parabolic GPH.
\begin{enumerate}
	\item  $H^1(Y, E)=0$ if 
$$d(E)> \ {\rm max} \ \{n(2g-2)+  (n-1) gpwt(E), \  
n(2g-2)+ (n-1) gpwt(E) + \frac{(n-1)^3}{n} d(L_0(D)) \}\, .$$  
	\item $E$ is globally generated if 
$$d(E)> \ {\rm max} \ \{n(2g-1)+  (n-1) gpwt(E), \  
n(2g-1)+ (n-1) gpwt(E) + \frac{(n-1)^3}{n} d(L_0(D)) \}\, .$$  	
	\item The evaluation map $e_{m'}: H^0(X,E) \to 
	\oplus_j (E\otimes {\mathcal O}_{D_j})$ is surjective if 
$$\mu(E) > \ {\rm max} \ \{2g-2+\sum_j d(D_j)+ \frac{1}{n} gpwt(E), $$ 
$$2g-2+\sum_j d(D_j)+  \frac{(n-1)}{n}  gpwt(E) + \frac{(n-1)^3}{n^2} d(L_0(D)) \}\, .$$	
\end{enumerate}
\end{corollary}    

\begin{proof}  
(1)     Assume that $h^1(Y, E) = h^0(Y, E^*\otimes \omega_Y) \neq 0$. Let $F$ be the kernel of the nonzero homomorphism $E \to \omega_Y$. Then $d(F) \ge d(E) - (2g - 2), r(F) = n-1$ i. e., 
\begin{equation}\label{be1}
     \mu(F) \ge \frac{n}{n-1}  \mu(E) -  \frac{2g -2}{n-1}.
     \end{equation}
 By Proposition \ref{bdd1},  
     $\mu(F) \le \ {\rm max} \ \{\mu(E)+ \frac{1}{n}\  gpwt(E), \  
\mu(E)+  \frac{1}{n} ( gpwt(E)+ \frac{(n-1)^2}{n}d(L_0(D))) \}\, .$ 
Suppose that the maximum is  $\mu(E)+ \frac{1}{n} (gpwt(E) + \frac{(n-1)^2}{n}d(L_0(D)))$.  Then from the inequality \eqref{be1} we get
$\mu(E) + \frac{1}{n} gpwt(E) + \frac{(n-1)^2}{n}d(L_0(D)) \ge \frac{n}{n-1}  \mu(E) -  \frac{2g -2}{n-1} $ i. e., $\mu(E) \le \frac{n-1}{n}(gpwt(E) + \frac{(n-1)^2}{n}d(L_0(D)))+ 2g - 2$. Hence if $d(E) > (n-1)(gpwt(E) + \frac{(n-1)^2}{n} d(L_0(D))) + n(2g -2)$, then $h^1(Y, E) = 0$. Similarly, if the maximum of $\mu(F)$ is $\mu(E)+ \frac{1}{n}\  gpwt(E)$, then we get the condition $d(E) >  (n-1) gpwt(E) + n(2g -2)$.\\
(2)   The vector bundle $E$ is globally generated if for all $x \in Y$, $H^1(Y, E(-x))= 0$.  Hence the condition for global generation follows from Part (1) by replacing $E$ by $E(-x)$.  \\   
(3) The evaluation map $e_{m'}$ is surjection if $H^1(Y, E(-\sum_j D_j)) = 0$.  Therefore  the condition for global generation follows from Part (1) by replacing $E$ by $E(-\sum_j D_j))$.           
\end{proof}

\subsection{The moduli space for parabolic GPH} \label{moduliparGPH} \hfill

            Generalising the constructions in \cite{Bh1} to the parabolic case, one can construct the moduli space 
$\mathcal{H}_{Y, par}(n,d)$ (respectively $\mathcal{H}_{par}(n,d, \chi)$) of parabolic Hitchin sheaves of rank $n$, degree $d$ (respectively with a "fixed determinant" $\chi$.) on $Y$ and the moduli space $H_{par}(n,d,L_0)$ (respectively $M_{par}(n,d,L_0)$) of parabolic GPH (respectively good parabolic GPH). We briefly sketch the construction in the following.             
               
             We note that it follows from  Proposition \ref{bdd1} that semistable parabolic GPH of rank $n$, degree $d$ with a fixed parabolic and GPB structure form a bounded family.   
We choose an ample line bundle ${\mathcal O}_Y(1)$ of degree $1$ on the irreducible curve $Y$ of genus $g$. 
Fix Cartier divisors $D_j \subset Y_{regular}\, , j=1, \cdots, m'$. Let $I \subset Y$ be a finite subset of distinct closed nonsingular points of $Y$ disjoint from all $D_j$.  
Let $P(t)$ denote the Hilbert polynomial of a vector bundle $E$ of rank $n\ge 2$, degree $d$ on $Y$ i.e., $P(t)= nt +d+n(1-g).$
 
Fix a rational number $\alpha\, , 0<\alpha \le 1$ and a line bundle $L_0$ on $Y$. 

     Let $S$ be the set of $\alpha$-semistable $L_0$-twisted GPH $(E,F(E),\phi)$ of fixed rank $n$, degree $d$, Euler characteristic $N= d+n(1-g)$, with GPH structure given by flags of type 
$E\otimes {\mathcal O}_{D_j} \supset F_j(E) \supset 0$ with $f_j(E)= n(d(D_j) - \delta(y_j))\, ,$  weights $(0,\alpha)$ and a homomorphism $\phi: E \to E \otimes L_0(D)$ and parabolic structure of type $\bar{n}_x$ for every $x \in I$. Let $q_j=q_j(E)$ for any $(E_*,F(E),\phi) \in S$. 
Tensoring $E$ by a line bundle, we may assume that $d$ is large enough so that there exists an integer $N_1$ such that for all $N \ge N_1$ and for all vector bundles $E$ underlying the elements of $S$, we have $H^1(Y,E(- \sum_j D_j - x)) = 0, x \in Y$ so that $E(- \sum D_j - x)$ is globally generated, $H^0(Y, E) \cong k^N$ and the evalutation map $e\colon  H^0(Y,E) \to \bigoplus_{j=1}^{m'} (E\otimes {\mathcal O}_{D_j})$ is surjective (by Proposition \ref{bdd2}). 
Fix $N \ge N_1$ and let 
                               $$V= k^{N}\, .$$  
Denote by $Q$ the quot scheme 
$$Q:= Quot(V \otimes_k {\mathcal O}_Y, P(t)),$$
the scheme of quotients of $V \otimes {\mathcal O}_Y$ with Hilbert polynomial $P(t)= nt+N$. Let 
$$q_Q: V \otimes {\mathcal O}_{Y\times Q} \rightarrow {\mathcal E}$$
be the universal quotient sheaf over $Q\times Y.$ For $q\in Q$, let 
$\E_q= \E\mid_{q\times Y}$.  

Let $Q_{lf} \subset Q$ be the subscheme 
$$ 
Q_{lf} := \{ q\in Q \ \mid \ \E_q \ {\rm is \ locally \ free}, \ 
H^0(V\otimes {\mathcal O}_{Y}) \to H^0({\mathcal E}_q) \ {\rm is \  isomorphism} \}\, .$$ 
Note that $H^1({\mathcal E}_q)= 0$ for $q \in Q_{lf}$.  
Over $Q_{lf}$, we have the locally free sheaf 
$${\mathcal F}= (End \ \E)  \otimes p_Y^*L_0(D)\, .$$ 

There exists a scheme (constructed using ${\mathcal F}$)  
$${\tilde Q} \rightarrow Q_{lf}$$
together with a local universal family of Hitchin pairs $({\mathcal E}, \phi)$ over ${\tilde Q}\times Y$ (\cite[Lemma 3.5]{N}), where $\mathcal{E}$ denotes the pull back of $\mathcal{E}\mid_{Q_{lf}\times Y}$ to ${\tilde Q}\times Y$.   

                      For $x\in I$, let $Flag_{\bar{n}(x)} (\mathcal{E}_x)$ denote the relative flag scheme of type $\bar{n}_x$ over ${\tilde Q}$ whose fibre over the point corresponding to $(E, \phi)$ is the flag variety  of type $\bar{n}_x$ on the fibre $E_x$. Let
$$Flag_{\tilde Q} := \times_{x \in I}  \ Flag_{\bar{n}(x)} (\mathcal{E}_x)\, ,$$
a fibre product over ${\tilde Q}$.  Let $p_{Fl}: \mathcal{E} \to Flag_{\tilde Q} \times Y$ be the pull back of $\mathcal{E}$ on ${\tilde Q}\times Y$ to $Flag_{\tilde Q} \times Y$.

            For each $j$, there is vector bundle ${\mathcal E}_{D_j}= (p_{Fl})_* ({\mathcal E}\otimes p_Y^*{\mathcal O}_{D_j})$ of rank $n d(D_j)$ over $Flag_{\tilde Q}$. Let ${\mathcal G}r_j$ be the Grassmannian bundle over $Flag_{\tilde Q}$ of $q_j$-dimensional quotients of ${\mathcal E}_{D_j}$ and 
$${\mathcal G}r = \times_j {\mathcal G}r_j,$$
the fibre product of ${\mathcal G}r_j$ over $Flag_{\tilde Q}$. It is a Grassmannian bundle over $Flag_{\tilde Q}$ with fibres isomorphic to the product of Grassmannians $\prod_j Gr(n d(D_j), q_j)$. A point of ${\mathcal G}r$ corresponds to a parabolic GPH $(E_*, F(E), \phi)$.

{\bf The scheme $Z$ and the morphism $\Psi \colon  Flag_{\mathcal{G}r} \to Z$} \hfill
        
      Let 
     $$W := H^0({\mathcal O}_Y (m))\, .$$   
 There exists an integer $M_1(N)$ such that for all $m \ge M_1(N)$ there is a closed embedding 
$$P \colon Q \hookrightarrow Gr(V\otimes W, P(m))\, .$$ 

      We choose $m \ge M_2(N)$ such that $(S^*_n(L_0^*(-D)))(m)$ is generated by global sections and $H^1((S^*_n(L_0^*(-D)))(m))=0$. Let 
      $$W_1:= H^0(Y, (S^*_n(L_0^*(-D)))(m) )\, .$$    
The Hitchin field $\tilde{\Phi}$ on $\E$ gives a homomorphism 
$$\Phi' \colon \E \otimes p_Y^*(S^*_n(L_0^*(-D))) \to \E$$
 over $\tilde{Q}\times Y$. Composing $q_{\tilde Q}\otimes 1$ with $\Phi'$ and tensoring  with $p_Y^*({\mathcal O}_Y(m))$ gives the homomorphism 
$$\Phi_m \colon V\otimes p_Y^*((S^*_n(L_0^*(-D)))(m)) \longrightarrow \E(m)\, .$$ 
Taking the direct image of $\Phi_m$ on ${\tilde Q}$, we have  
$$\Phi_V \colon (V\otimes W_1)\otimes {\mathcal O}_{\tilde Q} 
\rightarrow {p_{\tilde Q}}_* \E(m)\, .$$

One can show that $\Phi_V$ defines an embedding 
$${\tilde P}: 
{\tilde Q} \longrightarrow Gr(V\otimes W_1, P(m))\, ,$$
which lies over the embedding of $Q$ (see \cite[Subsection 3.3]{Bh2} for details).

     Let $V_{Flag \times Y} = V\otimes_k {\mathcal O}_{Flag_{\tilde{Q} \times Y}}$. 
We have the universal quotient 
$$q_{Flag}\otimes 1 \colon V_{Flag \times Y}\otimes p_Y^*{\mathcal O}_Y(m) \rightarrow \ \mathcal{E}(m)\, ,$$ 
where $\mathcal{E}(m)$ denotes the pull-back of $\mathcal{E}(m)$ (on $\tilde{Q}\times Y$) to 
$Flag_{\tilde{Q} \times Y}$. 
 
            For each $x \in I$,  tensoring $q_{Flag}\otimes 1$ with $p_Y^*\mathcal{O}_x$, one gets a surjective morphism $V_{Flag_{\tilde{Q}} \times x}  \rightarrow (\tilde{E}(m))_x$. By the universal property of the 
relative flag scheme $Flag_{\bar{n}_x} (\mathcal{E}_x)$ , this gives a morphism 
$P_{f,x}:  Flag_{\tilde Q} \rightarrow    Gr_x$ where 
$Gr_x := Gr_n(k^N) \times Gr_{r_1(x)}(k^N)\times \cdots \times Gr_{r_{l_x}}(k^N)$ (note $V= k^N$). Let $P_f := \times_{x\in I} P_{f,x}$.            

   Let $q_{\mathcal{G}r}$ be the pull back of $q_{Flag}$ to $\mathcal{G}r$. 
 Tensoring $q_{\mathcal{G}r}\otimes 1$ with $p_Y^*{\mathcal O}_{D_j}$ one gets a surjective map $V_{\mathcal{G}r \times Y}\otimes p_Y^*{\mathcal O}_{D_j} \ \rightarrow \ {\mathcal E(m)}\otimes p_Y^*{\mathcal O}_{D_j}$.  
Taking direct image to $\mathcal{G}r$ gives the morphism      
$$V\otimes_k H^0({\mathcal O}_{D_j}) \rightarrow (\tilde{E}(m))_{D_j} := (p_{\mathcal{G}r})_*(\mathcal{E}(m)\otimes p_Y^*{\mathcal O}_{D_j})\, .$$ 

Using this,  for each $j$, we have a morphism (see \cite[Subsection 3.3]{Bh2} for details)
$$P_j:  
\mathcal{G}r \longrightarrow Gr(V\otimes_k H^0({\mathcal O}_{D_j}), q_j)\, .$$

Thus we have the morphism 
$$\Psi:= 
({\tilde P},(P_j), P_f)\colon  {\mathcal G}r \longrightarrow Gr(V\otimes W_1, P(m)) \times 
\prod_j Gr(V\otimes_k H^0({\mathcal O}_{D_j}), q_j) \times \times_{x\in I} Gr_x\, ,$$
which is an embedding. 
Let 
$$Z = \Psi ({\mathcal G}r)\, .$$

       The group $G=SL(V)$ acts on $Q$, the action lifts to $\tilde{Q}$ and ${\mathcal G}r$, keeping $R$ invariant. $G$ also has an action on $Gr(V\otimes W_1, P(m))$ and on all $Gr(V\otimes_k H^0({\mathcal O}_{D_j}), q_j)$ and $Gr_x, x \in I$. The morphism $\Psi$ is a $G$-equivariant morphism and $Z$ is a $G$-invariant subset.        
   
\textit{Polarisation}: Recall that $0< \alpha \le 1$ is a fixed rational number and $m'$ is the number of divisors $D_j$. Let $K$  be an integer such that $K \alpha_i$ and $\frac{K \alpha q_j}{n}$ are integers for all $i$ and $j$. Let $d_i(x) =  \alpha_{i+1}(x) - \alpha_i(x)$. Then for
 any integers $(\beta_x)_{x \in I}  \in \mathbb{Z}^I_{\ge 0}$ and $\ell$, satisfying 
\begin{equation}\label{*}
n \sum_{x \in I} \beta_x + \sum_{x \in I} \sum_{i=1}^{l_x} K d_i(x) r_i(x) + n \ell = K N\, ,
\end{equation}
    on $Gr(V\otimes W_1, P(m)) \times \prod_j Gr(V\otimes H^0({\mathcal O}_{D_j}), q_j) \times \times_{x\in I} Gr_x$ we take the polarisation 
\begin{equation} \label{polarisation}
     (\ell - \alpha K \frac{\sum_j q_j}{n}) \ \times \prod_j \alpha K \ \times \ \prod_{x \in I} (\beta_x, d_1(x), \cdots d_i(x), \cdots  d_{\ell_{x}(x)}) )\, .
  \end{equation} 

 Let ${\mathcal G}r_{gpgh}$ (respectively ${\mathcal G}r^{str}_{gpgh}$) be the $G$-invariant closed subscheme of ${\mathcal G}r$ consisting of $(E_*, F(E), \phi)$ such that $\phi$ preserves the GPH structure and  the parabolic structure (respectively strongly preserves the parabolic structure). 
Our required moduli space of good parabolic GPH is the GIT quotient ${\mathcal G}r_{gpgh} // G$  by a linearised action of $G$ for this choice of  polarisation. It is constructed using the morphism $\Psi$. One shows that the GIT (semi)stable points of $Z \cong {\mathcal G}r$ correspond  precisely to the triples $(E_*, F(E), \phi)$ satisfying the parabolic (semi)stability condition \eqref{parGPHstability}.  We define 
$$U^h_{par}(n,d,L_0) \  := \ {\mathcal G}r // G\, ,$$  
the  GIT quotient. The $G$-invariant closed subscheme ${\mathcal G}r_{gpgh}$ determines a closed subscheme 
$$H_{par}(n,d,L_0) \ := \ {\mathcal G}r_{gpgh}//G\, ,$$ 
of  $U^h_{par}(n,d,L_0)$. The $G$-invariant closed subscheme ${\mathcal G}r^{str}_{gpgh}$ determines a closed subscheme 
$$H^{str}_{par}(n,d,L_0) \ := \ {\mathcal G}r^{str}_{gpgh} //G \  \subset H_{par}(n,d,L_0)\, ,$$ 
it is the moduli space of strongly parabolic (good) GPH.

\begin{theorem} \label{parGPHmoduli}
    There exists a coarse moduli space $H_{par}(n,d,L_0)$ of (equivalence classes of) semistable parabolic $L_0$-valued (good) parabolic GPH of rank $n$, degree $d$ on $Y$ containing the open subscheme $H^s_{par}(n,d,L_0)$ of stable (good) parabolic GPH.  There is a closed subscheme $H^{str}_{par}(n,d,L_0) \subset H_{par}(n,d,L_0)$ which is the moduli space of strongly parabolic good GPH.    
\end{theorem}

\subsection{The Theta line bundles} \hfill

      On $Flag_{\tilde Q}$, let $Det R\pi_{Flag_{\tilde Q}} \mathcal {E}$ denote the determinant of cohomology line bundle for the family $\mathcal{E}$. We define 
\begin{equation} \label{thetaf}
 \theta_F: = (Det R\pi_{(Flag_{\tilde Q})} \mathcal {E})^K \bigotimes \otimes_{x \in I} [det (\mathcal{E}_x)^{\beta_x} \otimes \otimes_{i=1}^{l_x} \ det (Q_{x, i})^{d_i(x)}] \bigotimes  (det \mathcal{E}_{x_0})^{\ell}.
\end{equation}
where $(Q_{x,i})_{(E_*, \phi)} = E_x / F'_i(E_x), r_i(x) = \ {\rm dim} \  E_x / F'_i(E_x)$.

Any $c \in k^*$ acts on $V$ and $V\otimes W_1$  by $c \  id$,  on the fibres of $det (\mathcal{E}_x)$ it acts by $c^n id$.  On the fibres of $Det R\pi_{(Flag_{\tilde Q})} \mathcal {E}, c \in k^*$ acts by $c^{- N} id$. On the  fibres of $det Q_{x, i}$ it acts by $c^{r_i(x)} id$. 
Hence on the fibres of $\theta_F$, it acts by $c^{wt \theta_F}$ where 
$$wt \ \theta_F = - NK + \sum_{x \in I} (n \beta_x + \sum_{1}^{l_x} r_i(x) d_i(x) ) + n \ell \, .$$
By condition \eqref{*}, $wt \theta_F = 0$, hence (the restriction of) the line bundle $\theta_F$ descends to an ample  line bundle on the moduli space of parabolic Hitchin sheaves, we denote this line bundle by $\theta_{par}$.  
 
        Let $p_{GF}: \mathcal{G}r \rightarrow Flag_{\tilde Q}$ be the canonical projection. Define 
\begin{equation} \label{thetaGr}
\theta_{\mathcal{G}r} := p_{GF}^* \theta_F \bigotimes \otimes_j (det Q_j)^{\alpha K} 
                                       \bigotimes \otimes (det \mathcal{E}_{x_0})^{- \alpha K \frac{\sum_j q_j}{n}}
\end{equation}

One has $wt \ (det Q_j) = q_j$ hence 
$$wt \ (\otimes_j (det Q_j)^{\alpha K}  \bigotimes \otimes (det \mathcal{E}_{x_0})^{- \alpha K \frac{\sum_j q_j}{n}}) = \sum_j \alpha K q_j  + n (- \alpha K  \frac{\sum_j q_j}{n}) = 0\, .$$
It follows that $wt \ (\theta_{\mathcal{G}r}) = 0$ and (the restriction of)
 $\theta_{\mathcal{G}r}$ descends to an ample line bundle $\theta_{pargph}$ on the parabolic GPH moduli space.

\section{The parabolic Hitchin map}

  Let the notations be as in Subsection \ref{goodparGPH}. In particular, $Y$ is in integral projective curve with the normalisation $X$. In this section, we first establish the relation between the $L_0$-twisted good parabolic GPH on $X$  and $L$-twisted torsion-free parabolic Hitchin sheaves on $Y$. We use it to define a parabolic Hitchin morphism on  the moduli space of $L$-twisted torsion-free parabolic Hitchin sheaves on $Y$.  
  
\subsection{The morphism $f_p$} \hfill

We now associate an $L$-twisted parabolic Hitchin sheaf $({E'}_*, \phi')$ to an $L_0 = p^* L$-twisted good parabolic GPH $(E_*, F(E), \phi)$.   \cite[Proposition 2.8]{Bh1} associates an $L$-twisted Hitchin sheaf $(E', \phi')$ to the good GPH $(E, F(E), \phi)$.  The sheaves $E$ and $E'$ fit in an exact sequence
\begin{equation}\label{E'}
0 \rightarrow E' \rightarrow p_*E \rightarrow p_*(\bigoplus_j \frac{E\otimes {\mathcal O}_{D_j}}{F_j(E)})\rightarrow 0\, .
\end{equation}
Since $p$ is an isomorphism over $Y \setminus \cup_j y_j$, $I \cong p(I)$ and $E' \cong p_*E$ outside $\cup_j y_j$. Hence the parabolic structure on $E$ induces a parabolic structure of the same type on $E'$ and conversely.  The Higgs field $\phi$ preserves (respectively strongly preserves) the parabolic structure on $E$ if and only if the Higgs field $\phi'$  preserves (respectively strongly preserves) the parabolic structure on $E'$.
Define $f_p(E_*, F(E), \phi) = (E_*', \phi')$. As in  \cite[Proposition 2.8(2)]{Bh1}, a parabolic Hitchin bundle $(E_*', \phi')$ determines a unique good parabolic GPH $(E_*, F(E), \phi)$ such that $f_p(E_*, F(E), \phi) = (E_*', \phi')$. The construction extends easily to families. 

      As in \cite[Theorem 1.9]{Bh1}, one can check that (for $\alpha = 1$),  $(E_*, F(E), \phi)$ is $1$-semistable if $(E_*', \phi')$ is semistable and if $(E_*', \phi')$ is stable then $(E_*, F(E), \phi)$ is $1$-stable. Therefore  we have the following theorem. 

\begin{theorem}  \label{mgphcorrespondence}
There exists a birational morphism (for $\alpha=1$ or for $\alpha <1$ and  close to $1$),
\begin{equation} \label{fp}
f_p \colon H_{par} (n,d,L_0) \longrightarrow {\mathcal H}_{par} (n,d,L)
\end{equation} 
from the moduli space of ($p^*L$-twisted) good parabolic GPH of rank $n$, degree $d$ on $X$ to the moduli space of ($L$-twisted) parabolic Hitchin sheaves of rank $n$, degree $d$ on $Y$. 

    Let $N_{par}(n,d,L)$ denote the image of $f_p$ in ${\mathcal H}_{par} (n,d,L)$.
   If $\alpha =1$, then $N_{par}(n,d,L)$ contains all semistable Hitchin bundles. If $\alpha <1$ and  close to $1$, then $N_{par}(n,d,L)$ contains all stable Hitchin bundles.
        
        The restriction of $f_p$ to the subset corresponding to the good GPH which map to stable Hitchin bundles is an isomorphism onto the open subset consisting of stable Hitchin bundles on $Y$. 
\end{theorem}

\subsection{Hitchin  morphism $h_p$ and Parabolic Hitchin morphism $h^{str}_p$} \hfill

            We  recall some notations from the subsection \ref{hitchinmap}.  
             We defined (see \eqref{bfA})
$${\bf A} := \bigoplus_{i= 1}^n \ {\bf H}^0(X, L_0^i)\, , \  {\bf A}' :=  {\bf H}^0(Y, L^i)\, , {\bf A}' \subset {\bf A}\, .$$

There is a Hitchin morphism 
$$h^H_p: H_{par}(n,d,L_0) \to {\bf A}$$ defined by $h_p (E_*, F(E), \phi) = (tr \ \phi, \cdots, tr \ \phi^n)$ where $tr$ denotes the trace and $\phi^i$ is defined inductively as 
$$\phi^i : E \to E \otimes L_0^i\, , \  \phi^i = (\phi \otimes Id) \circ \phi^{i-1}\, .$$ 

Let $H'_{par}(n,d,L_0) \subset H_{par}(n,d,L_0)$ be the subset corresponding to good parabolic GPH which give Hitchin bundles on $Y$. As in \cite[Lemma 2.10]{Bh1}, one shows that
$$h^H_p(H'_{par}(n,d,L_0)) \subset {\bf A}' \, .$$ 

{\bf Properness of the Hitchin morphisms $h^H_p, h_p, h^{str}_p$ : }    

       In fact, the Hitchin morphism is defined on ${\mathcal G}r$ (${\mathcal G}r$ defined in Section 5.1) and being $G$-equivariant defines a Hitchin morphism on $U^h_{par}(n,d, L_0) =  {\mathcal G}r// G$. As in \cite[Theorem 5.1]{Bh1}, one can show that $h_p: U^h_{par}(n,d, L_0) \to {\bf A}$ is proper. At the end of the proof, we need to use the properness of $Gr(n, \ d(D_j), q_j)$ as well as properness of the flag schemes $Flag_{\bar{n}(x)} k^n, x \in I$. Since $H_{par}(n, d, L_0)$ is closed in $U^h_{par}(n, d, L_0)$, it follows that $h^H_p: H_{par}(n, d, L_0) \to {\bf A}$ is proper. 
            
          As in \cite[Corollary 5.2]{Bh1}, one can show that $h^H_p$ (restricted to $H_{par}(n,d,L_0)$) goes down to the Hitchin morphism  
$$h_p: N_{par}(n, d, L) \longrightarrow {\bf A}$$
such that $N'_{par}(n,d,L) := N_{par}(n,d,L) \cap \mathcal{H}'_{par}(n, d, L)$ maps into ${\bf A}'$ under $h_p$.  Since the proper map  
$H_{par}(n, d, L) \to {\bf A}$ factors through $h_p: N_{par}(n, d, L) \to {\bf A}$, it follows that the latter is a proper map.                    
\begin{definition} \label{parhitchinmap}
The restriction of the Hitchin morphism  
$$h^{str}_p: N^{str}_{par}(n, d, L) \rightarrow {\bf A}$$
 is called the parabolic Hitchin morphism.
\end{definition}
Since $N^{str}_{par}(n,d, L)$ is a closed subscheme, the restriction $h^{str}_p$ of the Hitchin map to $N^{str}_{par}(n, d, L)$ is also proper.

\subsection{The parabolic Hitchin base ${\bf A}_p$}    \label{subAp}   \hfill

   For $(E_*, \phi) \in N_{par}(n,d,L)$, write
$$h_p(E_*, \phi) = (a_j(\phi))_j\, , \  a_j(\phi) \in {\bf H}^0(Y, (L(D))^j)\otimes p_*\mathcal{O}_X)\, .$$.  

        For $x \in D$, we reorder the multiplicities $(n_1(x), \cdots, n_{l_x(x)})$ in descending order as $(n'_1(x), \cdots,n'_{l_x(x)})\, , n'_1(x) \ge \cdots \ge n'_{l_x}(x) > 0$. This gives a partition of $n$. Let $(\mu_1(x) \ge \cdots \ge \mu_{l_x}(x) > 0)$ be the conjugate partition of $n$, given by 
        $$\mu_t(x) \ =  \  num( \{ c \ \vert \ n'_c(x) \ge t\, , 1 \le c \le l_x \})\ ,$$
where $num(S)$ denotes the number of elements in the set $S$.        
For $j= 1, \cdots, n$, define $\gamma_j(x)$ (level function) by 
$$ \gamma_j(x) = c \ {\rm if \ and \ only \ if} \ \sum_{t \le c-1}  \mu_t(x) < j \le \sum_{t \le c}  \mu_t(x)\, .$$

\begin{proposition} \label{thmSWW}
(1) For $(E_*, \phi) \in N^{' str}_{par}(n,d,L)$, one has 
$$ 
a_j(\phi) \in {\bf H}^0(Y, L^j \otimes \mathcal{O}_Y(\sum_{x \in D} (j - \gamma_j(x) x)))\, .$$
where $\gamma_j(x)$ are determined explicitly by the multiplicities $n_i(x)$ at $x \in D$.\\
(2) For $(E_*, \phi) \in N^{str}_{par}(n,d,L)$, one has 
$$ 
a_j(\phi) \in {\bf H}^0(Y, L^j \otimes \mathcal{O}_Y(\sum_{x \in D} (j - \gamma_j(x) x)) \otimes p_*\mathcal{O}_X)\, . $$
\end{proposition}

\begin{proof}
 (1)   This follows exactly as in \cite[Theorem 3.4]{SWW} as $E$ is locally free.  See the discussion before \cite[Theorem 3.4]{SWW} for details on $\gamma_j(x)$.  \\  
 (2)  We note that a torsion-free sheaf $E$ on $Y$ is locally free on the set $Y_0$ of nonsingular points of $Y$ and $I \subset Y_0$. Hence , for $(E_*, \phi) \in N^{str}_{par}(n,d,L)$, one has 
$$h_p(E_*, \phi) = (a_j(\phi))_j\, , \ 
a_j(\phi) \in H^0(Y, L^j \otimes \mathcal{O}_Y(\sum_{x \in D} (j - \gamma_j(x) x)) \otimes p_*\mathcal{O}_X)\, . $$
where $\gamma_j(x)$ are determined explicitly by the multiplicities $n_i(x)$ at $x \in D$.  
\end{proof}
  
Define 
\begin{equation} \label{bfAp}
{\bf A}_p: = \prod_{j= 1}^n  {\bf H}^0(Y, L^j \otimes \mathcal{O}_Y(\sum_{x \in D} (j - \gamma_j(x) x)) \otimes p_*\mathcal{O}_X)\, , 
\end{equation}
and 
\begin{equation} \label{bf A'p}
{\bf A}'_p: = \prod_{j= 1}^n  {\bf H}^0(Y, L^j \otimes \mathcal{O}_Y(\sum_{x \in D} (j - \gamma_j(x) x))\, . 
\end{equation}
By Proposition \ref{thmSWW}, $h_p(N^{' str}_{par}(n,d,L)) \subset {\bf A}'_p\, ; \ h_p(N^{str}_{par}(n,d,L)) \subset {\bf A}_p$. The parabolic Hitchin map   
$$h_p: N^{str}_{par}(n,d,L)) \rightarrow {\bf A}_p$$
is a proper morphism.

\subsection{Assumptions:}  \label{assumptions} \hfill

Henceforth we assume that the parabolic weights are general so that semistability is same as stability. 
We also assume that $Y$ is an integral projective curve with at most nodes as singularities.

             We note that, since $Y$ is an integral projective curve with at most nodes as singularities, the moduli space $U_{par}(n,d)$ of semistable parabolic sheaves on $Y$ of rank $n$ and degree $d$ with a fixed parabolic structure is irreducible (which can be proved using parabolic GPBs) and one has   
$$\ {\rm dim} \ U_{par}(n,d) = n^2(g-1) +1 + \sum_{x \in D} \ {\rm dim} \ (GL(n)/P_x) \, .$$

In case $L = \omega_Y$,
$$\ {\rm dim} \ {\mathcal H}_{par}(n,d, \omega_Y) = n^2(2g- 2 + deg(D)) +1\, ; $$
$$\ {\rm dim} \ {\mathcal H}^{str}_{par}(n,d, \omega_Y) =  2 \ {\rm dim} \ U_{par}(n,d)\, .$$

\subsection{General fibre of the parabolic Hitchin map $h_p$ over $a \in {\bf A'}$}  \hfill

     In the parabolic case, in subsection \ref{hitchinfibre}, we  replace $L$ by $L(D)$ where $D$ denotes the  parabolic divisor, ${\bf A}' =  \prod_{i= 1}^n {\bf H}^0(Y, (L(D))^i)$. Using the local equation of $Y_a$, one sees that the set of $a \in {\bf A}'$ such that $Y_a$ is integral and the projection map $\pi: Y_a \to Y$ is un-ramified over the nodes of $Y$ is a nonempty open set. Hence we may assume that for a general $a \in {\bf A}' $, the spectral curve $Y_a$ is an integral nodal curve un-ramified over nodes of $Y$.  
     
              For simplicity, assume that $D = x$, a single point and $y_i, i= 1, \cdots n$ the points of $Y_a$ lying over $x$.  In the notation of the proof of Proposition \ref{fibreainA'}, let $\mathbb{V}$ denote the direct image (under $id \times \pi$) on $\overline{J}(Y_a)\times Y$ of the universal sheaf $\mathcal{P}$ on $\overline{J}(Y_a)\times Y_a$. Note that the restriction $\mathbb{V}_x = \mathbb{V}\vert_{\bar{J}(Y_a)\times x}$ is a vector bundle (i.e., a locally free sheaf) of rank $n$.
Let $Fr(\mathbb{V}_x)$ be the ${\bf G}= GL(n)$-bundle (frame bundle) corresponding to the vector bundle $\mathbb{V}_x$. Let $P_x$ be the parabolic subgroup of ${\bf G}$ determined by the parabolic structure. Let $W_x$ be  the Coxeter subgroup corresponding to $P_x$. Define
   $$\mathcal{F}l := Fr(\mathbb{V}_x) \times_{\bf G} {\bf G}/P_x\, ,$$
   the associated flag bundle over $\overline{J}(Y_a)$.  

              Restricting the Higgs field on $\mathbb{V}$ to $\mathbb{V}_x$, we have
$$\theta_x: \mathbb{V}_x \rightarrow \mathbb{V}_x\, .$$
 Then $\theta_x = \oplus_i  \theta_{y_i}\, , \theta_{y_i} :=\theta_{\mathcal{P}_{y_i}, x}$, 
 it is regular semisimple everywhere. One has $\theta_x \in \mathcal{A}ut \mathbb{V}_x,$ the group scheme of local automorphisms of the vector bundle $\mathbb{V}_x$.  The centraliser of $\theta_x$ in $\mathcal{A}ut \mathbb{V}_x$ is a group scheme  $\mathcal{T}$ over $\overline{J}(Y_a)$ whose fibre over $N \in \overline{J}(Y_a)$ is a maximal torus $T$ in $Aut( (\mathbb{V}_x)_N)= GL(n)$. Being a subgroup scheme of  $\mathcal{A}ut \mathbb{V}_x$, the group scheme $\mathcal{T}$ acts on $\mathcal{F}l$. Since $({\bf G}/P_x)^T$ is bijective to $W_x$, it follows that the fixed point set $\mathcal{F}l^{\mathcal T}$ is a $W_x$-torsor over $\overline{J}(Y_a)$.

\begin{theorem} \label{parabolic Hitchin fibre}
For a general $a \in {\bf A}'$, there is an injective morphism 
$$h_{\mathcal F}: \mathcal{F}l^{\mathcal T} \to \mathcal{H}_{par}(n,d,L)\, .$$  The parabolic Hitchin fibre $(h_p)^{-1} a = h_{\mathcal F} (\mathcal{F}l^{\mathcal T}) \cap N_{par}(n,d,L)$.   

This morphism is an isomorphism from the restriction $\mathcal{F}l^{\mathcal T}\vert_{J(Y_a)}$ onto the intersection of $(h_p)^{-1}a$ with the parabolic Hitchin bundles.

\end{theorem}
\begin{proof}
      Let $\pi: Y_a \to Y$ be the spectral curve associated to $a$. As in Subsection \ref{hset} (replacing $L$ by $L(D)$), one has a set map $h_{p,set}: \mathcal{H}_{par}(n,d,L) \to A$. The parabolic Hitchin sheaf $(E_*, \theta)$ lies in $h_{p, set}^{-1} a$ if and only if there is a torsion-free sheaf $\mathcal{L}$ of rank $1$ on $Y_a$ such that $\pi_* {\mathcal L} = E$ and the filtration at $x \in D$ is compatible with $\theta_x$. 
           So $(E_*, \theta) \in h_{p,set}^{-1}a$ if and only if $(E_*, \theta)$ belongs to the (image of the) fibre of $\mathcal{F}l$ over ${\mathcal L}$ and the filtration at $x\in D$ is compatible with $\theta_x$. As the filtration of $\mathbb{V}_x$ is determined by a subbundle $P'$ of parabolic subgroups, $P' \subset {\mathcal A}ut \mathbb{V}_x$, the filtration is compatible with $\theta_x$ if and only if $\theta_x \in P'$. As $\theta_x$ is regular and semisimple, $\theta_x \in P'$ if and only if $\theta_x \in {\mathcal T}$ so that $(E_*, \theta)$ is in (the image of) $\mathcal{F}l^{\mathcal T}$. Thus $h_{\mathcal F}$ maps $\mathcal{F}l^{\mathcal T}$ bijectively onto its image $h_{\mathcal F}(\mathcal{F}l^{\mathcal T})$ whose underlying set of points is $h_{p,set}^{-1} (a)$. In particular it maps injectively into $\mathcal{H}_{par}(n,d,L)$ and $h_p^{-1} (a) = h_{\mathcal F}(\mathcal{F}l^{\mathcal T}) \cap N_{par}(n,d,L)$.
            
          The last statement follows as in the above correspondence, $E$ is locally free if and only if $\mathcal{L}$ is locally free.      
\end{proof}

\begin{corollary}\label{coro5.6SWW}
Let $a \in {\bf A}'$ be a general point.\\
(1)   The set of connected components $\pi_0((h_p)^{-1}(a))$ is in bijective correspondence to the Coxeter group $W_x$ associated with the parabolic subgroup $P_x$ determined by the parabolic structure at $x$.\\
(2)    For any $a$, dim $(h_p)^{-1} a \ge  \frac{n(n-1)}{2}(d(L) + deg(D)) + n(g-1) + 1$.    
\end{corollary}
\begin{proof}
(1)  Since $\mathcal{F}l^{\mathcal T}$ is a $W_x$-torsor over $\bar{J}(Y_a)$, the corollary follows from 
Theorem \ref{parabolic Hitchin fibre}. \\
(2)   The dimension of a fibre is at least equal to the dimension of the general fibre.  By Theorem \ref{parabolic Hitchin fibre} the general fibre has dimension equal to dim $J(Y_a) = g(Y_a)$. Since 
\begin{equation} \label{g(Ya)}
g(Y_a) =   \frac{n(n-1)}{2}(d(L) + deg(D)) + n(g-1) + 1,
\end{equation}
 \cite{BNR} the part (2) follows. 
 \end{proof}

     Let ${\bf H}_{par}(n,d,L)$ be the closure of $H'_{par}(n,d,L)$ in $N_{par}(n,d,L)$. Since $h_p(H'_{par}(n,d,L)) \subset {\bf A'}$ and ${\bf A'}$ is closed in ${\bf A}, h_p({\bf H}_{par}(n,d,L)) \subset {\bf A'}$. In view of 
Theorem \ref{parabolic Hitchin fibre} and Corollary \ref{coro5.6SWW}, the restrictions of $h_p$ to $H'_{par}(n,d,L)$ and ${\bf H}_{par}(n,d,L)$ are dominant. 

In case $Y$ has only one node, denote by $\mathcal{H}_{i, par}(n,d,L) \subset \mathcal{H}_{par}(n,d,L)$ the subset corresponding to Hitchin sheaves with local type $n - i$. Let $N_{i, par}(n,d,L) = N_{par}(n,d,L) \cap \mathcal{H}_{i, par}(n,d,L)$.

\subsection{The general fibre of the parabolic Hitchin map $h^{str}_p$ over $a \in {\bf A}'_p$.} \hfill

     In the strongly parabolic case, for $a \in {\bf A}'_p$ the  spectral curve $Y_a$ is totally ramified over points of $D$ (as is seen from its local equation). As in \cite[Lemma 7.1]{SWW}, we can show that the set of $a \in {\bf A}'_p$ such that $Y_a$ is integral, totally ramified over points of $D$ (and with nodes over the nodes of $Y$ outside $D$) is a nonempty open set.  We consider $a \in {\bf A}'_p$ general enough so that the spectral curve $Y_a$ is a singular curve which is an $n$-sheeted covering of $Y$ unramified over the nodes of $Y$ and totally ramified over the points of $D$. Let $Y^p_a$ be the partial normalisation of $Y_a$ obtained by resolving its singularities at the points lying over $x \in D$. Let $\pi: Y_a \to Y, \pi': Y^p_a \to Y_a$ be the natural maps, 
     $$\pi_p = \pi \circ \pi': Y^p_a \longrightarrow Y$$ 
be the composite.  

\begin{proposition} \label{Ypa}
      For a general $a \in {\bf A}'_p$, there are $n'_1(x)$ points in $Y^p_a$ lying over $x \in D$ with ramification degrees $(\mu_1(x), \cdots, \mu_{n'_1(x)})$ (with $n'_1(x)$ and $\mu_i(x)$ defined in Subsection  \ref{subAp}).
      
            The arithmetic genus of $Y^p_a$ is given by
$$g(Y^p_a) = \frac{n(n-1)}{2} d(L) + n(g-1) +1 + \sum_{x\in D} \ dim \ GL(n)/P_x\, .$$      
In particular for $L = \omega_Y$, 
$$g(Y^p_a) = n^2(g-1)+1 + \sum_{x\in D} \ dim \ GL(n)/P_x\, .$$          
\end{proposition}

\begin{proof}
    We note that since $D$ consists of nonsingular points of $Y$, the local analysis in \cite[section 4]{SWW} is valid in the nodal case also.  Hence the proposition can be proved exactly as \cite[Corollary 4.8]{SWW}.  
    
    We have $g(Y_a) - g(Y^p_a) = n(n-1) deg(D) - dim \ GL(n)/P_x$ (as in \cite[Ch. 5, Example 3.9.2]{Ha}) and $g(Y_a)$ is given by \eqref{g(Ya)} (computed in Proposition \ref{fibreainA'}).
 Notice that 
 \begin{equation} \label{g(Ypa)}
 g(Y^p_a) = \ {\rm dim} \ U_{par}(n,d) + \frac{n(n-1)}{2}(d(L) - (2g-2))\, .
 \end{equation}  
\end{proof}

\begin{lemma} \label{SWWL3.3}
Let $d(L) \ge 2g-2, D$ nonempty. Then
$${\rm dim} \ A'_p \ = \ {\rm  dim} \ U_{par}(n,d) + \frac{n(n+1)}{2} (d(L) - (2g-2)) - 1\, , \ {\rm if} \ L \neq \omega_Y\, ,$$
$${\rm dim} \ A'_p \ = \ {\rm  dim} \ U_{par}(n,d)  \, , \ {\rm if} \ L = \omega_Y\, ,$$ 
where 
$${\rm dim} \ U_{par}(n,d) \ = \ n^2(g-1) +1 +  \sum_{x} \ {\rm dim} \ GL(n)/ P_x\, .$$    
\end{lemma}
\begin{proof}
Suppose that $L \neq \omega_Y$. 
$$
\begin{array}{lll}
{\rm dim} \ {\bf A}'_p & := & \sum_{j=1}^n \ {\rm dim} \ H^0(Y, L^j \otimes_{\mathcal{O}_Y} ( \sum_{x \in D} (j - \gamma_j(x) x))) \\
{} & = & \sum_{j=1}^n  ( d(L^j) + \sum_{x \in D} (j - \gamma_j(x))  + 1- g) \\
{} & = & (\sum_{j=1}^n j) d(L)  - n(g-1) + \sum_{j=1}^n  \sum_{x \in D} (j - \gamma_j(x))  \\
{} & = & \frac{n(n+1)}{2} d(L) - n(g-1) +  \sum_{x \in D} \sum_{j=1}^n (j - \gamma_j(x)) \\
{} & = & \frac{n(n+1)}{2} d(L) - n(g-1) +  \sum_{x \in D} \ {\rm dim} \ GL(n)/P_x \\
{} & = & \ {\rm dim} \ U_{par}(n,d) + \frac{n(n+1)}{2} (d(L) - (2g-2))  - 1\, . 
 \end{array} 
  $$
 Since ${\rm dim} \ H^0(Y, \omega_Y) = g$, we similarly get  
 $${\rm dim} \ A'_p \ = \ {\rm  dim} \ U_{par}(n,d) \ {\rm if} \  L = \omega_Y\, .$$
\end{proof}

            For a general $a \in {\bf A}'_p$,  the Hitchin fibre has dimension equal to $g(Y^p_a)$. Since the Hitchin map  ${\mathcal H'}^{str}_{par}(n,d,L) \to {\bf A}'_p$ is dominant, one has 
${\rm dim} \ {\mathcal H'}^{str}_{par}(n,d,L) =  {\rm dim} \ A'_p + g(Y^p_a)$. Suppose first that $L \neq \omega_Y$. Then using Lemma \ref{SWWL3.3}, we have 
$${\rm dim} \ {\mathcal H'}^{str}_{par}(n,d,L) \ = \ 2 \ {\rm  dim} \ U_{par}(n,d) + n^2 (d(L) - (2g-2)) - 1\, .$$
 Similarly one has 
 \begin{equation} \label{dimNstrpar}
 {\rm dim} \ {\mathcal H'}^{str}_{par}(n,d,L) \ = \ 2 \ {\rm  dim} \ U_{par}(n,d) \, , \ {\rm for} \ L = \omega_Y\, .
\end{equation}

\subsection{Strongly Parabolic BNR Correspondence} \hfill
            
             We prove the following theorem generalising the BNR correspondence in the non-parabolic case to the strongly parabolic case.
                             
\begin{theorem} \label{parBNRcorrA'}
Let $a \in {\bf A}'_p$ be a general element and let 
$$\delta = d + \frac{(n^2 - n)}{2}d(L) + \sum_{x\in D} \ {\rm dim} \ (GL(n)/P_x)\ . $$
Then we have:\\
(1) The direct image map $(\pi_p)_*$ induces an injective morphism 
$$h_T: \bar{J}^{\delta}(Y^p_a) \longrightarrow \mathcal{H}^{str}_{par}(n,d,L)\, .$$  
(2) It induces an isomorphism $J(Y^p_a) \cong (h_p^{' str})^{-1} (a)$. \\
(3) The parabolic Hitchin fibre $(h^{str}_p)^{-1} (a) = h_T(\bar{J}^{\delta}(Y)) \cap N^{str}_{par}(n,d,L)$. 
\end{theorem}

\begin{proof}
       
       The proof is similar to the proof of \cite[Theorem 4.12]{SWW}.
          Given a torsion-free sheaf $N_p$ of rank $1$ and degree $\delta$ on $Y^p_a$, the direct image map gives a Hitchin sheaf $(E, \phi)$ on $Y$ where $E = (\pi_p)_* N_p$, the Hitchin field is induced by a section of $\pi_p^*L$ and $h_{p, set} (E, \phi) = a$.  We note that $E$ is locally free if and only if $N_p$ is so.  Since $\pi_p$ is a finite morphism, $\chi(N_p) = \chi((\pi_p)_*N_p)$ i.e. $\delta + (1- g(Y^p_a)) = d + n(1-g)$ by Riemann-Roch theorem.  Substituting for $g(Y^p_a)$ from \eqref{g(Ypa)}, one gets the expression for $\delta$ for $r(E) =n, d(E) =d$.
   By \cite[Proposition-Definition 2]{SWW}, $E= (\pi_p)_* N_p$ has a canonical Young diagram filtration with parabolic type $P_x$ for each $x \in D$ compatible with the Higgs field $\phi$  (see the proof of \cite[Theorem 4.12]{SWW} for details).  Thus we have a parabolic Hitchin sheaf $(E_*, \phi)$ given by $N_p$. This can be checked to define a morphism $h_T:  \bar{J}^{\delta}(Y) \longrightarrow \mathcal{H}^{str}_{par}(n,d,L)\, .$
                  
       Given a strongly parabolic Hitchin sheaf $(E_*, \phi)$ with $h_{p, set}(E, \phi) = a$, we have $E = \pi_*N$ for some rank $1$ torsion-free sheaf $N$ on $Y$.  The parabolic structure on $E$ over $x\in D$ corresponds to a filtration on $N_{\pi^{-1}x}$. Since $E$ is locally free at points of $x \in D$, $N$ is locally free at points of $\pi^{-1}x$ for all $x \in D$. Hence by \cite[Proposition 4.3]{SWW}, $N$ has a $\pi'_* \mathcal{O}_{Y^p_a}$-structure. Thus there exists a torsion-free sheaf $N_p$ of rank $1$ and degree $\delta$ on $Y^p_a$ such that $(\pi_p)_* N_p = E$. The $\pi'_* \mathcal{O}_{Y^p_a}$-structure on $(\pi_p)_* N_p$ is equivalent to giving $\phi$ on $E_*$ i.e. $\phi$ strongly parabolic. Thus $h_T$ maps bijectively onto its image in $\mathcal{H}^{str}_{par}(n,d,L)$. It follows that the parabolic Hitchin fibre $(h^{str}_p)^{-1} (a) = h_T(\bar{J}^{\delta}(Y)) \cap N^{str}_{par}(n,d,L)$. 
      
\end{proof}

\begin{corollary}
       For any $a \in {\mathcal U}$ and $(E_*, \phi)$ with $h_p (E_*, \phi) = a$, the Jordon blocks of $\phi$ mod $t_x$ are of size $(\mu_1(x), \cdots, \mu_{n'_1}(x))$ ($t_x$ being the local parameter at a parabolic point $x$).
\end{corollary}

\begin{proof}
This can be proved as \cite[Corollary 4.15]{SWW}.  
\end{proof}

\subsection{The fibre of the parabolic Hitchin map $h^{str}_p$ at a general $a \in {\bf A}_p$} \hfill

         Since ${\bf A}'_p$ is a closed subset of ${\bf A}_p$, we may assume that a general $a \in {\bf A}_p$ is in ${\bf A}_p \setminus {\bf A}'_p$. 
           By subsection \ref{hitchinfibreA}, for a general $a \in {\bf A} \setminus {\bf A}'$ there is a spectral curve $Y_a$ totally ramified over nodes of $Y$. There is an injective morphism from the compactified Jacobian $\bar{J}(Y_a)$ into $\mathcal{H}(n,d,p_*L_0)$.  Its restriction to  the generalised Jacobian $J(Y_a)$ is an isomorphism onto the image.  For a general $a \in {\bf A}_p \setminus {\bf A}'_p$, the curve $Y_a$ is also totally ramified at all $x \in D$. Let $Y^p_a$ denote the partial normalisation of $Y_a$ obtained by normalising it at points over $x \in D$. 
 
           By arguments similar to those in \cite[Chapter 5, Ex. 3.9.2]{Ha}, for a general $a \in {\bf A}_p \setminus {\bf A}'_p$ we have 
           $$g(Y_a) - g(Y^p_a) = \frac{n(n-1)}{2} deg(D) -  \sum_x \ {\rm dim} \ GL(n)/P_x\, .$$
By Proposition \ref{fibreainA},  $g(Y_a) = \frac{n(n-1)}{2}  d(p_*(L_0(D))) + n(g-1) +1$ and $d(p_*(L_0(D))) = d(L_0(D)) + m'$ ($m'$ being the number of nodes of $Y$), we get  
$$g(Y^p_a) = \frac{n(n-1)}{2}  d(p_*L_0) + n(g-1) +1 + \sum_x \  {\rm dim} \ GL(n)/P_x\, .$$.

Since $\pi_p$ is finite, for a line bundle $N$ on $Y^p_a$ we have $\chi(N) = \chi((\pi_p)_* N)$. Hence if $d(N) =  \delta_A, d((\pi_p)_* N) = d$, then  
$\delta_A + (1- g(Y^p_a)) = d + n(1- g)$ so that $\delta_A = d + \frac{n(n-1)}{2}(d(L) + m') + \sum_x \  {\rm dim} \ GL(n)/P_x\, .$
         
                    By arguments similar to those in the proof of Theorem \ref{parBNRcorrA'}, one has the following:
\begin{theorem} \label{paraBNRcorrA}

Let $a \in {\bf A}_p$ be a general element and $\delta_A = d+ \frac{(n^2 - n)}{2}(d(p_*L_0) + \sum_{x\in D} \ dim \ (GL(n) /P_x) $. Then we have:\\
(1) The direct image map $(\pi_p)_*$ induces an injective morphism 
$$h_T: \bar{J}^{\delta}(Y^p_a) \longrightarrow \mathcal{H}^{str}_{par}(n,d,p_*L_0)\, .$$  
(2) It induces an isomorphism $J(Y^p_a) \cong (h_p^{' str})^{-1} (a)$. \\
(3) The parabolic Hitchin fibre $(h^{str}_p)^{-1} (a) = h_T(\bar{J}^{\delta}(Y)) \cap N^{str}_{par}(n,d,L)$. 

\end{theorem}

\subsection{The spectral curves $X^p_a$ and $Y^p_a$.}  \hfill

Define 
${\bf A}_X = \sum_{i=1}^n {\bf H}^0(X, L_0^i)\, ,  {\bf A}_X \cong {\bf A}$ and let ${\bf A}_{X,p}$ be the corresponding parabolic Hitchin base for the parabolic structure over the divisor $p^{-1} D \subset X$.    
Let $a \in {\bf A}_p \setminus {\bf A}'_p \subset {\bf A}'_{1, p}$ general enough and 
$$a_X := p^*a \in {\bf A}_{X, p} \subset {\bf A}'_{X, 1, p} \subset 
\prod_{i=1}^n {\bf H}^0(X, p^*(L_1)^i)\, , p^*L_1 = L_0(\sum_{j=1}^{m'} (x_j+z_j))\, .$$              
There is an integral spectral curve $X_{a_X}$ associated to  $a_X$. By \cite[Lemma 5.3]{Bh1}, there is a birational surjective morphism $p_P: X_{a_X} \to Y_a$ making the following diagram commutative. 
$$
\begin{array}{ccc}
X_{a_X} & \stackrel{p_P}{\longrightarrow} & Y_a\\
{} & {} & {}\\
\downarrow \pi_X & {} & \downarrow \pi \\
{} & {} & {}\\
X & \stackrel{p}{\longrightarrow} & Y\, .
\end{array}
$$
Moreover, $p_P$ is an isomorphism over the nonsingular points of $Y_a$. 
Define 
$$X^p_a = X_{a_X}  \times_{Y_a} Y^p_a\, , \pi'_X: X^p_a \to X_{a_X}\, , p'_P: X^p_a \to Y^p_a $$
the projections and let $\pi_{X,p}: X^p_a \to X$ denote the composite $\pi_X \circ \pi'_X$.
Then $p'_P$ is an isomorphism over $Y^p_a \setminus \pi^{-1} (\cup_j y_j)\, .$  Since $Y^p_a$ is nonsingular over the points of $D$, it follows that $X_a$ is nonsingular over points of $p^{-1}D$.      
Hence $X^p_a$ is the parabolic spectral curve associated to $a_X := p^*a \in {\bf A}_{X, p}$. 
We have a commutative diagram
$$
\begin{array}{ccc}
X^p_{a} & \stackrel{p'_P}{\longrightarrow} & Y^p_a\\
{} & {} & {}\\
\downarrow \pi_{X, p} & {} & \downarrow \pi_p \\
{} & {} & {}\\
X & \stackrel{p}{\longrightarrow} & Y\, .
\end{array}
$$

\begin{proposition} \label{paraontoA}
Let $\mathcal{H}^{str}_{0, par}(n,d,L)$ be the subset of $\mathcal{H}^{str}_{par}(n,d,L)$ consisting of Hitchin sheaves with local type $b_j= n$ at all the nodes $y_j$.  Then the Hitchin map is well defined on $\mathcal{H}^{str}_{0, par}(n,d,L)$ and it maps $\mathcal{H}^{str}_{0, par}(n,d,L)$ onto ${\bf A}_p$. 
$$(h^{str}_p)^{-1} a \cap  \mathcal{H}_{0, par}(n, d, L) \cong  (h^{str}_{X , p})^{-1} a\, .$$
\end{proposition}

\begin{proof}
     Using the parabolic GPH, as in \cite[Proposition 4.3]{Bh2}, we see that $\mathcal{H}^{str}_{0, par}(n,d,L) = N^{str}_{0, par}(n,d,L)$ so that the parabolic Hitchin map $h^{str}_p$ is well defined on $\mathcal{H}^{str}_{0, par}(n,d,L)$. We have a commutative diagram 
$$
\begin{array}{ccccc}
\mathcal{H}^{str}_{X, par} (n, d - m' n, L_0) & \stackrel{p_*}{\cong}  & \mathcal{H}^{str}_{0, par}(n, d, L) & =  &  N^{str}_{0, par}(n, d, L)\\
{} &        {}& {}  & {}   & {} \\
\downarrow h^{str}_{X , p} & {} & \downarrow h^{str}_p & {} & \downarrow h^{str}_p\\   
{} &        {}& {}  & {}   & {} \\
{\bf A}_p & {=} & {\bf A}_p  & {=}   & {\bf A}_p 
\end{array}
$$     
Since $h^{str}_{X , p}$ is surjective, the restriction $h^{str}_p: \mathcal{H}_{0, par}(n, d, L) \to {\bf A}_p$ is surjective. 

              It follows that $(h^{str}_p)^{-1} a \cap  \mathcal{H}_{0, par}(n, d, L) \cong  (h^{str}_{X , p})^{-1} a$.  Since $(h^{str}_{X , p})^{-1} a$ consists of $(\pi_{X,p})_* N, N$ a torsion-free sheaf of rank $1$ (and of a suitable degree) on $X^p_a$, $(h^{str}_p)^{-1} a \cap  \mathcal{H}_{0, par}(n, d, L)$ consists of $(p\circ \pi_{X,p})_* N, N$ a torsion-free sheaf of rank $1$ (and of a suitable degree) on $X^p_a$. 
\end{proof}

\section{Nilpotent cones of parabolic Hitchin maps} \label{nilconenonfixedcase}

In this section, we assume that the Higgs field has values in $\omega_Y(D)$.
\begin{definition} \label{nilp}
     Define
\begin{equation}      Nil_p := h_p^{-1}(0) \subset {N}_{par}(n,d,L)\ ,
\end{equation}
\begin{equation} \label{nil'p}
       Nil'_p = h_p^{-1}(0) \cap \mathcal{H}'_{par}(n,d,L) = h_p^{-1}(0) \cap \mathcal{H}'_{par}(n,d,L) \cap N_{par}(n,d,L)\, ,
\end{equation}       
\begin{equation}      Nil^{str}_p := (h^{str}_p)^{-1}(0) \subset {N}^{str}_{par}(n,d,L)\, , 
\end{equation}
\begin{equation} \label{nil'strp}      
Nil^{' str}_p = (h^{str}_p)^{-1}(0) \cap \mathcal{H}'_{par}(n,d,L) = (h^{str}_p)^{-1}(0) \cap \mathcal{H}'_{par}(n,d,L) \cap {N}^{str}_{par}(n,d,L)\, .
\end{equation}

\end{definition}  
    
  The subset $Nil'_p$ (respectively $Nil^{' str}_p$) is an open subset of $Nil_p$ (respectively $Nil^{str}_p$).
 For a general $a \in {\bf A}'_p$,  the Hitchin fibre has dimension equal to $g(Y^p_a)$ (by Proposition \ref{Ypa}), hence the special fibre has dimension at least $g(Y^p_a)$.  By Equation (\ref{g(Ypa)}), one has 
 $${\rm  dim} \  (h^{str}_p)^{-1}(0) \ge \ {\rm dim} \  U_{par}(n,d)\, .$$  
 
\subsection{The dimension of $Nil^{' str}_p$ and $Nil^{str}_p$}  \label{dimNilstrp}  \hfill

        We shall determine the dimension of $Nil^{' str}_p$  using infinitesimal deformations following \cite[Section 6]{SWW}. We briefly sketch the proof, see  \cite[Section 6]{SWW} for details.

        There exists a scheme $Flag_{\tilde Q}$ with a local universal family $({\mathcal E}_*, \Phi)$ for Hitchin sheaves $(E_*,\phi)$ where $E_*$ is a parabolic vector bundle of the fixed type and $\phi$ is a Hitchin field on $E$ (see Subsection \ref{moduliparGPH} for details). Let $\tilde{Q}^p$ (respectively $\tilde{Q}^{strp}$ be the closed subset of $Flag_{\tilde Q}$ consisting of Hitchin sheaves $(E_*,\phi)$ such that $\phi$ preserves (respectively strongly preserves) the parabolic structure. The moduli space $\mathcal{H}_{par}(n,d,L)$ of (weakly) parabolic Hitchin sheaves  (respectively $\mathcal{H}^{str}_{par}(n,d,L)$ of strongly parabolic Hitchin sheaves) is a GIT quotient of $\tilde{Q}^p$ (respectively $\tilde{Q}^{strp}$) by $SL(V)$. For generic parabolic weights, the moduli space of strongly parabolic Hitchin bundles $\mathcal{H}^{' str}_{par}(n,d,L)$ is a  geometric quotient of $\tilde{Q}^{strp}$ by $SL(V)$. Let 
\begin{equation}   \label{q}       
          q: \tilde{Q}^{strp} \to \mathcal{H}^{str}_{par}(n,d,L)
\end{equation}          
be the quotient map. One has a surjective homomorphism $V \otimes_k \mathcal{O}_{Y \times \tilde{Q}^{strp}} \to \mathcal{E}_*$. Restricting it to $q^{-1} Nil^{' str}_p$ we get a universal family $\mathbb{U}_N:= (V \otimes_k \mathcal{O}_{Y \times q^{-1} Nil^{' str}_p} \to \mathcal{E}_*, \Phi)$ satisfying the following universal property. For a scheme $S$ and  a flat family of quotients $(V \otimes_k \mathcal{O}_{Y \times  S} \to \mathcal{E}_S, \phi_S)$ of strongly parabolic Hitchin bundles with $h^{str}_p(\mathcal{E}_S, \phi_S) = 0$, there is a morphism $f_S: S \to Nil^{' str}_p$ such that $(Id \times f_S)^*(\mathbb{U}_N) \cong (V \otimes_K \mathcal{O}_{Y \times  S} \to \mathcal{E}_S, \phi_S)$. 
          
          It suffices to determine the dimension of each irreducible component of  $Nil^{' str}_p$ with reduced structure. Let $\eta$ be the generic point of an irreducible component of $Nil^{' str}_p$ and let $\phi_{\eta}: \Phi\vert_{q^{-1} (\eta)^{red}}$ be the restriction. Then $\phi_{\eta}$ induces the  canonical filtration $\{Ker(\phi^i_{\eta})\}$ of  $\mathcal{E} \vert_{q^{-1} (\eta)^{red}}$. We note that  torsion-free  sheaves on nodal curves are limits of vector bundles. We may assume that graded terms in the canonical filtration are vector bundles. Then  there is an open subset $U$ of $q^{-1} (\eta)^{red}$ with generic point $\eta$ such that $\phi_{U} := \Phi\vert_{q^{-1} (U)^{red}}$ induces a canonical decreasing filtration of $E_U:= \mathcal{E} \vert_{q^{-1} (U)^{red}}$ given by  $K^j_U= Ker \ \phi^{r-j}_U$ i. e.,
$$K^{\bullet} _U: E_U  = K^0_U \supset K^1_U \supset \cdots \supset K^{r}_U = 0\, ,$$ 
with  $K^j_U$ locally free for all $j= 0, \cdots, r$.       
One can check that at any $x \in D$ the canonical filtration is coarser than the parabolic filtration. 

For a parabolic vector bundle $E_*$, let $ParEnd(E_*)$ (respectively $SParEnd(E_*)$) denote the sheaf of parabolic (respectively strongly parabolic) endomorphisms of $E_*$. For filtered vector bundles $E_1, E_2$ with decreasing filtrations, let $Hom^{filt}(E_1, E_2)$ (respectively $Hom^{sfilt}(E_1, E_2)$) denote the sheaf of local homomorphisms preserving (respectively strongly preserving) the filtrations. For $E_1=E=E_2$, these sheaves will be denoted by $End^{filt}(E,E)$ and $End^{sfilt}(E,E)$ respectively. The following theorem can be proved exactly as \cite[Theorem 6.4]{SWW}.

\begin{theorem} \label{deformations}
The space of infinitesimal deformations of a nilpotent strongly parabolic Hitchin bundle $(E_*,\phi)$  in $U$ is canonically isomorphic to $\mathbb{H}(Y, \mathcal{A}^{\bullet})$. Here $\mathcal{A}^{\bullet}$ denotes the complex of sheaves
$$\mathcal{A}^{\bullet}: 0 \to ParEnd (E_*) \cap End^{filt} E \stackrel{ad(\phi)}{\longrightarrow} (SParEnd (E_*))\cap End^{sfilt}(E)) \otimes \omega_Y(D) \to 0\, .$$
\end{theorem}
  
           The canonical filtration on $E$ determines a reduction $E_P$ (of the frame bundle) of $E$ to a  parabolic subgroup $P \subset GL(n)$. Let ${\bf p} \subset {\bf g}$ be the Lie algebras of $P$ and $GL(n)$. Let ${\bf n}$ be the Lie algebra of the unipotent radical of $P$. Then $\check{\bf n} \cong {\bf g}/ {\bf p}$ as $P$-linear representation \cite[Lemma 6.5]{SWW}. 
                
              For $x \in D$, the parabolic structure determines a parabolic subgroup $P_x \subset GL(n)$, let ${\bf p}_x$ be its Lie algebra.  Since the parabolic filtration is finer than the canonical filtration, $P_x \subset P$,  one has an exact sequence 
$$0 \to ParEnd(E_*) \cap End^{filt}(E) \to End^{filt}E \to \bigoplus_{x \in D} \to {i_x}_* \ {\bf p} /\ {\bf p}_x\, ,$$
where $i_x: x \to Y$ is the inclusion.  
It gives an exact sequence 
$$0 \to \mathcal{A}^{\bullet} \to \mathcal{A}^{' \bullet} \to {i_x}_* \ {\bf p} /\ {\bf p}_x \to 0\, ,$$
where 
$$\mathcal{A}^{' \bullet}: 0 \to End^{filt} E \stackrel {ad(\phi)}{\longrightarrow}   End^{sfilt} E \otimes \omega_Y(D) \to 0\, .$$

\begin{theorem}\label{P6.7SWW}
(1) The dimension of the nilpotent cone $Nil'^{str}_p$ at any generic point is given by 
$${\rm dim} \ (T_u U) = {\rm dim} \ \mathbb{H}(Y, \mathcal{A}^{\bullet}) \ = \ {\rm dim} \ U_{par}(n,d)\, .
$$
(2) If $\mathcal{H}^{' str}_{par}(n,d,L)$ is smooth, then the Hitchin map $h^{str}_p$ is flat and dominant.
\end{theorem}
\begin{proof}
(1) This is proved as \cite[Proposition 6.7]{SWW}, giving the dimension of $Nil^{' str}_p$.\\
(2) The proof is similar as that of \cite[Corollary 1]{G} and  \cite[Theorem 6.8]{SWW}. 
Take any  $s \in {\bf A}'_p$. Then the closure of $\mathbb{G}_m$-orbit of $s$ contains $0$. Hence
$${\rm dim} \  (h^{str}_p)^{-1}(s) \le {\rm dim} \  (h^{str}_p)^{-1}(0) = \ {\rm dim} \ U_{par} (n,d)\, .$$
By Equation \eqref{g(Ypa)} in the proof of Proposition \ref{Ypa}, the dim $(h^{str}_p)^{-1}(s) \ge \ {\rm dim} \  U_{par}(n,d)$.  It follows that dim $(h^{str}_p)^{-1}(s) = \ {\rm dim} \  U_{par}(n,d)$ for all $s \in {\bf A}'_p$. 
The strongly parabolic Hitchin map is dominant as its fibre is nonempty (Theorem \ref{parBNRcorrA'}) over an open set. 
\end{proof}

\begin{theorem} \label{T6.9SWW}
\begin{enumerate}
\item    For the (weakly) parabolic nilpotent cone $Nil'_p$ we have 
\begin{equation} \label{nilpdim}
{\rm  dim} \ Nil'_p =  n^2(g-1) +1 + \frac{n(n-1)}{2} {\rm deg}(D)\, .
\end{equation}
\item If $\mathcal{H}'_{par}(n,d,L)$ is smooth, then the Hitchin map $h_p: \mathcal{H}'_{par}(n,d,L) \to {\bf A}'$ is flat and dominant.
\end{enumerate}
\end{theorem}

\begin{proof}
(1) This can be proved exactly as \cite[Theorem 6.9]{SWW}. One uses the fact that if the parabolic structure is a Borel structure $(B, \beta)$, then a nilpotent $\phi$ is strongly nilpotent so that $Nil'_p(B, \beta) = Nil^{' str}_p(B, \beta)$ has dimension given by Theorem \ref{P6.7SWW}(1). Then the theorem is proved by showing that there is a dominant map from a finite union of  $nil_p(B, \beta_i)$ ($\beta_i$ Borel subgroups) to $Nil'_p = Nil'_p(P,\alpha)$. \\
(2) The map $h_p$ is dominant since the general fibre is nonempty (by Theorem \ref{parabolic Hitchin fibre}).

           By Corollary \ref {coro5.6SWW} (2), the fibre $h_p^{-1}(s)$ for any $s \in {\bf A}'$, has dimension at least $n^2(g-1) +1 + \frac{n(n-1)}{2} {\rm deg}(D)$. Since dim $h_p^{-1}(s)  \le \ {\rm dim} \ h_p^{-1}(0)$, by Part (1), we have dim $h_p^{-1}(s) = \ {\rm dim} \ h_p^{-1}(0)$ for all $s$ proving flatness if $\mathcal{H}'_{par}(n,d,L)$ is smooth.    
\end{proof}

We have already proved that $h_p: N^{str}_{par}(n,d,L) \to {\bf A}_p$ is surjective. In fact, its restriction $h_p: N^{str}_{0, par} (n,d,L) \to {\bf A}_p$ is surjective (Proposition \ref{ontoA}).  The restriction is flat if $N^{str}_{0, par}(n,d,L)$ is smooth (by \cite[Theorem 6.9]{SWW}).

\subsection{Very stable parabolic bundles and the restriction of the Hitchin map} \hfill
 
\begin{definition}  \label{VstrE}
         For a quasi-parabolic sheaf $E$, define  
$$V^p_E \ := \ {\bf H}^0(Y, ParEnd(E)\otimes \omega_Y(D))\, ,$$
$$V^{str}_E \ := \ {\bf H}^0(Y, SParEnd(E)\otimes \omega_Y(D))\, ,$$

    For a quasi-parabolic torsion-free sheaf $E$,  a section $\phi \in H^0(Y, ParEnd(E)\otimes \omega_Y(D))$
  is called nilpotent if $(E, \phi) \in Nil_p$. 
  
    For a parabolic torsion-free sheaf $E_*$,  a section $\phi \in H^0(Y, SParEnd(E)\otimes \omega_Y(D))$
  is called strongly nilpotent if $(E_*, \phi) \in Nil^{str}_p$. 

Let $V^{nilp}_E$ (respectively $V^{str nilp}_E$) be the sets of nilpotent (respectively strongly nilpotent) Higgs fields in $V_E$. 
\end{definition} 

\begin{definition}

A quasi-parabolic torsion-free sheaf $E$ is called very stable (respectively strongly very stable) if it has no non-zero nilpotent (respectively strongly nilpotent) section.

\end{definition}

In case $Y$ is smooth, this definition has appeared in \cite{Pe}.  
Note that a very stable parabolic torsion-free sheaf is strongly very stable. 

\begin{lemma} \label{parvstable}
         Let $E_*$ be a parabolic torsion-free sheaf (with fixed parabolic weight $\alpha$). Assume that $g \ge d(D)+ 1$.

 If $E_*$ is not parabolic stable, then there is a nonzero $\phi \in H^0(Y, SParEnd(E_*)\otimes \omega_Y(D))$ such that  $\phi^r =0$ for some positive integer $r$. 
 
 If $E_*$ is a very stable parabolic vector bundle, then $E_*$ is parabolic stable (i.e., $\alpha$-stable). \\

\end{lemma}

\begin{proof}
       If  the parabolic sheaf $E_*$ is not parabolic stable, then there is an exact sequence of parabolic torsion-free sheaves 
\begin{equation} \label{vs0}        
        0 \longrightarrow N_* \stackrel{i}{\longrightarrow} E_* \stackrel{q}{\longrightarrow} Q_* \longrightarrow 0\, ,
\end{equation}
with $p\mu(N) \ \ge \ p\mu(Q)$. 
Let
$$F_* := SParHom(Q_*, N_*\otimes \omega_Y(D))\, .$$
 A nonzero element in $H^0(Y, F_*)$  can be considered as a nonzero strongly parabolic homomorphism $u: Q_* \to N_*\otimes \omega_Y(D)$.  Then $\phi: E_* \to E_*\otimes \omega_Y(D)$ defined by $\phi = (i \otimes Id_{\omega_Y}) \circ u \circ q$  is a nonzero strongly parabolic nilpotent Hitchin field. If $E$ is a vector bundle, this will imply that $(E_*, \phi) \in Nil_p$ so that $E_*$ is not a very stable parabolic bundle. Hence any very stable parabolic bundle is parabolic stable. If $E$ is not a vector bundle, then $(E_*, \phi)$ thus constructed may not be in $N_{par}(n,d,L)$.

   We now check that for $g \ge d(D)+1$, we have $h^0(Y, F_*) > 0 \, .$ 
      We keep the notation of \cite{BhBi}, see \cite[p. 501]{BhBi} for more details. The set of weights of $E$ at $x$ is the union of  the sets of weights of $N$ and $Q$ at $x$. For a weight $\alpha_i(x)$ in this we write $n^N_i(x)$ for the multiplicity of  $\alpha_i(x)$ in the fibre of $N$ at $x$, we have the convention that $n^N_i(x) =0$ if $\alpha_i$ is not a weight of $N$ at $x$, similarly we have $n^Q_i(x)$.

    The sum $\sum_{i \ge \ell} n^N_i(x) n^Q_{ \ell}(x) \le \ \sum_{i, \ell}  n^N_i(x) n^Q_{ \ell}(x) = r(N) r(Q)$. Hence  
    $$\sum_{x\in I} \ \sum_{i \ge \ell} n^N_i(x) n^Q_{ \ell}(x) \le r(N) r(Q) d(D)\, .$$
Note that since the weights belong to the interval $[0, 1)$, at any $x \in D$, the weight of $N$ at $x$ satisfies $0 \le wt_x(N) < r(N),$ similarly for $Q$, so that $ \frac{wt(Q)}{r(Q)} \ge 0$ and $d(D) - \frac{wt(N)}{r(N)} > 0$. 
Hence 
$$- \frac{wt N}{r(N)} + \frac{wt(Q)}{r(Q)} +d(D) >0\, .$$           
It follows from \cite[Equations (2.1) and (2.2)]{BhBi} that $\chi(F_*) = \chi(F) - \sum_{x\in I} \ \sum_{i \ge \ell} n^N_i(x) n^Q_{ \ell}(x)$. Hence 
$$\chi(F_*) = r(Q) d(N \otimes \omega(D)) - r(N) d(Q) +r(Q)r(N) (1- g) - \sum_{x\in I} \ \sum_{i \ge \ell} n^N_i(x) n^Q_{ \ell}(x) + \sum_j b_j(N) b_j(Q)\, ,$$ 
where  $b_j(N), b_j(Q)$ denote the local types of $N$ and $Q$. Hence we have
$$
\begin{array}{lll}
h^0(Y, F_*) & \ge & \chi(F_*) \\
    {}             & = & r(Q) r(N) (\mu(N) - \mu(Q) + d(D) +g -1) \\
    {} & {} & - \sum_{x\in I} \ \sum_{i \ge \ell} n^N_i(x) n^Q_{ \ell}(x) + \sum_j b_j(N) b_j(Q)\\
    {}     & \ge & r(Q)r(N)(p\mu(N_*) - p\mu(Q_*) - \frac{wt(N)}{r(N)} + \frac{wt(Q)}{r(Q)} + d(D) +g -1) \\
   {} & {} & - \sum_{x\in I} \ r(Q) r(N) + \sum_j b_j(N) b_j(Q) \\
   {}    & > & r(Q) r(N)(p\mu(N_*) - p\mu(Q_*) + g-1) \\
   {} & {} & - r(Q) r(N) d(D) + \sum_j b_j(N) b_j(Q)\\
   {} & {} \ge & r(Q) r(N) (g -d(D) -1) +  \sum_j b_j(N) b_j(Q)\, .
\end{array}
$$ 
Hence $h^0(Y, F_*) > 0$ if $g \ge d(D)+1$.  \\
 
\end{proof}

\begin{remark} 
       In \cite[Lemma 3.5]{Pe}, Lemma \ref{parvstable} is proved for a nonsingular curve $Y$ of genus $g \ge 2$. The proof uses many results on orbifold bundles and root stacks.  
\end{remark}

       Note that $\mathcal{H}^{' str}_{par}(n,d,L)$ is an open dense subset of $(h^{str}_{par})^{-1} {\bf A'}_p$ and $(h^{str}_{par})^{-1} {\bf A'}_p$ is a closed set as the Hitchin map is proper.

\begin{theorem} \label{verystableexist}
      In the moduli space of stable parabolic sheaves $U_{par}(n,d)$, the subset consisting of very stable sheaves contains a nonempty (Zariski) open subset $S$.  

\end{theorem} 

\begin{proof}
       
 Let $\mathcal{H}^{ps} \subset (h^{str}_{par})^{-1} {\bf A'}_p$ be the open dense subset consisting of parabolic Hitchin sheaves with the underlying parabolic sheaf stable. Then there is a morphism ${Fg}: \mathcal{H}^{ps} \to  U^s_{par}(n,d)$ defined by $(E_*, \phi) \mapsto E_*$.  Let $Z_1 \subset \mathcal{H}^{ps}$ be the subset consisting of $(E_*, \phi)$ with $\phi$ nilpotent and nonzero. There is a $\mathbb{G}_m$ action on $\mathcal{H}^{ps}$ given by $t (E_*, \phi) = (E_*, t \phi)$, $Z_1$ is  $\mathbb{G}_m$-invariant and $Fg$ is $\mathbb{G}_m$-equivariant.     

       Since $Z_1 \subset Nil^{str}_p$, all non very stable sheaves in $U^s_{par}(n,d)$ are contained in $Fg(Z_1)$. For $(E_*, \phi) \in Z_1$, $E_*$ is stable and $\phi$ is nonzero, hence $\mathbb{G}_m$ acts freely on $Z_1$. Thus 
$${\rm dim} \ {Fg} (Z_1) \le {\rm dim} \ Z_1 -1 \le \ {\rm dim} \ Nil^{str}_p -1 = {\rm dim} \ U^s_{par}(n,d) -1\, .$$ 
It follows that the set of very stable parabolic sheaves contains a nonempty open set. 

\end{proof}
In case $Y$ is nonsingular,  the result of Theorem \ref{verystableexist} was proved in \cite[Theorem 6.14]{SWW}, our proof uses a similar strategy as in \cite{SWW}.
 
\begin{corollary} \label{forgetfuldominant}
           For a general $a \in {\bf A'}_p$, the restriction of the forgetful map 
           $${Fg}_a: \ (h^{str}_{par})^{-1} (a) \dasharrow U^s_{par}(n,d)$$ 
is a dominant rational map.  
\end{corollary}

\begin{proof}
   Since the subset $U^{' s}_{par}(n,d)$ is an open dense subset of $U^s_{par}(n,d)$, it suffices to check that  for a general $a \in {\bf A'}_p$, the restriction of the forgetful map ${Fg}_a : (h^{str}_{p})^{-1} (a) \cap \mathcal{H}^{' str}_{par}(n,d) \dasharrow U^{' s}_{par}(n,d)$ is a dominant rational map.
   
       Consider the composite map $\rho = {Fg}_a \times h^{str}_p$ i.e.,  
 $$\rho:  \mathcal{H}^{' str}_{par}(n,d) \dasharrow U^{' s}_{par}(n,d) \times {\bf A'}_p\, ; (E_*, \phi) \mapsto (E_*, h^{str}_p(E_*, \phi))\, .$$    
Since dim $\mathcal{H}^{' str}_{par}(n,d) =$ dim $U^{' s}_{par}(n,d) + $ dim ${\bf A'}_p$, the morphism $\rho$ is generically finite. By Theorem \ref{verystableexist}, there is an open subset $S \subset U^{' s}_{par}(n,d)$ such that for $\rho^{-1}(E_*, 0) = (E_*, 0)$ for $(E_*,0) \in S \times \{0\}$.  Hence $ 
(h^{str}_{p})^{-1} (a) \cap \mathcal{H}^{' str}_{par}(n,d) \dasharrow  U'_{par}(n,d)$ is dominant for a general $a \in {\bf A'}_p$.    
\end{proof}

\begin{lemma} \label{Fgcodim2}
        The rational forgetful map ${Fg}_a$ is defined on a big open subset of  $(h^{str}_{p})^{-1} (a) \cap \mathcal{H}^{' str}_{par}(n,d)$ i.e. the open subset where it is defined has complement of codimension at least $2$.
\end{lemma}
\begin{proof}
     This can be proved using the method similar to that in \cite{VLP} , for example as in \cite[Corollary 4.8]{WW}.
\end{proof}

    We remark that there is a (non-canonically defined) parabolic theta line bundle $\mathcal{L}_p$ over $U^s_{par}(n,d)$. In view of Lemma \ref{Fgcodim2}, $({Fg}_a)^* \mathcal{L}_p$ can be extended  to a line bundle over $(h^{str}_{p})^{-1} (a) \cap \mathcal{H}^{' str}_{par}(n,d)$.  

\begin{lemma} \label{dimVstrE}
For $E_*$ stable, 
$$\ {\rm dim} \ V^{str}_E \ = \ {\rm dim} \ U_{par}(n,d)\, .$$
\end{lemma}
\begin{proof}
We have $V^{str}_E = H^0(Y, SParEnd(E)\otimes \omega_Y(D))$. 
Since $E_*$ is stable, $h^0(Y, ParEnd(E)) = 1$. By Serre duality, this gives 
$h^1(Y, SParEnd(E)\otimes \omega_Y(D)) = 1$. From the exact sequence 
$$0 \to SParEnd(E)\otimes \omega_Y(D) \to  End(E)\otimes \omega_Y(D) \to \times_{x \in D} {\bf p}_x \to 0\, ,$$ we have $\chi(SParEnd(E) \otimes \omega_Y(D)) = \chi(End(E)\otimes \omega_Y(D)) - \sum_{x\in D} \ {\rm dim} \ {P_x}$. Therefore,
$$
\begin{array}{lll}
h^0(Y, SParEnd(E)\otimes \omega_Y(D)) & = & 1 + \chi(SParEnd(E) \otimes \omega_Y(D))\\
{} & = & 1 + \chi(End(E)\otimes \omega_Y(D)) - \sum_{x\in D} \ {\rm dim} \ {P_x} \\
{} & = & 1+ n^2(2g-2 + d(D)) + n^2(1-g) - \sum_{x\in D} \ {\rm dim} \ {P_x} \\                               
{} & = & 1 + n^2(g-1)  + \sum_{x\in D} \ (n^2  - {\rm dim} \  P_x) \\
{} & = & n^2(g-1)+1 + \sum_{x \in D} \ {\rm dim} \ (GL(n)/ P_x)  \\
{} & = & {\rm dim} \ U_{par}(n,d)
\end{array}{ll}
$$
\end{proof}

\begin{theorem}
Let $E_*$ be a stable parabolic bundle. Let 
$$h^{str}_{E} = h^{str}_p\vert_{V^{str}_E}: V^{str}_E \to {\bf A'}_p$$
 be the restriction of the strongly parabolic Hitchin morphism $h^{str}_p$  ($V^{str}_E$ is defined in Definition \ref{VstrE}). Then the following are equivalent:\\
(1) E is strongly very stable.\\
(2) $h^{str}_E$ is finite.\\
(3) $h^{str}_E$ is quasi-finite.\\
(4) $h^{str}_E$ is proper. \\
(5) The inclusion map $i^{str}_V: V^{str}_E \rightarrow N_p^{str}(n,d,\omega_Y)$ is proper. 
\end{theorem}
\begin{proof}
$(1) \Rightarrow (2)$ : Since both $V^{str}_E$ and ${\bf A'}_p$ are affine spaces and very stability of $E$ means  that $(h^{str}_{E})^{-1} (0) = 0$, it follows that $h^{str}_E$ is finite (\cite[Lemma 3.1]{Pe}, 
\cite[Lemma 1.3]{Z}). \\
$(2) \Rightarrow (3), (4)$ : Any finite map is quasi-finite and proper. \\
$(3) \Rightarrow (1)$ : Enough to show that if $E$ not strongly very stable, then  $h^{str}_E$ is not quasi-finite. 
If $E$ is not strongly very stable, then $E$ has a strongly nilpotent vector field $\phi$. Since $E_*$ is stable, the $\mathbb{C}$-orbit of $\phi$ gives a positive dimension fibre of $h^{str}_E$, contradicting the quasi-finiteness of $h^{str}_E$.\\
$(4) \Rightarrow (3)$ and $(2)$: By Lemma \ref{SWWL3.3}, dim ${\bf A'}_p =$ dim $U_{par}(n,d)$. Then by Lemma \ref{dimVstrE}, we have dim $V^{str}_E = $ dim ${\bf A'}_p$. Thus $h^{str}_E$ is a proper map of affine spaces of same dimension, hence it is quasi-finite and therefore finite. \\
$(4) \Rightarrow (5)$ : One has $h^{str}_E = h^{str}_p \circ i^{str}_V$.  If $(4)$ holds i.e., the composite $h^{str}_p \circ i^{str}_V$ is proper, then $i^{str}_V$ is proper.\\
$(5) \Rightarrow (4)$ because the composite of two proper maps is proper. 
\end{proof}

\section{The fixed determinant case} \label{fixdet} \hfill

        Fix a line bundle $\xi$ of degree $d$ on the nodal curve $Y$. We assume that the parabolic weights are general so that semistability is same as stability for parabolic (Hitchin) bundles.

\subsection{The moduli space of parabolic bundles with a fixed determinant} \label{Thetapar} \hfill 

Let $U'_{par}(n, \xi)$ denote the moduli space of semistable parabolic vector bundles of rank $n$ with a fixed determinant $\xi$ on $Y$.  Let $U_{par}(n, \xi)$ be the closure of $U'_{par}(n, \xi)$ in $U_{par}(n, d)$. One has   
\begin{equation}
 {\rm dim} \ U_{par}(n,\xi) = (n^2 - 1)(g-1) + \sum_{x \in D} \ {\rm dim} \ (SL(n)/P_x)\, ,
\end{equation}
where $P_x$ is the parabolic subgroup of $SL(n)$ determined by the parabolic structure at $x$. 

For positive integers $K, \ell$ and non negative integers $b_x, x \in D$ satisfying 
\begin{equation}  \label{*}   
\sum_{x \in d} \sum_{i= 1}^{l_x} K (\alpha_{i+1}(x) - \alpha_i(x)) r_i(x) + n \sum_{x\in D} \ b_x + n \ell = K \chi\, ,     
\end{equation}       
where $\chi$ denotes the Euler characteristics of $E$ i.e. $\chi = d+n(1-g)$, there is an ample line bundle $\Theta_{par}$ on $U_{par}(n,\xi)$ having the following universal property. Let $\mathcal{E}_*$ be a family of parabolic vector bundles on $Y$ parametrised by a scheme $T$ and $\pi_T: T \times Y \to T$ the projection. 
    
        Define 
$$\theta_T := (Det R_{\pi_T} \mathcal{E} )^K \otimes \bigotimes_{x \in D} [det  (\mathcal{E}_{T \times x})^{b_x} \otimes \bigotimes_{i = 1}^{l_x}  (det  (\mathcal{E}_{T\times x}/F'_i(\mathcal{E}_{T\times x}))) ] \otimes det (\mathcal{E}_{T\times t})^{\ell} \, .$$

Then there is a morphism $f_T : T \to  U_{par}(n,\xi)$ such that 
$$\theta_T = f_T^*(\Theta_{par})\, .$$

\subsection{The Hitchin base and the Hitchin map}  \label{Apxi}  \hfill

       We have the fixed determinant moduli space ${\mathcal H}^{str}_{par}(n, \xi, L) \subset  {\mathcal H}^{str}_{par}(n, d, L)$ defined as the closure of ${\mathcal H}^{' str}_{par}(n, \xi, L)$, the moduli space of strongly parabolic Hitchin bundles of rank $n$ and a fixed determinant $\xi$ with the Hitchin field $\phi$ having values in $L(D), D$ being the parabolic divisor. For $(E_*, \phi) \in {\mathcal H}^{str}_{par}(n, \xi, L)$, one has $tr \phi = 0$. Hence the parabolic Hitchin base in this case is 
$${\bf A}_{p, \xi} := \prod_{j=2}^n \  {\bf H}^0(Y, L^j \otimes \mathcal{O}_Y(\sum_{x \in D} (j - \gamma_j(x) x) )\otimes p_*\mathcal{O}_X\, .$$     
Let
$${\bf A}'_{p, \xi} := \prod_{j=2}^n \  {\bf H}^0(Y, L^j \otimes \mathcal{O}_Y(\sum_{x \in D} (j - \gamma_j(x) x) )\subset {\bf A}_{p, \xi}\, .$$
We remark that both ${\bf A}_{p, \xi}$ and ${\bf A'}_{p, \xi}$ are independent of $\xi$.  
\begin{lemma} \label{A'pxi}
For $L = \omega_Y$, 
$${\rm dim} \  {\bf A'}_{p, \xi} \ = \ {\rm dim} \ U'_{par}(n,\xi)\, .$$
\end{lemma}
\begin{proof}
One has  ${\rm dim} \  {\bf A'}_{p, \xi} \ = \ {\rm dim} {\bf A'}_p - h^0(Y, \omega_Y(\sum_{x\in D}  (1- \gamma_1(x)) x ) )\, .$ Note that $\gamma_1(x) = 1$, so that 
$$
\begin{array}{lll}
{\rm dim} \  {\bf A'}_{p, \xi} & = &  {\rm dim} \ {\bf A'}_p - h^0(Y, \omega_Y) \\
                       {}              & =  &   {\rm dim} \ {\bf A'}_p - g\\
                       {}                & = &  (n^2 -1) (g-1) + \sum_{x\in D} \  {\rm dim} \ G/P_x\\
                       {}                 & = & \ {\rm dim} \ U'_{par} (n, \xi) 
\end{array}
$$ 
\end{proof}

Let $N^{str}_{par}(n,\xi,L) \ = \ N_{par}(n,d,L) \cap \mathcal{H}^{str}_{par}(n,\xi,L)$ and let 
$$h^{str}_{p, \xi} :  N^{str}_{par}(n, \xi, L) \longrightarrow   {\bf A}_{p, \xi} $$ 
be the parabolic Hitchin map. There are commutative diagrams 
$$
\begin{array}{ccc}
N^{str}_{par}(n, \xi, L) & \stackrel{h^{str}_{p, \xi}}{\longrightarrow} &  {\bf A}_{p, \xi} \\
{} & \stackrel{h^{str}_p}{\searrow} & \downarrow\\
{} &  {} &  {\bf A}_{p} 
\end{array}
$$ 
and 
$$
\begin{array}{ccc}
{\mathcal H}^{' str}_{par}(n, \xi, L) & \stackrel{h^{' str}_{p, \xi}}{\longrightarrow} &  {\bf A'}_{p, \xi} \\
{} & {\stackrel{h^{' str}_p}{\searrow}} & \downarrow  \\
{} & {} &  {\bf A'}_{p} 
\end{array}
$$ 

Since $N^{str}_{par}(n, \xi, L)$ is a closed subset of $N^{str}_{par}(n, d, L)$, the restriction of $h^{str}_p$ to $N^{str}_{par}(n, \xi, L)$ is a proper map. This restriction factors through $h^{str}_{p, \xi}$, hence $h^{str}_{p, \xi}$ is proper. 

\subsection{The Prym variety $P_{\delta}$} \label{Pdelta}   \hfill
 
      Assume that $L= \omega_Y$.
     
     It follows from the parabolic BNR correspondence (see Theorem \ref{parBNRcorrA'}) that, for a general $a \in {\bf A}'_p$  the fibre of the strongly parabolic Hitchin map $(h^{' str}_p)^{-1} (a)$ is isomorphic to 
$J^{\delta} (Y^p_a)$. 
   Recall that on $\bar{J}(Y^p_a) \times Y^p_a$, we have a universal family $\mathcal{N}$ of torsion-free sheaves of rank $1$ and degree $\delta$. Taking the direct image under the morphism $Id \times \pi_p$, we get a family of Hitchin sheaves of rank $n$ and degree $d$ on $Y$ parametrised by $\bar{J}(Y^p_a)$ giving an injective morphism $\psi: \bar{J}(Y^p_a) \to  \mathcal{H}^{str}_{par}(n,d, \omega_Y)$.  One has $(h^{str}_{p})^{-1}(a) = \psi(\bar{J}(Y^p_a)) \cap N^{str}_{par}(n,d,L)$. Moreover $\psi$ restricts to the isomorphism 
 $J^{\delta} (Y^p_a) \cong (h^{' str}_p)^{-1} (a)$. 

             For a general $a \in {\bf A}^{' \xi}_p$, the fibre $(h^{str}_{p, \xi})^{-1} (a) = (h^{str}_{p})^{-1}(a) \cap \mathcal{H}^{str}_{par}(n,\xi, \omega_Y)\, .$  
 
   Define 
 $$P_{\delta}: =  \{ N \in  \overline{J}^{\delta} (Y^p_a) \ \vert \ \psi(N) \in  (h^{str}_{p, \xi})^{-1}(a) \, .\}$$
 Then $P_{\delta}$ contains an open subset $P'_{\delta}$ consisting of $N \in J^{\delta} (Y^p_a)$ such that $\psi(N) \in  (h^{str}_{p})^{-1}(a) \cap \mathcal{H}^{' str}_{par}(n,\xi, \omega_Y)\, .$ 
   
                 Thus for a general $a \in {\bf A}^{' \xi}_p$, there is an injective morphism from the Prym variety $P_{\delta}$ to $\mathcal{H}^{str}_{par}(n,\xi, \omega_Y)$, whose restriction to $P'_{\delta}$ is an isomorphism onto $(h^{str}_{p, \xi})^{-1} (a) \cap \mathcal{H}^{' str}_{par}(n,\xi, \omega_Y)\, .$

   Assume that $Y$ is defined over an algebraically closed field of characteristic $0$. Then  $P'_{\delta}$ is smooth.
  Let $\widetilde{P}_{\delta}$ be the pull back of $P_{\delta}$ to $\tilde{J}(Y^p_a)$, the desingularisation (normalisation) of  $\bar{J}(Y^p_a)$. Then $\widetilde{P}_{\delta}$ is smooth and it has an open subset mapping isomorphically onto $P'_{\delta}$.  We identify this open subset with $P'_{\delta}$.                  

Given $z \in Y \setminus D$, we may assume that $\pi_p$ is unramified over $z$.

\begin{lemma} \label{pdeltacodim2}
For $g \ge 2$, 
$${\rm Codim} \ (P_{\delta} - P'_{\delta}, P_{\delta}) \ge 2\, .$$  
\end{lemma}

\begin{proof}
 For each node $y_j$ of $Y$, let $\{y'_{j, i} \}, i= 1, \cdots ,n$ be the nodes of $Y^p_a$ lying over the node $y_j$. Let $\tilde{p}: \widetilde{Y^p_a} \to Y^p_a$ be the normalisation and let $x'_{j, i}, z'_{j, i}$ be the points of $\widetilde{Y^p_a}$ lying over $y'_{j, i}$. For $N \in \bar{J}^{\delta}(Y^p_a)$, let  $\psi(N) = (E_*, \phi)$ where the underlying sheaf $E = \pi_{p *} N$. Since $\pi_p$ is a finite map, $R^1\pi_{p *} N =0$ and  the fibre $E_{y_j} = \oplus_i N_{y_{j, i}} = k^{a_j} \oplus k^{2 b_j}, a_j + b_j = n$, here $b_j$ is the number of $y_{j, i}, i = 1, \cdots, n$ where $N$ is not locally free. It follows that the stalk $E_{(y_j)} \cong   \mathcal{O}_{(y_j)}^{a_j} \oplus m_{(y_j)}^{b_j}$ i.e., $E$ is of local type $b_j$ at $y_j$. 
 
       If $E$ is of local type $b_j = 1$ at all $y_j$ where $E$ is not locally free (or equivalently $N$ is not locally free), then $E$ has a determinant isomorphic to a nonlocally free sheaf and hence $(E_*, \phi) \notin  \mathcal{H}^{str}_{par}(n,\xi, \omega_Y)$.  Hence $N \in P_{\delta} \setminus P'_{\delta}$ if and only if $N$ is not locally free at two or more nodes lying over $y_j$. 
 
         We claim that the set $S = \{N \in P_{\delta}$ such that $N$ is not locally free at minimum two nodes lying over some $y_j$\} has codimension at least $2$.  One has a normalisation $\tilde{J}(Y^p_a)$ of $\bar{J}(Y^p_a)$ which is a fibre bundle $\pi_p: \tilde{J}(Y^p_a) \to J(\widetilde{Y^p_a})$.  Let $\widetilde{P_{\delta}}$ (resp. $\widetilde{P'_{\delta}}$) be the inverse image of $P_{\delta}$ (resp. $P'_{\delta}$) in $\tilde{J}(Y^p_a)$.  Then for $(N', F(N'))  \in \widetilde{P_{\delta}}$, $N'$ is a line bundle belonging to $\pi^p(\widetilde{P'_{\delta}})$ and $F(N') = (F_{j, i}(N'))_{j, i}, F_{j, i}(N') \subset N'_{x'_{j, i}} \oplus N'_{z'_{j, i}}$ is a $1$-dimensional subspace. Thus the fibre over $N'$ is contained in a product of $\mathbb{P}^1$s and contains a product of $\mathbb{C}^*$ contained in each of those $\mathbb{P}^1$s, the elements in $\widetilde{P'_{\delta}}$ correspond to the latter. Suppose that $(N' , F(N'))$ maps to $N \in S$. Then $F_i(N') = N'_{x_{j, i}}$ or $F_i(N') = N'_{z_{j, i}}$ for at least two $i$, so that $F(N')$ belongs to a subspace of codimension at least $2$ in the fibre. Since such $(N', F(N'))$ are precisely those which map to $N$ which are not locally free at two of the $y_{j, i}$s, it follows  that 
codim $(P_{\delta} \setminus P'_{\delta}, P_{\delta}) \ge 2)$.   
\end{proof}

           Let $N' (= m' n)$ be the number of $\{y'_{j, i} \}, i= 1, \cdots ,n$ with $j$ varying over nodes of $Y$. The singular  set of $P_{\delta} \subset \bar{J}(Y^p_a)$ is a union of $N'$ number of irreducible closed subsets $S_{j,i}$. The inverse image of this singular set in $\tilde{P}_{\delta}$ is the disjoint union of $2 N'$ irreducible closed subsets $S^1_{j,i}$ and $S^2_{j,i}$ and $S_{j,i}$ is obtained by a simple identification of $S^1_{j,i}$ and $S^2_{j,i}$.

\subsection{Nilpotent cone in the fixed determinant case} \hfill
         
          For a filtered vector bundle $E$, let $End^{ 0, filt}(E)$ be the sheaf of local trace zero endomorphisms of $E$  which preserve the filtration. For a parabolic vector bundle $E_*$, let $ParEnd^{ 0}(E)$ and $SParEnd^0(E)$ denote respectively the sheaf of parabolic endomorphisms with trace zero and the sheaf of strongly parabolic endomorphisms with trace zero.   

Define
      $$Nil^{str}_{p, \xi} := (h^{str}_{p, \xi})^{-1}(0) \subset \mathcal{N}^{str}_{par}(n,\xi,L)\, , Nil^{' str}_{p, \xi} = (h^{str}_{p, \xi})^{-1}(0) \cap \mathcal{H}'_{par}(n,\xi,L)\, .$$
As in the non-fixed determinant case, using Lemma \ref{A'pxi} we see that dim $(h^{str}_{p, \xi})^{-1}(0) \ge \ {\rm dim} \  U_{par}(n,\xi)$.
  
  The dimension of $Nil^{' str}_{p, \xi}$ can be computed using infinitesimal deformations as in the non-fixed determinant case (see  Subsection \ref{dimNilstrp} for details). In the notations of Subsection \ref{dimNilstrp}, restricting $q$ (defined in the equation \eqref{q}) to $q^{-1}Nil^{' str}_{p, \xi}$ we get a universal family $\mathbb{U}_{N, \xi}:= (V \otimes_k \mathcal{O}_{Y \times q^{-1} Nil^{' str}_{p, \xi}} \to \mathcal{E}_*, \Phi)$ satisfying the following universal property. For a scheme $S$ and  a flat family of quotients $(V \otimes_k \mathcal{O}_{Y \times  S} \to \mathcal{E}_S, \phi_S)$ of strongly parabolic Hitchin bundles with a fixed determinant $\xi$ and satisfying $h^{str}_{p, \xi}(\mathcal{E}_S, \phi_S) = 0$, there is a morphism $f_S: S \to Nil^{' str}_{p, \xi}$ such that $(Id \times f_S)^*(\mathbb{U}_{N, \xi}) \cong (V \otimes_K \mathcal{O}_{Y \times  S} \to \mathcal{E}_S, \phi_S)$. 
          
          It suffices to determine the dimension of each irreducible component of  $Nil^{' str}_{p, \xi}$ with reduced structure.  Let $\eta$ be the generic point of an irreducible component of $Nil^{' str}_{p, \xi}$ and let $\phi_{\eta}: \Phi\vert_{q^{-1} (\eta)^{red}}$ be the restriction. Then $\phi_{\eta}$ induces the  canonical filtration $\{Ker(\phi^i_{\eta})\}$ of  $\mathcal{E} \vert_{q^{-1} (\eta)^{red}}$. There is an open subset $U_0$ of $q_{\xi}^{-1} (\eta)^{red}$ with generic point $w$ such that the graded terms in the canonical filtration are vector bundles.   

\begin{theorem} \label{dimnullcone}
\begin{enumerate}
\item The space of infinitesimal deformations in $U_0$ of a strongly nilpotent Hitchin bundle $(E_*,\phi)$ with $det E = \xi$, is canonically isomorphic to the hypercohomology $\mathbb{H}^1(Y, \mathcal{A}^{\bullet}_{\xi})$ of the following complex of sheaves on $Y$ 
\begin{equation} \label{Axibullet}
\mathcal{A}^{\bullet}_{\xi}: 0 \longrightarrow ParEnd^0(E) \cap End^{ 0, filt} E \stackrel{ad \phi}{\longrightarrow} 
SParEnd ^0(E) \cap End ^{ 0, s-filt} E \otimes \omega_Y(D) \rightarrow 0
\end{equation}
\item The dimension of  $Nil^{' str}_{p, \xi}$ is given by 
$${\rm dim}_k \ (T_w U_0) = \ {\rm dim}_k \ \mathbb{H}^1(Y, \mathcal{A}^{\bullet}_{\xi}) = \ {\rm dim} \  U_{par}(n,\xi)\, .$$ 
\item The Hitchin map $h^{str}_{p, \xi}$ is flat and surjective if $\mathcal{H}^{' str}_{par}(n,\xi, \omega_Y)$ is smooth. 
\end{enumerate}
\end{theorem} 
\begin{proof}
(1) This can be proved as in Theorem \ref{deformations} or \cite[Theorem 6.4]{SWW}, replacing endomorphisms by trace zero endomorphisms.\\
(2) We prove Part (2) on the lines of the proof of Theorem \ref{P6.7SWW} replacing $GL(n)$ by $SL(n)$, we briefly sketch the proof.  An $SL(n)$-bundle $E'$, can be interpreted as a vector bundle $E$ with determinant the trivial line bundle $I_1$ (more generally a fixed line bundle) together with an isomorphism $det E \cong I_1$. Let ${\bf g}, {\bf p}$ and ${\bf n}$ be the Lie algebras of  $SL(n)$,  a parabolic subgroup $P$  and  its  nilpotent subgroup $N$ respectively. 

      The sequence \eqref{Axibullet} gives an exact sequence 
$$0 \to \mathcal{A}^{\bullet}_{\xi} \to \mathcal{A}^{' \bullet}_{\xi} \to {i_x}_* \ {\bf p} /\ {\bf p}_x \to 0\, ,$$
where 
$$\mathcal{A}^{' \bullet}_{\xi}: 0 \to End^{ 0 filt} E \stackrel {ad(\phi)}{\longrightarrow}   End^{ 0 \ s-filt} E \otimes \omega_Y(D) \to 0\, .$$

    Let $E'_P$ be the reduction of the structure group of $E'$ to the parabolic subgroup $P$. 
One has $ad E'_p({\bf g}) = End^0(E), ad E'_p({\bf p}) = End^{ 0, filt}(E), ad E'_p({\bf n}) = End^{ 0, s-filt}(E)$.  One has $\check{\bf n} \cong {\bf g}/ {\bf p}, \chi(End^0(E)) = \chi(End(E) - \chi(I_1) = (n^2 -1)(1- g)$. 
$$,
\begin{array}{ll}
\chi(\mathcal{A}^{' \bullet}_{\xi})  & =  \chi(End^{ 0, filt}(E)) - \chi(End^{ 0, s-filt}(E) \otimes \omega_Y(D))\\
{}  & = \chi(ad E'_p({\bf p})) - \chi(ad E'_p({\bf n}) \otimes \omega_Y(D))\\
{} & = \chi(ad E'_p({\bf p})) - \chi(ad E'_p({\bf n}) \otimes \omega_Y) - d(D) \ {\rm dim} \ {\bf n}\\
{} & = \chi(ad E'_p({\bf p})) + \chi(ad E'_p(\check{\bf n})) - d(D) \ {\rm dim} \ {\bf n}\\
{} & = \chi(ad E'_p({\bf g}))  - d(D) \ {\rm dim} \ {\bf n}\\
{} & = (n^2 -1)(1- g) - d(D) \ {\rm dim} \ {\bf n}
\end{array}
$$ 
Hence 
$$\chi(\mathcal{A}^{\bullet}_{\xi}) = \chi(\mathcal{A}^{' \bullet}_{\xi})  - \sum_{x \in D} \ {\rm dim} \ {\bf p}/ {\bf p}_x 
= (n^2 -1)(1- g) -  \sum_{x\in D} \ {\rm dim} \ {\bf g}/{\bf p}_x\, .$$

The stability of $E_*$ implies that $H^0(Y, End^0(E_*)) = 0$ so that $h^0(Y, \mathcal{A}^{\bullet}_{\xi}) = 0$. One has $h^2(Y, \mathcal{A}^{\bullet}_{\xi}) = 0$ (see the proof of \cite[Proposition 6.7]{SWW}). Hence 
$${\rm dim}_k \ (T_w U_0) =  \ {\rm dim}_k \ \mathbb{H}^1(Y, \mathcal{A}^{\bullet}_{\xi}) = (n^2 -1)(g - 1) + \sum_{x\in D} \ {\rm dim} \ {\bf g}/{\bf p}_x = \ {\rm dim} \ U_{par}(n, \xi)\, .$$  

(3) This can be proved as in Theorem \ref{P6.7SWW} using Parts (1) and (2), replacing $J(Y^p_a)$ by $P'_{\delta}$. 
\end{proof}

Let $Fg^{\xi}: \mathcal{H}^{str}_{par}(n,\xi,  \omega_Y) \dasharrow U_{par}(n,\xi)$ be the surjective forgetful rational map which forgets the parabolic structure. 

Let $\mathbb{H}(n, \xi)$ be the closure of $\mathcal{H}^{' str}_{par}(n,\xi,  \omega_Y)$ in $\mathcal{H}^{str}_{par}(n,\xi,  \omega_Y)$. Then for a general $a \in {\bf A'}_{p, \xi}$, the Hitchin fibre $(h^{str}_{p, \xi})^{-1} (a)$ is contained in $\mathbb{H}(n, \xi)$.

\begin{corollary} \label{forgetfuldominant}
           For a general $a \in {\bf A'}_{p, \xi}$, the restriction of the forgetful map $Fg^{\xi}$ to 
 $(h^{str}_{p, \xi})^{-1} (a)$,
 $$Fg^{\xi}_a: (h^{str}_{p, \xi})^{-1} (a) \dasharrow U_{par}(n, \xi)$$ 
 is a dominant rational map.  
\end{corollary}

\begin{proof}   
       Consider the composite map $\rho^{\xi} = Fg^{\xi}_a \times h^{str}_{p, \xi}$ i.e.,  
 $$\rho^{\xi}: H^{str}_{par}(n,\xi) \dasharrow U_{par}(n,\xi) \times {\bf A'}_{p, \xi}\, ; (E_*, \phi) \mapsto (E_*, h^{str}_{p, \xi}(E_*, \phi))\, .$$    
Note that $\rho^{\xi}$ is a dominant map. Since dim $\mathbb{H}^{str}_{par}(n, \xi) =$ dim $U_{par}(n,\xi) + $ dim ${\bf A'}_{p, \xi}$, the morphism $\rho^{\xi}$ is generically finite. The forgetful map preserves determinants of vector bundles. As in Theorem \ref{verystableexist}, one sees that there is an open subset $S^{\xi} \subset U_{par}(n,\xi)$ such that for $(\rho^{\xi})^{-1}(E_*, 0) = (E_*, 0)$ for $(E_*,0) \in S^{\xi} \times \{0\}$.  It follows that $Fg^{\xi}_a:  
(h^{str}_{p, \xi})^{-1} (a) \dasharrow  U_{par}(n,\xi)$ is dominant for a general $a \in {\bf A'}_{p, \xi}$.    
\end{proof}

       We note that one can similarly see that 
$$Fg^{' \xi}_a: (h^{' str}_{p, \xi})^{-1} (a) \dasharrow U'_{par}(n, \xi)$$ 
 is a dominant rational map. 

\section{Semistability of the Poincar\'e sheaf} \hfill
  
       In this final section, we give an interesting application of the results of the previous sections.  We recall the following definitions from \cite[Section 2]{ABhS}.

\subsection{Degree of a sheaf on a big open subset}\hfill

        Let $W$ be a projective variety of dimension $m$  with an ample line bundle $H$. 
Let $\mathcal{E}$ denote a torsion free sheaf on the projective variety $W$. 
Let $C$ denote a general complete intersection curve rationally equivalent to  $H^{m-1}$. 
 If the singular set of $W$ has codimension at least $2$, then the general complete intersection curve 
 can be chosen to lie on the set of nonsingular points of $W$.

We define the degree of $\mathcal{E}$ with respect to the polarisation $H$ as,
\begin{equation*}
deg~ \mathcal{E}=deg ~\mathcal{E}\vert_C.
\end{equation*}

\begin{remark}
Let $U$ be an open subset of $W$ such that codimension of $S=W - U$ in $W$ is at least 2.
Since $\dim ~S \leq \dim W-2$, for any torsion-free coherent sheaf $\mathcal{F}$ on the open subset 
$U$, the direct image $i_*(\mathcal{F})$ is a torsion-free coherent sheaf on $W$.
\end{remark}

\begin{definition}
Let $U$ be an open subset of $W$ such that codimension of $W-U$ in $W$ is at least 2. Let $\mathcal{F}$ be a torsion free sheaf on $U$. We define $deg(\mathcal{F})= deg (i_*(\mathcal{F}))$, where $i:U\to W$ is the inclusion map and $i_*(\mathcal{F})$ is the direct image sheaf of $\mathcal{F}$ on $W$.

\end{definition}

        We remark that with the above definitions, the proof of \cite[Lemma 2.1]{BBrN} works in the following generality.

\begin{lemma} \label{L2.1BBrN}
   Let $X$ and $Y$ be smooth quasi projective varieties of same dimension $m$. Suppose that $X$ and $Y$ are open subsets of projective varieties $\bar{X}$ and $\bar{Y}$ respectively such that codim $(\bar{X} \setminus X, \bar{X}) \ge 2$, codim $(\bar{Y} \setminus Y, \bar{Y}) \ge 2$. Let $f : X \dasharrow Y$ be a dominant rational map defined outside a closed subset $Z \subset X$ with codim $(Z, X) \ge 2$. Let $D_X \subset X, D_Y\subset Y$ be restrictions of ample divisors on $\bar{X}$ and $\bar{Y}$ respectively. Assume that  $f^*D_Y \vert_{(X\setminus Z)} = D_X \vert_{(X \setminus Z)}$. 
   
        Let $E$ be a vector bundle on $Y$ such that $f^*E$ extends to a vector bundle $F$ on $X$. If $F$ is $D_X$-semistable (resp. stable) on $X$ then $E$ is $D_Y$-semistable (resp. stable) on $Y$.   
 \end{lemma}

\subsection{The Poincar\'e sheaf $\mathcal{U}_*$.}    \hfill

   Henceforth in this section $Y$ denotes an integral projective curve with at most nodes as singularities and $D$ is a finite set of distinct points of $Y$ away from nodes.  We assume that the greatest common divisor of $n, d$ and $n_i(x), x \in D$ is $1$. Here $n_i(x), x \in D$ are the multiplicities in the parabolic structure. Then there exists a universal family of stable parabolic torsion-free sheaves on $Y$, namely  the Poincar\'e sheaf $\mathcal{U}_*$ on  $U_{par}(n, \xi) \times Y$.  We note that $U_{par}(n, \xi)$ is a projective variety with an ample line bundle $\Theta_{par}$. The singular set $U_{par}(n,\xi) \setminus U'_{par}(n,\xi)$ is of codimension at least two in $U_{par}(n, \xi)$ \cite[Theorem 3.2]{BhBi2}.  
   
          Let $z$ be a closed point of  $Y$ distinct from the nodes and the points in $D$. Let 
$\mathcal{U}_z$ be the  restriction of the universal sheaf to $U_{par}(n, \xi) \times z$.             
Note that $\mathcal{U}_z$ has no parabolic structure. 

\begin{theorem}  \label{Uzsemistable}
$\mathcal{U}_z$ is semistable with respect to the ample line bundle $\Theta_{par}$ on $U_{par}(n, \xi)$.
\end{theorem}                        

\begin{proof}
       Our proof of this theorem is on the lines of the proofs of \cite[Theorem 2.5]{BBrN} and \cite[Theorem 3.4]{BaBiD}. 
        
      Let $a \in {\bf A}'_{p, \xi}$ be a general point (where ${\bf A}'_{p, \xi}$ is defined in the subsection \ref{Apxi}) and $\pi_p: Y^p_a \to Y$ the associated parabolic spectral curve.  
      
     Let $\mathcal{N}$ denote the restriction of the Poincare sheaf on $\bar{J}(Y^p_a) \times Y^p_a$ to $P_{\delta} \times Y^p_a$, where the Prym variety $P_{\delta}$ is defined in Subsection \ref{Pdelta}.  Let $Id \times \pi_p : P_{\delta} \times Y^p_a  \to P_{\delta} \times Y$ be the product map. Then the direct image $(Id \times \pi_p) _* \mathcal{N}$ is a Hitchin sheaf of rank $n$ with a Hitchin field 
 $$\Phi \in H^0(P_{\delta} \times Y, End((Id \times \pi_p)_* \mathcal{N}) \otimes p_Y^* \omega_Y(D))\, .
 $$
 and a parabolic structure over $P_{\delta} \times D$ induced by the $(Id \times \pi_p)_* \mathcal{O}_{P_{\delta} \times Y^p_a}$-module structure on it (by Parabolic BNR correspondence in Theorem  \ref{parBNRcorrA'}). 
 This family induces a morphism $P_{\delta} \to \mathcal{H}^{str}_{par}(n,\xi)$ , composing it with the forgetful map $\mathcal{H}^{str}_{par}(n,\xi) \dasharrow U_{par}(n,\xi)$ gives the (forgetful) rational map 
 $$Fg^{\xi}_a: P_{\delta} \dasharrow U_{par}(n,\xi)$$
  defined by mapping $N$ to the parabolic torsion-free sheaf $(\pi_p)_*(N)$. Let $T_{\delta} \subset P_{\delta}$ be the subset where $Fg^{\xi}_a$ is well defined. Then $Fg^{\xi}_a$ is dominant by Corollary \ref{forgetfuldominant}.
  Define 
  $$\mathcal{E}:= ((Id \times \pi_p)_* \mathcal{N})\vert_{T_{\delta} \times Y}\, .$$
  By the universal property of the fine moduli space $U_{par}(n,\xi)$ we have the isomorphism (on $T_{\delta}\times Y$)
 $$(Fg^{\xi}_a \times Id_Y)^* \mathcal{U} \cong \mathcal{E} \otimes L_0\, ,$$
 for some line bundle $L_0$ on $T_{\delta}$. Restricting to $U_{par}(n,\xi) \times z$, we get 
 $(Fg^{\xi}_a)^* \ \mathcal{U}_z  \cong \mathcal{E}\vert_{(T_{\delta} \times \pi_p^{-1}(z))} \otimes L_0$ so that on $T_{\delta}$ we have 
 $$(Fg^{\xi}_a)^* \ \mathcal{U}_z \ \cong \ (\oplus_{i=1}^{n} \ \mathcal{N}_{z_i}) \otimes L_0\, ,$$
 where $(\pi_p)^{-1}z = \{z_1, \cdots, z_n\}$ is a set of distinct points of $Y^p_a$. Since $z \in Y$ is a nonsingular point,  $\mathcal{N}_{z_i}$ are line bundles. In fact they are the restrictions of line bundles on $\bar{J}(Y^p_a)$ for all $i$.

      We first show that $(Fg^{\xi}_a)^* \ \mathcal{U}_z$ extends to a vector bundle on $P'_{\delta} \cup T_{\delta}$ which is semistable with respect to the restriction of any ample line bundle on $P_{\delta}$. The moduli space $U_{par}(n,\xi)$ is a GIT quotient of $R_{par, \xi}$ by PGL$(N)$ for some large $N$. The subset of $R_{par, \xi}$ corresponding to parabolic torsion-free sheaves with underlying sheaf stable is an open subset with complement of codimsion at least $2$ in $R_{par, \xi}$ \cite[Theorem 3.2]{BhBi2}.  Using this and the local universal property of $R_{par, \xi}$, one sees that codim $(P_{\delta} \setminus T_{\delta}, P_{\delta}) \ge 2$. 
 Since $P'_{\delta}$ is a smooth open subset of $P_{\delta}$,  codim $(P'_{\delta} \setminus P'_{\delta} \cap T_{\delta}, P'_{\delta}) \ge 2$, the line bundle $L_0$ uniquely extends to a line bundle $L'_0$ on $P'_{\delta}$ (and hence to $P'_{\delta} \cup T_{\delta}$). Since the line bundles $\mathcal{N}_{z_i}$ are all algebraically equivalent, the vector bundle $\oplus _i  \mathcal{N}_{z_i} \otimes L'_0$ is semistable with respect to the restriction of any ample line bundle on $P_{\delta}$.  

      We now show that $(Fg^{\xi}_a)^* \Theta_{par}$ extends to a restriction of an ample line bundle to $P'_{\delta}$. For the family $\mathcal{E}$ of parabolic torsion-free sheaves on $Y$ parametrised by $T_{\delta}$, we define 
$$\Theta_{T_{\delta}} := (Det R_{(\pi_{T_{\delta}})} \mathcal{E} )^K \otimes M'$$
where 
$$M' = \bigotimes_{x \in D} [(det  \mathcal{E}_{T_{\delta} \times x})^{b_x} \otimes \bigotimes_{i = 1}^{l_x}  (det  (\mathcal{E}_{T_{\delta} \times x}/F'_i(\mathcal{E}_{T_{\delta}  \times x}))) ] \otimes det (\mathcal{E}_{T_{\delta} \times t})^{\ell} \, .$$  
Then by the universal property of $\Theta_{par}$, on $T_{\delta}$ we have 
$$ \Theta_{T_{\delta}} \ = \ (Fg^{\xi}_a)^*(\Theta_{par})\, .$$

           Let $\tilde{P}_{\delta}$ be the pull back of $P_{\delta}$ to $\tilde{J}(Y^p_a)$, the desingularisation (normalisation) of  $\bar{J}(Y^p_a)$. Then $\tilde{P}_{\delta}$ is smooth and it has an open subset mapping isomorphically onto $P'_{\delta}$.  We identify this open subset with $P'_{\delta}$. Let $\theta_{P'_{\delta}}$ (respectively $\theta_{P_{\delta}}$) be the restriction of the theta line bundle $\theta_{\bar{J}}$, an ample line bundle on $\bar{J}(Y^p_a)$, to $P'_{\delta}$ (respectively to $P_{\delta}$). 
           The family $\mathcal{N}$ is the restriction of the family $\mathcal{N}$ on $\bar{J}(Y^p_a) \times Y^p_a$. Hence the family $\mathcal{E}$ is the restriction of the family $\mathcal{E}$ on $\bar{J}(Y^p_a) \times Y^p_a$. We denote by $\tilde{\mathcal{E}}$ its pullback to $\tilde{J}(Y^p_a)$. 

           By \cite[Theorem 1.2]{Bh5}, we have 
$$(Det R_{\pi_{\tilde{J}(Y^p_a)}} \tilde{\mathcal{E}} ) = (n^2/\bar{a}) \theta_{\tilde{J}(Y^p_a)}\vert_{\tilde{J}(Y^p_a)} \otimes M'_1\, ,$$ 
for some line bundle $M'_1$ on $\tilde{J}(Y^p_a)$ and $\bar{a}$ is the greatest common divisor of $n$ and $ d$.          
                 
Hence on $T'_{\delta} = P'_{\delta} \cap T_{\delta}$ we have  
$$ (Fg^{\xi}_a)^*(\Theta_{par}) =  \Theta_{T'_{\delta}} =    
\theta^{m_0 K}_{P'_{\delta}} \otimes M' \
{\rm for \  some} \  m_0 \in \mathbb{Z}_{> 0} \ {\rm and \ for \ some \ line \ bundle} \ M'$$  
on $T'_{\delta}$.          
We showed that codim $(P'_{\delta} \setminus T'_{\delta}, P'_{\delta}) \ge 2$. By Lemma \ref{pdeltacodim2}, codim $(P_{\delta} \setminus P'_{\delta}, P_{\delta}) \ge 2$. Hence  codim $(\tilde{P}_{\delta} \setminus T'_{\delta}, \tilde{P}_{\delta}) \ge 2$, the (restriction of the) line bundle $M'$ extends uniquely to a line bundle $M$ on $\tilde{P}_{\delta}$.    Hence 
$(Fg^{\xi}_a)^*(\Theta_{par})$ extends to a line bundle $H, H :=  \theta^{m_0 K}_{\tilde{P}_\delta} \otimes M$
on $\tilde{P}_{\delta}$.
Since the extension to $P'_{\delta}$ of $(Fg^{\xi}_a)^* \ \mathcal{U}_z$ is semistable with respect to $H$, the vector bundle $\mathcal{U}_z\vert_{U'_{par}(n,\xi)}$ is semistable with respect to the restriction of $\Theta_{par}$ by Lemma \ref{L2.1BBrN} (taking $X = P'_{\delta}, \bar{X} = \tilde{P}_{\delta}, Y = U'_{par}(n,\xi), \bar{Y} = U_{par}(n, \xi)$ in Lemma \ref{L2.1BBrN}). It follows that $\mathcal{U}_z$ is $\Theta_{par}$-semistable.
\end{proof}            

\begin{theorem}   \label{Uparstable} 
     The parabolic torsion-free sheaf $\mathcal{U}_*$ on $U_{par}(n, \xi) \times Y$ is parabolic stable with respect to any integral ample divisor of the form $a \Theta_{par} + b D_Y$ where $a, b >0$ and $D_Y$ is an ample divisor on $Y$. 
\end{theorem}          
\begin{proof}
     This follows from Theorem \ref{Uzsemistable} and \cite[Lemma 3.1]{BaBiD}. 
\end{proof}

\end{document}